\documentclass[a4paper, 10pt]{article}
\usepackage[DIV=14,BCOR=2mm,headinclude=true,footinclude=false]{typearea}
\usepackage[utf8]{inputenc}
\usepackage{lmodern}

\usepackage[left=3.5cm,right=3.5cm,top=3.5cm,bottom=4cm]{geometry}

\RequirePackage{amssymb}
\RequirePackage{amsmath}
\RequirePackage{tikz}
\RequirePackage{scalerel}
\usepackage{xfrac}
\usetikzlibrary{calc}
\newcommand\ivee{\mathrel{\rotatebox[origin=c]{270}{$\geqslant$}}}
\newcommand\dna{\mathtt{DNA}}

\newcommand\cpc{\mathtt{CPC}}
\newcommand\ipc{\mathtt{IPC}}

\newcommand\at{\mathtt{AT}}
\newcommand\MP{\mathtt{MP}}
\newcommand\US{\mathtt{US}}

\newcommand\langInt{\mathcal{L}_{\ipc}}
\newcommand\langIntsta{\mathcal{L}_{\mathtt{Cl}}}
\newcommand\langInqI{\mathcal{L}_{\ipc}^\otimes}
\newcommand\langInqst{\mathcal{L}^\otimes_{\mathtt{Cl}}}

\newcommand\inq{\mathtt{Inq}}
\newcommand\inqint{\mathtt{Inq\Lambda}}
\newcommand\inqI{\mathtt{InqI}}
\newcommand\inqB{\mathtt{InqB}}

\tikzstyle{dot}=[inner sep=1pt, fill, black, circle, draw, minimum size = 5pt]
\tikzstyle{highlight}=[inner sep=1pt, draw=black, fill=white, circle, minimum size = 9pt]

\newcommand\dep{\mathop{=\!}\xspace}

\tikzstyle{world}=[inner sep=2pt, black, circle, fill=lightgray]
\tikzstyle{worldBlack}=[inner sep=2pt, white, circle, fill=black]
\tikzstyle{worldWhite}=[inner sep=2pt, black, circle, fill=white, draw]

\tikzstyle{modelbox}=[rounded corners]
\tikzstyle{info}=[fill=gray,fill opacity=.1, rounded corners]
\tikzstyle{counter}=[info, densely dashed]

\makeatletter
\newcommand\xLongLeftRightArrow[2][]{%
	\ext@arrow 0099{\LongLeftRightArrowfill@}{#1}{#2}}
\def\LongLeftRightArrowfill@{%
	\arrowfill@\Leftarrow\Relbar\Rightarrow}
\makeatother

\usepackage{enumitem}
\usepackage{csquotes}
\usepackage{syllogism}
\usepackage{bussproofs} 
\usepackage{lscape}
\usepackage{amsmath}
\usepackage{footnote}
\usepackage{amsthm}
\usepackage[hang]{footmisc}
\usepackage{setspace}
\usepackage{tikz}
\usetikzlibrary{cd}
\usepackage{ stmaryrd }
\usepackage{bussproofs}

\newtheorem{theorem}{Theorem}[section]

\newtheorem{proposition}[theorem]{Proposition}
\newtheorem{corollary}[theorem]{Corollary}
\newtheorem{lemma}[theorem]{Lemma}
\newtheorem*{claim}{Claim}

\theoremstyle{definition}

\theoremstyle{definition}
\newtheorem{definition}[theorem]{Definition} % definition numbers are dependent on theorem numbers
\allowdisplaybreaks

\renewcommand\qedsymbol{$\square$}

\newtheoremstyle{TheoremNum}
{\topsep}{\topsep}              %%% space between body and thm
{\itshape}                      %%% Thm body font
{}                              %%% Indent amount (empty = no indent)
{\bfseries}                     %%% Thm head font
{.}                             %%% Punctuation after thm head
{ }                             %%% Space after thm head
{\thmname{#1}\thmnote{ \bfseries #3}}%%% Thm head spec
\theoremstyle{TheoremNum}

\usepackage[backend=biber, isbn=false, doi=false, url=false, style=numeric, citestyle=numeric, giveninits=true]{biblatex}
\bibliography{bibliografia.bib}
\title{On Intermediate Inquisitive and Dependence Logics:\\ An Algebraic Study \thanks{I am very grateful to Fan Yang for many helpful conversations and for her constant support and advice. I would also like to thank Georgi Nakov, Gianluca Grilletti, Tommaso Moraschini and Nick Bezhanishvili for helpful discussions. This research was supported by grant 336283 of the Academy of Finland and Research Funds of the University of Helsinki. }}

\author{Davide Emilio Quadrellaro}
\date{ }
\begin{document}
	\maketitle
	\noindent
	%ADD REFERENCE TO ALLEN MANN PAPER
	
	\begin{abstract}
		This article provides an algebraic study of intermediate inquisitive and dependence logics. While these logics are usually investigated using team semantics, here we introduce an alternative algebraic semantics and we prove it is complete for all intermediate inquisitive and dependence logics. To this end, we define inquisitive and dependence algebras and we investigate their model-theoretic properties.  We then focus on finite, core-generated, well-connected inquisitive and dependence algebras: we show they witness the validity of formulas true in inquisitive algebras, and of formulas true in well-connected dependence algebras. Finally, we obtain representation theorems for finite, core-generated, well-connected, inquisitive and dependence algebras and we prove some results connecting team and algebraic semantics.
	\end{abstract}

\section*{Introduction}
In this work, we pursue an algebraic study of inquisitive and dependence logic. Although the connection between these two logical systems is now firmly established, inquisitive and dependence logic were introduced in different contexts and with different research interests in mind. 

Inquisitive logic was formally developed by Ciardelli, Groenendijk and Roelofsen in a series of articles, most notably in \cite{Ciardelli:09salt,Ciardelli2011-CIAIL}, where they introduced the so-called ``support semantics''. Atomic formulas are assigned, under this semantics, to sets of possible worlds in a Kripke model. Inquisitive logic was developed hand-in-hand with inquisitive semantics -- a linguistic framework that aims at providing a uniform formal characterisation of both questions and statements in natural languages. In particular, polar questions expressing ``\textit{whether p holds or not}'' are represented by an operator $?p$ defined using the inquisitive disjunction as $?p:=p\ivee \neg p$. The interested reader can refer to \cite{Ciardelli2018-CIAIS} for a general introduction to inquisitive semantics.

Dependence logic, on the other hand, was introduced by V\"a\"an\"anen \cite{Vaananen2007-VNNDLA} as an extension of first-order logic with dependence atoms. The underlying motivation of dependence logic was to provide a logical framework able to capture several  relations of \textit{dependencies} between variables. In its standard formulation, dependence logic is defined via team semantics, originally introduced in \cite{Hodges}, which generalises standard Tarski's semantics by teams, which are sets of assignments that map first-order variables to elements of the domain. In its propositional version, a team is a set of valuations mapping propositional atoms to either 1 or 0. Propositional dependence logic has been extensively studied by Yang and V\"a\"an\"anen in \cite{Yang2016-YANPLO}, while in \cite{Yang2017-YANPTL} they considered several extensions of classical logic using team semantics. Intuitively, the dependence atom $\dep(\vec{p},q)$ expresses the fact that the value of the variable $q$ is uniquely determined by the values of the variables $\vec{p}$. The constancy atom $\dep(p)$ can then be seen as a special case of the dependency atom, saying that the value of a variable is constant in the underlying team.	

It was soon noticed that, at the propositional level, the team semantics of dependence logic and the state semantics of inquisitive logic are in fact equivalent -- states of possible worlds are nothing but teams of propositional assignments. The close connection between these two approaches was pointed out and developed e.g. in \cite{Ciardelli2016,Yang2016-YANPLO}. As a matter of fact, both dependence and inquisitive logic are  expressively equivalent and they are both complete with respect to the class of all downward closed team properties.

The connection between inquisitive and dependence logic was recently pushed further by \textcite{ciardelli2020}, which introduced versions of propositional inquisitive and dependence logics which are based on intuitionistic, rather than classical logic. Here we generalise their approach -- also drawing on \cite{10.1093/logcom/exw021} -- and we provide an algebraic study of intermediate inquisitive logics $\inqint$, which extend the intuitionistic inquisitive logic $\inqI$, and  intermediate dependence logics $\inqint^\otimes$, which extend the intuitionistic dependence logic $\inqI^\otimes$.

The interest for algebraic semantics of inquisitive and dependence logic is not new and some works in the literature already consider the issue. The study in the algebraic interpretations of dependence logic was initiated by \textcite{abr} and later developed by \textcite{Luck}. An early work on inquisitive logic from an algebraic perspective is \cite{Roelofsen2013-ROEAFF}. More recently,  \textcite{grilletti} have introduced an algebraic and topological semantics for the system $\inqB$ of classical propositional inquisitive logic, while Bezhanishvili, Grilletti and Quadrellaro \cite{Quadrellaro.2019B,grillq,Quadrellaro.2019} have further developed this approach and extended it to other logics. In a similar fashion,  \textcite{puncochar_inquisitive_2021} has introduced a semantics for several extensions of intuitionistic inquisitive logic.

In the general framework of abstract algebraic logic, as described e.g. in \cite{Font.2016}, logics are described as consequence relations which are additionally closed under uniform substitution. However, both inquisitive and dependence logic are not closed under uniform substitution and they are thus not logics in this strict sense of the word. As a consequence of this fact, we cannot directly apply the standard framework of abstract algebraic logic. In fact, the algebraic semantics for versions of inquisitive logics described in \cite{grilletti} and \cite{puncochar_inquisitive_2021} are quite non-standard: inquisitive logics are shown to be complete with respect to some classes of Heyting algebras which are not axiomatisable by means of (quasi-)equations, i.e. which do not form algebraic (quasi-)varieties. 	

In the present work, we aim at providing an algebraic study of inquisitive and dependence logics which is closer in spirit to the usual approach of abstract algebraic logic. To this end, we introduce both inquisitive and dependence algebras in terms of classical model theory and we show that they form an elementary class axiomatised by universal Horn sentences. Therefore, even if they do not form a variety, inquisitive and dependence algebras make for an interesting class which is suitable of further model-theoretic investigations. Interestingly, the classes of algebras investigated in \cite{grilletti} and in \cite{puncochar_inquisitive_2021} can be then seen as special collections of representatives -- respectively for the classical and the intuitionistic case.

The main results of this article are the following. Firstly, we prove a full completeness result for intermediate inquisitive and dependence logics, which states that intermediate inquisitive logics are complete with respect to their corresponding class of intermediate inquisitive algebras and that intermediate dependence logics are complete with respect to their corresponding class of intermediate dependence algebras. We then investigate several model-theoretic properties of inquisitive and dependence algebras. In particular, we show several results concerning finite, core-generated, well-connected inquisitive and dependence algebras and we prove two representation theorems for these subclasses of algebras. Finally, we extend these results to algebraic models and we obtain some bridge principles between teams and algebraic semantics.

%Also, we qualify one sense in which intermediate inquisitive and dependence logics can be regarded as \textit{algebraisable} logics and we investigate several model-theoretic properties of inquisitive and dependence algebras. In particular, we show several results concerning finite, core-generated, well-connected inquisitive and dependence algebras and we prove two Representation Theorems for these subclasses of algebras. Finally, we extend these results to models and we obtain some bridge principles concerning teams and algebraic models.

The structure of the present article is the following. In Section \ref{one} we introduce the syntax and the usual team semantics of inquisitive and dependence logic. In Section \ref{two},  we define inquisitive and dependence algebras, we introduce so-called core semantics and we prove using the method of free algebras that this semantics is complete for all intermediate inquisitive and dependence  logics.  In Section \ref{three}, we focus on the model-theoretic properties of inquisitive and dependence algebras and prove several results concerning finite, core-generated, well-connected inquisitive and dependence algebras. Later, in Section \ref{four}, we provide representation theorems for finite, core-generated, well-connected inquisitive and dependence algebras and we prove some results concerning the team and the algebraic semantics of inquisitive and dependence logic.  Finally, in Section \ref{ten}, we review our results and highlight some future line of research.

\section{Inquisitive and Dependence Logic}\label{one}

In this section we introduce inquisitive and dependence logics in axiomatic terms and we recall their standard team semantics.

\subsection{Axiomatic Systems}

We fix at the outset some propositional signatures for different systems of logic. Throughout this paper we shall always denote by $\at$ a fixed set of atomic propositional variables, which we will always assume to be countable. We denote by $\langInt$ the standard signature of intuitionistic logic $\langInt=\{ \land, \lor, \rightarrow, \bot \}$. With slight abuse of notation, we denote by $\langInt$ also the set of formulas built recursively from $\at$ in this signature: We let $\phi\in \langInt$ if and only if $\phi$ is generated by the following grammar: 
\begin{align*}
\phi ::= p \mid  \bot  \mid  \phi \land \phi \mid \phi \lor \phi \mid \phi\rightarrow\phi;
\end{align*}
\noindent where $p\in\at$ is any propositional variable. Negation is treated as a defined operation and can be introduced by letting $\neg \phi := \phi \to \bot$.

The \textit{intuitionistic propositional calculus} $\ipc$ -- \textit{intuitionistic logic} for short -- is the set of formulas of $\langInt$ which contains the usual axioms, it is closed under \textit{modus ponens} $ (\MP) $ and uniform substitution $ (\US) $. An \textit{intermediate logic} $\mathtt{L}$ is a consistent set of formulas of $\langInt$ which contains $\ipc$ and is closed under \textit{modus ponens} and uniform substitution. Intermediate logics are known to form a lattice structure, whose maximal element is the \textit{classical propositional calculus} $\cpc$. We refer the reader to \cite{Zakharyaschev.1997} for more on intermediate logics.

In this paper we shall formulate inquisitive logic in the language $\langInt$ -- hereby following Ciardelli's original presentation in \cite{Ciardelli.2009} -- and we will adopt the standard disjunction symbol $\vee$ in place of the more common $\ivee$ symbol to denote the inquisitive disjunction operation. The reason of this choice is that we want to stress that inquisitive disjunction is nothing but intuitionistic disjunction. In fact, as it is often remarked in the literature, inquisitive logic is very close to intermediate logics, as it contains $\ipc$ and it is contained in $\cpc$. However, the fact that inquisitive logic is not closed under uniform substitution means that it is not an intermediate logic. In fact, for this very same reason, inquisitive logic does not fit the framework of \textit{abstract algebraic logic} (see e.g. \cite{Font.2016}), where logics are defined as consequence operators closed under uniform substitutions.

Although inquisitive logic does not admit full substitution, it is closed under a restricted version of substitution, namely substitution of $\vee$-free formulas. This reflects the fact that in inquisitive semantics $\vee$-free formulas correspond to \textit{sentences}, while formulas containing $\vee$ are intended to model \textit{questions}. A formula of $\langInt$ is said to be \textit{standard} if it is $\vee$-free. We write  $\langIntsta$ for the set of all standard formulas and also for the signature $\langIntsta=\{ \land, \to, \bot \}$.

Inquisitive logic is usually presented in semantical terms, as the logic of states of possible worlds. However, it is also possible to define it in more syntactical terms. Here we adapt the natural deduction system presented in \cite{ciardelli2020} and we present it in a Hilbert-style fashion. We use Greek letters $\phi,\psi,\dots$ as meta-variables for arbitrary inquisitive formulas and $\alpha,\beta,\dots$ as meta-variables for arbitrary standard inquisitive formulas. We then define intuitionistic inquisitive logic in the following way.

\begin{definition}[Intuitionistic Inquisitive Logic] \label{inqI}
	
	The system $\inqI$ of \textit{intuitionistic inquisitive logic} is the smallest set of formulas of $\langInt$ such that, for all $\phi,\psi,\chi \in \langInt$ and for all $\alpha\in \langIntsta$, $\inqI$ contains the following formulas:
	\begin{align*}
	\text{(A1) }	& \phi\rightarrow(\psi\rightarrow \phi)\\		 
	\text{(A2) }	& (\phi\rightarrow (\psi\rightarrow \chi))\rightarrow (\phi\rightarrow \psi) \rightarrow (\phi\rightarrow \chi) \\
	\text{(A3) }	&\phi \land \psi \rightarrow \phi\\
	\text{(A4) }	& \phi \land \psi \rightarrow \psi\\
	\text{(A5) }	& \phi\rightarrow (\psi\rightarrow \phi \land \psi)\\
	\text{(A6) }	& \phi \rightarrow \phi\lor \psi\\
	\text{(A7) }	& \psi \rightarrow \phi\lor \psi\\
	\text{(A8) }	& (\phi\rightarrow \chi) \rightarrow ((\psi\rightarrow \chi)\rightarrow (\phi \lor \psi \rightarrow \chi))\\
	\text{(A9) }	&\bot \rightarrow \phi\\
	\text{(A10) }	& (\alpha \rightarrow (\phi \lor \psi)) \rightarrow ((\alpha \rightarrow \phi)\lor (\alpha \rightarrow \psi))		
	\end{align*} 
	\noindent and in addition it is closed under the rule of \textit{modus ponens} $ (\MP) $.
	
\end{definition}

\noindent Since the axioms schemas (A1)--(A9) plus \textit{modus ponens} axiomatise intuitionistic logic, we can think of inquisitive logic as a \textit{theory} extending intuitionistic logic:
$$ \inqI:= \mathtt{MP}(\ipc + (\text{A10}));   $$
\noindent meaning that inquisitive logic is a theory over intuitionistic logic which contains every admissible instance of the schema (A10) and is closed under \textit{modus ponens}. We will often refer to (A10) as the \textit{Split} axiom.

Similarly to the case of $\ipc$, we can define several extensions of intuitionistic inquisitive logic. We say that a set $\Lambda\subseteq \langInt$ is closed under \textit{standard} substitution if it is closed under every substitution assigning standard formulas to atomic formulas. 

\begin{definition}
	An \textit{intermediate inquisitive logic} is any set of formulas $\inqint$ such that: $\inqint =  \mathtt{MP}(\inqI \cup \Lambda ) $, where $\Lambda\subseteq \mathcal{L}_{\mathtt{CL}}$ is any set of standard formulas closed under standard substitution.
\end{definition}

\noindent We notice that this definition differs from the so-called \textit{inquisitive superintuitionistic logic}$^*$ defined by \textcite{puncochar_inquisitive_2021}, for we do not require intermediate inquisitive logics to satisfy the disjunction property and we do not allow non-standard formulas in $\Lambda$. It is easily verified that, if $\inqint$ and $\mathtt{inq\Delta}$ are intermediate inquisitive logics, then also $\inqint \cap \mathtt{inq\Delta}$ and $\MP(\inqint \cup \mathtt{inq\Delta}   )$ are intermediate inquisitive logics. Hence, intermediate inquisitive logics form a lattice whose least element is $\inqI$. An important example of an intermediate inquisitive logic is the classical version of inquisitive logic $\inqB$. This is defined as the extension of $\inqI$ by the axiom scheme $\neg \neg \alpha \to \alpha$, where $\alpha\in \langIntsta$.

\begin{definition}
	The system $\inqB$ of classical inquisitive logic is defined as:  
	$$\inqB:= \mathtt{MP}(\inqI \cup \{\neg \neg \alpha \rightarrow \alpha\}_{\alpha \in \langIntsta}    ).$$ 
\end{definition}

\noindent It is then easy  to see that $\inqB$ is the maximal element in the lattice of intermediate inquisitive logics.

Dependence logics extend inquisitive logics in an expanded syntax. Here we shall take a slightly non-standard approach and formulate dependence logic in the vocabulary $\langInqI$, which expands $\langInt$  by adding the \textit{tensor disjunction} operator $\otimes$. Intuitively, the tensor disjunction is meant to be as much as a \enquote{classical} disjunction as possible in the given intuitionistic framework. We fix  the signature of dependence logic $\langInqI=\{ \land, \lor, \rightarrow, \otimes, \bot \}$ and with slight abuse of notation we denote by $\langInqI$ also the set of formulas  defined by induction in this signature over $\at$, i.e.	$\phi\in \langInqI$ if and only if $\phi$ is generated by the following grammar: 
\begin{align*}
\phi &::= p\mid \bot  \mid  \phi \land \phi \mid \phi \lor \phi \mid \phi\rightarrow\phi\mid \phi\otimes\phi.	
\end{align*}
\noindent	where $p\in\at$ is an arbitrary atomic variable. As it is the case for inquisitive logics, we are often interested in formulas which are $\vee$-free -- i.e. which do not contain the $\vee$ symbol -- and we refer to such formulas as \textit{standard} formulas. We write $\langInqst$ for the set of standard dependence formulas and also for the restricted signature $\langInqst=\{ \land, \rightarrow, \otimes,  \bot\}$. Clearly $\langInt\subseteq \langInqI$ and $\langIntsta\subseteq\langInqst $.

Negation is defined over $\langInqI$ as in intuitionistic logic, by letting $\neg\phi:= \phi\rightarrow\bot$.	More interestingly, we can define the so-called \textit{constancy} and \textit{dependency atoms}, as partial operations defined only on atomic formulas. Intuitively, the former says that the value of an atomic formula is constant, while the latter says that the value of an atomic formula is functionally determined by the value of a tuple of other atomic formulas. Let $\vec{p}=p_0,\dots,p_n$, we define them as follows:
\begin{align*}
\dep(p) &:= p \lor \neg p ;\\
\dep(\Vec{p},q) &:= (\bigwedge_{i\leq n} \dep(p_i) )\rightarrow \dep(q).
\end{align*} 

\noindent Notice that, since $\otimes$ does not occur in the definitions above, these operators can be introduced both in $\langInt$ and $\langInqI$. What is specific of dependence logic is the presence of the tensor disjunction rather than the dependence operator itself. In fact, it was shown By Barbero and Ciardelli in \cite{10.1007/978-3-662-60292-8_3} that although inquisitive and dependence logics are both expressively complete with respect to the class of all downward closed teams, the tensor disjunction is not definable in inquisitive logic.

We define intuitionistic dependence logic analogously to how we defined intuitionistic inquisitive logic. 

\begin{definition}[Intuitionistic Dependence Logic]\label{inqIotimes}
	The system $\inqI^\otimes$ of \textit{intuitionistic dependence logic} is the smallest set of formulas of $\langInqI$ such that, for all $\phi,\psi,\chi,\tau \in \langInqI$ and for all $\alpha,\beta,\gamma\in \langInqst$, $\inqI^\otimes$ contains the formulas (A1)--(A10) of Definition \ref{inqI} and the following: 
	\begin{align*}
	\text{(A11) }	& \alpha\rightarrow(\alpha \otimes \beta)\\	 
	\text{(A12) }	& (\alpha \otimes \beta) \rightarrow (\beta \otimes \alpha)\\
	\text{(A13) }	&  \phi \otimes (\psi \lor \chi) \rightarrow (\phi\otimes \psi) \lor (\phi\otimes \chi) \\
	\text{(A14) }	&(\phi\rightarrow \chi) \rightarrow ((\psi\rightarrow \tau)\rightarrow (\phi \otimes \psi \rightarrow \chi \otimes \tau))\\
	\text{(A15) }	&(\alpha\rightarrow \gamma) \rightarrow ((\beta\rightarrow \gamma)\rightarrow (\alpha \otimes \beta \rightarrow \gamma)).
	\end{align*} 
	\noindent And in addition it is closed under \textit{modus ponens}.
\end{definition}

\noindent We often refer to (A13) as the \textit{Dist} axioms and to (A14) as the \textit{Mon} axiom. We can think of $\inqI^\otimes$ in the following way:
$$ \inqI^\otimes:= \MP(\ipc + (A10 - A15)).   $$
\noindent Hence, dependence logic is a theory extending $\ipc$ in the language $\langInqI$. 

Intermediate dependence logics are defined as follows. We say that a set $\Lambda\subseteq \langInqI$ is closed under \textit{standard} substitution if it is closed under every substitution assigning standard formulas to atomic formulas. 

\begin{definition}		
	An \textit{intermediate dependence logic} is any set of formulas $\inqint^\otimes$ such that: $\inqint =  \MP(\inqI^\otimes \cup \Lambda ) $, where $\Lambda\subseteq \langInqst$ is any set of standard formulas closed under standard substitution.
\end{definition}

\noindent If we consider the subset of standard formulas in an intermediate dependence logic $\inqint^\otimes$, what we obtain is an intermediate logic with tensor as the disjunction operator. Let  $\ipc^\otimes$ refer to the intuitionistic propositional calculus in the signature $\langInqst$, i.e. in a syntax where the usual disjunction is replaced by the tensor, then for any set $\Lambda\subseteq \langInqst$  we have that  $\ipc^\otimes \subseteq \inqint^\otimes {\upharpoonright} \langInqst $ and that $ \inqint^\otimes {\upharpoonright} \langInqst $ is closed under \textit{modus ponens} and uniform substitution, hence $\inqint^\otimes {\upharpoonright} \langInqst$ is an intermediate logic.  For any intermediate dependence logic $\inqint^\otimes$ we denote by $\ipc^\otimes +\Lambda$ the intermediate logic $\inqint^\otimes {\upharpoonright} \langInqst$.

An important example of intermediate dependence logics is the classical version of dependence logic $\inqB^\otimes$. This is defined as follows.

\begin{definition}
	The system $\inqB^\otimes$ of classical dependence logic is defined as:  
	$$\inqB^\otimes:= \mathtt{MP}(\inqI^\otimes + \{\neg \neg \alpha \rightarrow \alpha\}_{\alpha \in \langInqst}    ).$$ 
\end{definition}
\noindent If we denote by $\cpc^\otimes$ the classical propositional calculus in the signature $\langInqst$, i.e. in a syntax where the usual disjunction is replaced by the tensor, then it can be seen that $\cpc^\otimes \subseteq \inqB^\otimes$, which shows the sense in which dependence logic is an extension of classical propositional logic. As in the case of intermediate inquisitive logics, intermediate dependence logics form a bounded distributive lattice, of which $\inqI^\otimes$ is the least  and $\inqB^\otimes$ the greatest element.

\subsection{Semantics via Teams}

Inquisitive and dependence logics are usually introduced via some version of team semantics. In particular, \textcite{ciardelli2020} have defined a version of team semantics based on Kripke models, while the classical version of team semantics dates back to  \textcite{Hodges} and was already used in \cite{Vaananen2007-VNNDLA} and \cite{Ciardelli.2009}. We describe here the team semantics on Kripke models and we explain how standard team semantics can be seen as a special case of it.

%We introduce teams over Kripke models to provide a semantics to every formulas in the language $\langInqI$. Clearly, since $\langInt\subseteq \langInqI$, it follows that this semantics allows us to interpret both inquisitive and dependence formulas. Notice in particular that, since $\langInt\subseteq \langInqI$, every result which holds for all dependence formulas applies also to inquisitive formulas.	

Firstly, we recall that an \textit{intuitionistic Kripke frame} is a partial order $\mathfrak{F}=(W, R)$, where $W$ is a set of possible worlds and $R$ a partial ordering, i.e. a reflexive, transitive, and antisymmetric relation. An \textit{intuitionistic Kripke model} is  a pair $\mathfrak{M}=(\mathfrak{F},V)$, where $\mathfrak{F}$ is an intuitionistic Kripke frame and $V:W\rightarrow \wp(\at)$ a valuation of atomic formulas such that, if $p\in V(w)$ and $wRv$, then $p\in V(v)$. In this article Kripke frames and Kripke models are always meant to be intuitionistic Kripke frames and intuitionistic Kripke models. A world in a model can be viewed as a label for a subset of atomic formulas -- hence we shall write $w(p)=1$ if and only if $p\in V(w)$. In this sense a world $w$ corresponds to a \textit{classical assignmen}t $w:\at \rightarrow 2$. The notions of team and extension of a team are defined as follows.

\begin{definition}
	Let $\mathfrak{M}=(W,R,V)$ be an intuitionistic Kripke model. A \textit{team} is any subset $t\subseteq W$ of the set of possible worlds. A team $s$ is an \textit{extension} of a team $t$ if $s\subseteq R[t]$.
\end{definition}

\noindent A team is a set of possible worlds, hence, by our previous considerations, a team can be  considered as a set of assignments. The \textit{team semantics} (or \textit{support semantics}) of the logics $\inqI$ and $\inqI^\otimes$ is defined as follows.

\begin{definition}[Kripke Team Semantics]
	Let $\mathfrak{M}=(W,R,V)$ be an intuitionistic Kripke model. The notion of a formula $\phi\in\langInqI$ being \textit{true in a team} $t\subseteq W$ is defined as follows: 
	\begin{equation*}
	\begin{array}{l @{\hspace{1em}\Longleftrightarrow\hspace{1em}} l}
	\mathfrak{M},t\vDash p & {}\forall w\in t \ ( w(p)=1)  \\
	\mathfrak{M},t\vDash \bot &  t=\varnothing\\
	\mathfrak{M},t\vDash \psi \lor \chi & \mathfrak{M},t\vDash \psi \text{ or } \mathfrak{M},t\vDash \chi\\
	\mathfrak{M},t\vDash \psi \land \chi & \mathfrak{M},t\vDash \psi \text{ and } \mathfrak{M},t\vDash \chi\\
	\mathfrak{M},t\vDash \psi \otimes \chi & \exists s,r\subseteq t \text{ such that } s\cup r = t \text{ and } \mathfrak{M},s\vDash \psi, \mathfrak{M},r\vDash \chi \\
	\mathfrak{M},t\vDash \psi \rightarrow \chi &  \forall s \ ( \text{if }s\subseteq R[t] \text{ and } \mathfrak{M},s\vDash\psi \text{ then } \mathfrak{M},s\vDash \chi ).\\
	\end{array}
	\end{equation*}
\end{definition}

\noindent We write $\mathfrak{M}\vDash \phi$ if $\mathfrak{M},t\vDash\phi$ for all $t\subseteq W$ and $\mathfrak{F}\vDash \phi$ if $(\mathfrak{F},V)\vDash \phi$ for all valuations $V$. If $\mathcal{C}$ is a class of Kripke frames, we write $\mathcal{C}\vDash \phi$ if, for all $\mathfrak{F}\in \mathcal{C}$, we have that $\mathfrak{F}\vDash \phi$. We write $\mathfrak{M}\vDash \Gamma$ if $\mathfrak{M}\vDash \phi$ for all $\phi\in \Gamma$ and  $\mathfrak{F}\vDash \Gamma$ if $\mathfrak{F}\vDash \phi$ for all $\phi\in \Gamma$. For $\Gamma\cup\{\phi\}\subseteq \langInqI$, we write $\Gamma \vDash \phi$ if, for all Kripke models $\mathfrak{M}$, $\mathfrak{M}\vDash \Gamma$ entails $\mathfrak{M}\vDash \phi$. We write $\phi \equiv \psi$ if $\phi \vDash \psi$ and $\psi \vDash \phi$.

Let $\mathfrak{F}=(W,R)$ and $\mathfrak{G}=(W',R')$. We recall that a function $p: \mathfrak{F}\to\mathfrak{G}$ is said to be a \textit{p-morphism} if (i) $xRy$ entails $f(x)R'f(y)$ and (ii) if $f(x)R'y$ then there is $z\in W$ such that $f(z)=y$ and $xRz$. We denote by \textsf{KF} the category of intuitionistic Kripke frames with p-morphisms.

The following theorem was essentially proved in \cite{ciardelli2020} and it shows that $\inqI$ is sound and complete with respect to this version of team semantics.

\begin{theorem}[Ciardelli, Iemhoff, Yang]\label{complteam}
	For any formula $\phi\in \langInt$ and any formula $\psi\in \langInqI$, we have that:
	\begin{align*}
	\phi\in \inqI & \Longleftrightarrow  \mathfrak{F}\vDash \phi \text{ for all  intuitionistic Kripke frames } \mathfrak{F}; \\
	\psi\in \inqI^\otimes & \Longleftrightarrow  \mathfrak{F}\vDash \psi \text{ for all  intuitionistic Kripke frames } \mathfrak{F}.	
	\end{align*}
\end{theorem}

Let us now consider the special case of classical inquisitive and classical dependence logic. If a Kripke frame $\mathfrak{F}$ is such that $\mathfrak{F}\vDash \inqB$ or $\mathfrak{F}\vDash \inqB^\otimes$ then, for any standard formula $\alpha\in \langInqst$:
$$ \mathfrak{F}\vDash \neg\neg\alpha \rightarrow\alpha .$$ 
\noindent From this it follows (by Proposition \ref{classformula} and  \cite{Zakharyaschev.1997}) that $\mathfrak{F}=(W,R)$ is a classical frame, meaning that its underlying order trivialises, i.e. it follows that $R=id$.

As a consequence of this fact, we can give a simpler description of team semantics in the case of classical inquisitive and dependence logics. Since an assignment (also \textit{valuation}) is a function $w:\at\to 2$, then $2^{\at}$ is the set of all classical assignments. A \textit{team} is then a set of assignments $t\subseteq 2^\at$ and $\wp(2^\at)$ is the set of all teams over $\at$. Classical Kripke frames can be thus simply seen as sets of possible worlds or, equivalently, as sets of classical assignments. Therefore, a classical Kripke frame is simply a team. We then define as follows the classical team semantics for $\inqB$ and $\inqB^\otimes$.

\begin{definition}[Team Semantics of $\inqB^\otimes$]
	The notion of a formula $\phi\in\langInqI$ being \textit{true in a team} $t\in \wp({2^\at})$ is defined as follows: 
	
	\begin{equation*}
	\begin{array}{l @{\hspace{1em}\Longleftrightarrow\hspace{1em}} l}
	t\vDash p & {}\forall w\in t \ ( w(p)=1)  \\
	t\vDash \bot &  t=\varnothing \\
	t\vDash \psi \lor \chi & t\vDash \psi \text{ or } t\vDash \chi\\
	t\vDash \psi \land \chi & t\vDash \psi \text{ and } t\vDash \chi\\
	t\vDash \psi \otimes \chi & \exists s,r\subseteq t \text{ such that } s\cup r = t \text{ and } s\vDash \psi, r\vDash \chi \\
	t\vDash \psi \rightarrow \chi &  \forall s \ ( \text{if }s\subseteq t \text{ and } s\vDash\psi \text{ then } s\vDash \chi ).
	\end{array}
	\end{equation*}
\end{definition}

\noindent The notions of truth and the related ones are defined as in the more general case above.

The following result was proven by  Ciardelli and Roelofsen  \cite{Ciardelli2011-CIAIL} for $\inqB$ and extended by Yang and Väänänen  \cite{Yang2016-YANPLO} to $\inqB^\otimes$.

\begin{theorem}[Ciardelli, Roelofsen, Yang, Väänänen]\label{classcomplteam}
	For any formula $\phi\in \langInt$ and any formula $\psi\in \langInqI$, we have that:
	\begin{align*}
	\phi\in \inqB &\Longleftrightarrow  \wp(2^\at) \vDash  \phi; \\
	\psi\in \inqB^\otimes & \Longleftrightarrow  \wp(2^\at) \vDash  \psi.	
	\end{align*}
\end{theorem}

\subsection{Properties of Inquisitive and Dependence Logic}

We recall some important properties of team semantics over Kripke models and their special formulation in the classical setting. We omit the proofs of these results and refer the interested reader to \cite{ciardelli2020, Ciardelli2011-CIAIL, Yang2016-YANPLO, Yang2017-YANPTL}.

\begin{proposition}[Downward Team Property]
	For every Kripke model $\mathfrak{M}=(W,R,V)$, for every team $t\subseteq W$ and for every inquisitive or dependence formula $\phi\in\langInqI$ we have that $\mathfrak{M},t\vDash \phi$ and $s\subseteq R[t]$ entail $M,s\vDash \phi$. 	For every classical team $t$ and $s\subseteq t$  we have that $t\vDash \phi$ entails $s\vDash \phi$.
\end{proposition}

\noindent The next corollary allows us to conclude that a formula $\phi$ is satisfiable if and only if it is satisfied by some upward-closed team, i.e. by some team $t$  such that $R[t]\subseteq t$.

\begin{corollary}[Up-Set property]\label{upsetprop}
	For every Kripke model $\mathfrak{M}$ and for every inquisitive or dependence formula $\phi\in\langInqI$ we have that:
	\begin{align*}
	\mathfrak{M},t\vDash \phi \Longleftrightarrow \mathfrak{M},R[t]\vDash \phi.
	\end{align*}
\end{corollary}

A second key property of inquisitive and dependence logic is the Empty Team Property, which states that every formula is true in the empty team.

\begin{proposition}[Empty Team Property] 
	For all $\phi\in\langInqI$ and for every Kripke model $\mathfrak{M}$  we have that $\mathfrak{M},\varnothing\vDash\phi$.  For all $\phi\in\langInqI$ we have  that $\varnothing\vDash\phi$.
\end{proposition}

We also recall that the logics $\inqI,\inqI^\otimes,\inqB,\inqB^\otimes$ satisfy the finite model property. It is a non-trivial problem whether this property can be extended also to other intermediate inquisitive and dependence logics. We will use this fact in Section \ref{four} to give a completeness proof for $\inqI,\inqB,\inqB^\otimes$ that does not make use of free algebras.

\begin{theorem}[Finite Model Property]\label{FMP}
	$\;$
	\begin{itemize}
		\item[(i)] For all $\phi\in \langInt$, if $\phi\notin\inqI$ then there is a finite Kripke model $\mathfrak{M}$ and a finite team $t$ such that $\mathfrak{M},t\nvDash \phi$. If $\phi\notin\inqB$, then there is a finite team $t\in \wp(2^\at)$ such that $t \nvDash  \phi$.
		
		\item[(ii)] 
		For all $\phi\in \langInqI$, if $\phi\notin\inqI^\otimes$ then there is a finite Kripke model $\mathfrak{M}$ and a finite team $t$ such that $\mathfrak{M},t\nvDash \phi$. If $\phi\notin\inqB^\otimes$, then there is a finite team $t\in \wp(2^\at)$ such that $t \nvDash  \phi$.
	\end{itemize}
\end{theorem}

Finally, we recall the two following results. The next proposition gives an important characterisations of standard formulas in team semantics \cite[Prop. 3.10]{ciardelli2020}.
\begin{proposition} \label{classformula}	
	Let $\phi \in \langInqI$, then there is some $\alpha\in\langInqst$ such that $\phi\equiv \alpha$ if and only if the following condition holds, for all Kripke model $\mathfrak{M}$:
	$$ \mathfrak{M},t \vDash \phi \Longleftrightarrow \mathfrak{M},\{w\}\vDash \phi \text{ for all } w\in t. $$
\end{proposition}

\noindent The following  Disjunctive Normal Form Theorem  \cite[Thm. 4.9]{ciardelli2020} allows us to express every inquisitive and dependence formulas as a disjunction of standard formulas.

\begin{proposition}[Disjunctive Normal Form]\label{disjnf}
	Let $\phi \in \langInqI$, then there are standard inquisitive formulas $\alpha_0,\dots,\alpha_n\in \langInqst$ such that $\phi \equiv \bigvee \{ \alpha _i  \}_{i\leq n}$.
\end{proposition}

%\dnote{should I put here a remark about failure of uniform substitution? or maybe it is clear already}

%\noindent Finally, notice that the previous propositions also hold for all inquisitive formulas in $\langInt$ as we have that $\langInt\subseteq\langInqI$. 

\section{Algebraic Semantics for Inquisitive and Dependence Logic}\label{two}

We introduce in this section algebraic semantics of intermediate inquisitive and dependence logics and we prove its soundness and completeness. We first define inquisitive and dependence algebras -- $\inqI$-algebras and $\inqI^\otimes$-algebras -- and we show they are elementary structures axiomatised by universal Horn formulas. We then introduce so-called core semantics over such algebras and we prove using free algebras that our semantics is complete with respect to every intermediate inquisitive and dependence logic. 

%We conclude this section with some remarks on the algebraisability of inquisitive and dependence logics and we try to connect our work to the standard setting of abstract algebraic logic. 

\subsection{Inquisitive Algebras and Dependence Algebras}

Algebras are usually defined as a set together with some operations, i.e. as structures in an exclusively functional signature \cite{Burris.1981}. In order to provide a semantics to inquisitive and dependence logics we need to part ways from this definition and make space for a less restricted notion of algebras. In particular, we define inquisitive and dependence algebras in an expanded signature, consisting of functional symbols together with a unary predicate. Inquisitive and dependence algebras should be thus understood, from a model-theoretic perspective, as structures interpreting an algebraic language expanded by a unary predicate symbols.

We use $\mathcal{L}$ to refer to an arbitrary first-order language, and we use calligraphic letters $\mathcal{A}, \mathcal{B},\dots$ to denote first-order structures, in particular we shall use calligraphic letters to refer to standard and inquisitive algebras. If $\mathcal{A}$ is a structure, then we write $\mathsf{dom}(\mathcal{A})$ to refer to its underlying domain or universe. However, we shall often use the same symbol to denote a structure and its underlying universe. For all functional symbols $f\in\mathcal{L}$ and all relational symbols $R\in\mathcal{L}$, we write $f^{\mathcal{A}}$ and $R^{\mathcal{A}}$ for their interpretation in $\mathcal{A}$. However, when it is not confusing, we abide with the usual conventions and use the same notation for symbols and their interpretation.

An $\mathcal{L}$-structure $\mathcal{B}$ is a \textit{substructure} of an $\mathcal{L}$-structure $\mathcal{A}$ if (i) $\mathsf{dom}(\mathcal{B})\subseteq \mathsf{dom}(\mathcal{A})$, (ii) for all n-ary functional symbols $f\in \mathcal{L}$ and n-tuples $(b_1,\dots,b_{n})\in \mathcal{B}^n$, $f^{\mathcal{A}}(b_1,\dots,b_n)\in\mathcal{B}$, and (iii) for all relational symbols $R\in \mathcal{L}$ of arity $n$, $R^{\mathcal{B}}= R^{\mathcal{A}} \cap \mathsf{dom}(\mathcal{A})^n$.

If $ \mathcal{A} $ is an $\mathcal{L}$-structure and $X\subseteq \mathsf{dom}(\mathcal{A})$, then we denote by $\langle X\rangle $  the smallest substructure of $ \mathcal{A} $ containing $X$, i.e. the closure of $X$ under all the functional operations of $\mathcal{L}$. If we are interested in the closure of $X$ in $\mathcal{A}$ only with respect to some specific operations $f_0,\dots f_n \in \mathcal{L}$, then we write  $\langle X \rangle_{(f_0,\dots f_n)}$. 

Before defining inquisitive algebras, let us recall some well-known algebraic structures. A \textit{Brouwerian semilattice} $ \mathcal{B} $ is a bounded join-semilattice lattice with an extra-operation $\to$ such that for all $a,b,c\in \mathcal{B}$:
$$ a\land b \leq c \Longleftrightarrow a \leq b\to c. $$
\noindent A \textit{Heyting algebra} $\mathcal{H}$ is a bounded distributive lattice with an extra operation $\to$ satisfying the former equivalence. Given an element $a\in \mathcal{B}$, where  $ \mathcal{B} $ is a Brouwerian semilattice, we define its \textit{pseudocomplement} $\neg a$ as $\neg a := a \rightarrow 0$. Pseudocomplements of Heyting algebras are defined analogously. If $\mathcal{H}$ is a Heyting algebra such that for all $a\in \mathcal{H}$ it is the case that $a \land \neg a = 0$ and 	$a \lor \neg a = 1$, then we say that $\mathcal{H}$ is a \textit{Boolean algebra}. We define inquisitive algebras as follows.

\begin{definition}[Inquisitive Algebra]\label{inqui}
	An  (intuitionistic) \textit{inquisitive algebra} (or $\inqI$-algebra) is a structure $\mathcal{A}=(A, A_c, \land, \lor,\to, 0)$ in the vocabulary $\langInt\cup \{A_c \}$, such that:
	\begin{itemize}
		\item $A_c \subseteq A$;
		
		\item $(\langle A_c\rangle, \lor, \land, \rightarrow, 0 )$ is a Heyting algebra, where $\langle A_c\rangle$ is the closure of $A_c$ under the operations $\{\lor, \land, \rightarrow, 0\}$;
		
		\item  $(A_c, \land, \rightarrow, 0)$ is a Brouwerian semilattice;
		
		\item For all $x,y,z\in \langle A_c \rangle$ and $a\in A_c$, the following equation holds:
		\begin{align*}
		(\textit{Split})&  \hspace{10pt} a \rightarrow (x \lor y) = (a \rightarrow x)\lor (a \rightarrow y).
		\end{align*}
	\end{itemize}
\end{definition}

As we have remarked above, inquisitive algebras are algebras in a slightly non-standard sense: while algebras are usually defined as first-order structures in a purely algebraic signature, here we are expanding the signature by a unary predicate $A_c$, which we interpret as a \enquote{signed} subset of the algebra. The motivation for the addition of this predicate is that it captures at a semantical level the syntactic difference between standard formulas, which can be substituted freely, and non-standard formulas, for which uniform substitution fail.

Notice that, with slight abuse of notation, we write $A_c$ both for the predicate symbol in the language and for the corresponding subset of $A$. Given an inquisitive algebra $\mathcal{A}$, we generally refer to  this signed subset as the \textit{core} of $\mathcal{A}$, and we also denote it by $\mathsf{core}(\mathcal{A})$. By our definition, the core of a $\inqI$-algebra forms a Brouwerian semilattice in the  signature $\{\land, \rightarrow, 0 \}$. We write $\mathcal{A}_c$ for the Brouwerian semilattice $(A_c, \land, \rightarrow, 0)$ and we also write $\langle \mathcal{A}_c \rangle$ in place of $\langle A_c \rangle$.

Since $\langle \mathcal{A}_c \rangle$ is the closure of $A_c$ under all operations in $\{\land, \rightarrow, \lor, 0 \}$, it follows that $\mathcal{A}_c$ is a subalgebra of both  $\mathcal{A}$ and $\langle \mathcal{A}_c \rangle$ with respect to the reduct $\{\land, \rightarrow, 0 \}$. Negation is defined as $\neg x := x\rightarrow 0$ and the top element is $1:=0\to0$. Therefore, $\mathcal{A}$ and $\langle \mathcal{A}_c \rangle$ also agree on their interpretation of negation and $1$. 

It is important to stress that in our definition of inquisitive algebra every requirement and equation has a limited scope, i.e. they refer to elements of $A_c$ or $\langle A_c\rangle$ and not to arbitrary elements of $A$. This means that inquisitive algebras are somehow underspecified in their structure. Although this might seem as a downside of our definition, it is meant to reflect the distinction between standard and non-standard formulas in inquisitive and dependence logics, and the fact that every formula is inductively obtained from standard ones.

%fact that inquisitive logics do not allow for arbitrary substitution. %As we shall see, this also explains the key role that so-called core-generated algebras play in the semantics of inquisitive logic.

Dependence algebras are defined in the expanded signature $\langInqI\cup\{A_c\}$ in a similar fashion.

\begin{definition}[Dependence Algebra]\label{depe}
	An  (intuitionistic) \textit{dependence algebra} (or $\inqI^\otimes$-algebra) is a structure $\mathcal{A}=(A, A_c, \land, \otimes,  \lor, \to 0)$ in the language $\langInqI\cup\{A_c\}$, such that:
	\begin{itemize}
		\item $A_c \subseteq A$;
		
		\item $(\langle A_c\rangle, \lor, \land, \rightarrow, 0 )$ is a Heyting algebra, where $\langle A_c\rangle$ is the closure of $A_c$ under the operations $\{\lor, \otimes, \land, \rightarrow, 0\}$;
		
		\item  $(A_c, \otimes, \land, \rightarrow, 0)$ is a Heyting algebra;
		
		\item For all $x,y,z,k\in \langle A_c \rangle$ and $a\in A_c$, the following equations hold:
		\begin{align*}
		(\textit{Split})&  \hspace{10pt} a \rightarrow (x \lor y) = (a \rightarrow x)\lor (\alpha \rightarrow y);				\\
		(\textit{Dist})& \hspace{10pt}  x \otimes (y \lor z) = (x\otimes y) \lor (x\otimes z); \\
		(\textit{Mon})& \hspace{10pt}  (x\to z) \to (y \to k) = (x\otimes y) \to (z\otimes k).
		\end{align*}
	\end{itemize}
\end{definition}

\noindent It is clear from our definition that $\inqI^\otimes$-algebras are $\inqI$-algebras with an extra tensor operator $\otimes$, which satisfies the axioms \textit{Dist} and \textit{Mon} and whose core forms a Heyting algebra. We use the same conventions as for inquisitive algebras to refer to the underlying universe and to the core of a dependence algebra.

Given our previous considerations on the non-standard definition of such structures, one may wonder whether inquisitive and dependence algebras are structures in the first-order meaning of the word. To see that this is the case, it suffices to notice that we can use the predicate $A_c$ to express quantification over core elements, and we can use first-order-terms over $A_c$ to encode quantification over elements of $\langle A_c\rangle$. In this way it is straightforward to translate the definitions above into a list of first-order sentences and see that inquisitive and dependence algebras are elementary classes. To avoid confusion with the algebra operators, we use \& and $\supset$ as the first-order symbols of  conjunction and implication. We use $\tau,\sigma,\rho,\tau',\sigma'$ to denote arbitrary terms in the vocabulary $\langInt$. We use the abbreviations $\vec{x}:=x_0,\dots, x_n$ and $A_c(\vec{x}):=\bigwedge\nolimits_{i\leq n}A_c(x_i) $.

\begin{proposition}\label{firstorderaxiomat}
	
	(i) A structure $\mathcal{A}=(A, A_c, \land, \rightarrow, \lor, 0)$ is an inquisitive algebra if and only if it satisfies the following axioms and axiom schemas:
	\begin{align*}
	(1)	&\hspace{10pt}\forall x \forall y\:[A_c(x) \: \& \:A_c(y) \supset  A_c(x\land y)] \\   
	(2)	&\hspace{10pt}\forall x \forall y\:[A_c(x) \: \& \: A_c(y) \supset A_c(x\rightarrow y)] \\
	(3)	&\hspace{10pt} A_c(0) \\
	(4)	&\hspace{10pt}\forall \vec{x} \: [A_c(\vec{x})\supset  (\tau(\vec{x}) \rightarrow \tau(\vec{x})=1)]   \\
	(5)	&\hspace{10pt}\forall \vec{x}\:[A_c(\vec{x})\supset  (\tau(\vec{x}) \land 0 =0 )]   \\
	(6)	&\hspace{10pt}\forall \vec{x}\:[A_c(\vec{x}) \supset ( \tau(\vec{x}) \rightarrow (\sigma (\vec{x})  \land \rho(\vec{x} ))) = (\tau(\vec{x}) \rightarrow \sigma (\vec{x}) ) \land (\tau(\vec{x}) \rightarrow \rho (\vec{x}) )    ]   \\
	(7)	&\hspace{10pt}\forall \vec{x}\:[A_c(\vec{x})\supset  (\tau( \vec{x}) \land (\tau(\vec{x})\to \sigma(\vec{x}))) =( \tau( \vec{x}) \land \sigma(\vec{x}))]  \\
	(8)	&\hspace{10pt}\forall \vec{x}\:[A_c(\vec{x})\supset  (\sigma( \vec{x}) \land (\tau(\vec{x})\to \sigma(\vec{x}))) =  \sigma(\vec{x})]  \\
	(9)	&\hspace{10pt}\forall \vec{x}\:[A_c(\vec{x})\supset  ( \sigma(\vec{x}) \lor \tau(\vec{x}))   =  ( \tau(\vec{x}) \lor \sigma(\vec{x}))]  \\
	(10)	&\hspace{10pt}\forall \vec{x}\:[A_c(\vec{x})\supset  ( \tau(\vec{x}) \lor (\sigma(\vec{x}) \lor \rho(\vec{x})))   =   (\tau(\vec{x}) \lor \sigma(\vec{x}))   \lor \rho(\vec{x})]  \\
	(11)	&\hspace{10pt}\forall \vec{x}\:[A_c(\vec{x})\supset  ( \tau(\vec{x}) \lor \tau(\vec{x})   =   \tau(\vec{x})  )]  \\
	(12)	&\hspace{10pt}\forall \vec{x}\:[A_c(\vec{x})\supset  ( \tau(\vec{x}) \lor (\tau(\vec{x}) \land \sigma(\vec{x})))   =   \tau(\vec{x})  ]  \\
	(13)	&\hspace{10pt}\forall \vec{x}\:[A_c(\vec{x})\supset  ( \tau(\vec{x}) \land (\sigma(\vec{x}) \lor \rho(\vec{x})) )  =   (\tau(\vec{x}) \land \sigma(\vec{x}))   \lor (\tau(\vec{x}) \land \rho(\vec{x}))]  \\
	(14)	&\hspace{10pt}\forall \vec{x}  \:\forall y[(A_c(\vec{x}) \: \& \: A_c(y)	  \supset ( y \rightarrow (\sigma (\vec{x})  \lor \rho(\vec{x} ))) = (y \rightarrow \sigma (\vec{x}) ) \lor (y \rightarrow \rho (\vec{x}) )    ] . 
	\end{align*}
	(ii) A  structure $\mathcal{A}=(A, A_c, \land, \lor, \otimes, \to 0)$ is a dependence algebras if and only if it satisfies the axioms and axiom schemas $(1)-(14)$ above and the following:		
	\begin{align*}
	(15)	&\hspace{10pt}\forall \vec{x} \: \forall y\:[A_c(\vec{x})	  \supset  \tau(\vec{x}) \otimes (\sigma (\vec{x})  \lor \rho(\vec{x} )) = (\tau(\vec{x}) \otimes \sigma (\vec{x}) ) \lor (\tau(\vec{x}) \otimes \rho (\vec{x}) )    ]  \\ 
	(16)	&\hspace{10pt}\forall \vec{x} \: \forall y\:[A_c(\vec{x})	  \supset    (\tau(\vec{x}) \to \tau'(\vec{x})) \to(\sigma (\vec{x})  \to \sigma'(\vec{x} )) = (\tau(\vec{x}) \otimes \sigma (\vec{x})) \to (\tau'(\vec{x}) \otimes \sigma' (\vec{x}) )    ] . 
	\end{align*}
\end{proposition}
\begin{proof}
	Firstly we notice that the lists of axioms above are just an immediate translation of Definitions \ref{inqui} and \ref{depe} into first-order logic.
	
	(i) $(\Rightarrow)$	Axioms 4--13 guarantee that $\langle A_c \rangle$ is a Heyting algebra. Axioms 1--3 make sure that $A_c$ is closed under $\{\land, \to, 0\}$, hence it is a Bouwerian semilattice.  Finally, Axiom 14 corresponds to the \textit{Split} axiom. $(\Leftarrow)$ Analogous.	
	
	(ii) Immediate by (i) and the correspondence of Axiom 15 to \textit{Dist} and Axiom 16 to \textit{Mon}.
\end{proof}

\noindent	Since every formula in the lists above is a Horn formula, it follows that inquisitive and dependence algebras are elementary classes axiomatised by Horn formulas. Notice, however, that since axioms 4 -- 16 are schemas, with $\tau,\sigma,\rho,\tau',\sigma'$ being arbitrary terms, it follows that such axiomatisation is not finite.

Since $\inqI$-algebras and $\inqI^\otimes$-algebras are first-order structures, we can apply to our context the usual model-theoretic definitions of embedding, isomorphism, etc. In particular, we will often be interested in homomorphisms between inquisitive or dependence algebras.  We say that a function $h: \mathcal{A}\rightarrow \mathcal{B}$ between two  $\inqI$-algebras is an \textit{inquisitive homomorphism}, or a $\inqI$-homomorphism,  if $h$ commutes with the operators $\land,\lor,\rightarrow,\otimes,0 $ and, in addition, we have that $h[\mathcal{A}_c]\subseteq \mathcal{B}_c$. A function $h: \mathcal{A}\rightarrow \mathcal{B}$ between two  $\inqI^\otimes$-algebras is a \textit{dependence homomorphism}, or a $\inqI^\otimes$-homomorphism, if it is a inquisitive homomorphism which also preserves the tensor operation, i.e. $h(x\otimes y)= h(x)\otimes h(y)$ for all $x,y\in \mathcal{A}$.

We  denote by $ \mathsf{InqAlg} $ the category of inquisitive algebras with $\inqI$-homomorphisms  and by $ \mathsf{InqAlg}^\otimes $ the category of inquisitive algebras with $\inqI^\otimes$-homomorphisms.

\subsection{Core Semantics and Translation into Horn Formulas}

Now that we have defined inquisitive and dependence algebras we can use them in order to give suitable algebraic semantics to inquisitive and dependence logics. To this end, we introduce \textit{core semantics}, where valuations are restricted to range over a subset of a structure.

Let $\mathcal{A}$ be an inquisitive or dependence algebra, we say that a function $\mu$ is a \textit{core valuation} over $\mathcal{A}$ if it assigns atomic formulas from $\at$ to elements in $A_c$, i.e. $\mu:\at\to\mathcal{A}_c$. Similarly, if $\mathcal{A}$ is a dependence algebra, a \textit{core valuation} is a function $\mu:\at\to\mathcal{A}_c$.

\begin{definition}[Algebraic Model]
	An \textit{inquisitive algebraic model} is a pair $\mathcal{M}=(\mathcal{A},\mu)$ where  $\mathcal{A}$ is an inquisitive algebra and $\mu:\at\rightarrow \mathcal{A}_c$ is a core valuation.
	A \textit{dependence algebraic model} is an inquisitive algebraic model $\mathcal{M}=(\mathcal{A},\mu)$ where $\mathcal{A}$ is also a dependence algebra.
\end{definition}

The interpretation of arbitrary formulas in an algebraic model $\mathcal{M}$ is  defined recursively as follows. Notice that this definition is standard, besides for the fact that atomic formulas can be assigned only to core elements of the underlying algebra. 

\begin{definition}[Interpretation of Arbitrary Formulas]
	Given an inquisitive algebraic model $\mathcal{M}$ and a formula $\phi\in\langInt$, its \textit{interpretation} $\llbracket \phi \rrbracket^{\mathcal{M}}$ is defined as follows:
	\begin{equation*}
	\begin{array}{r@{\hspace{.3em}}c@{\hspace{.3em}}l  @{\hspace{1.5em}}  r@{\hspace{.3em}}c@{\hspace{.3em}}l  @{\hspace{1.5em}}  r@{\hspace{.3em}}c@{\hspace{.3em}}l}
	\llbracket p \rrbracket^\mathcal{M} &= &\mu(p)
	&\llbracket \bot \rrbracket^\mathcal{M} &= &0
	&\llbracket \phi \lor \psi \rrbracket^M &= &\llbracket \phi \rrbracket^{\mathcal{M}} \lor \llbracket \psi \rrbracket^{\mathcal{M}}\\[.5em]
	\llbracket \phi \land \psi \rrbracket^\mathcal{M} &= &\llbracket \phi \rrbracket^\mathcal{M} \land \llbracket \psi \rrbracket^\mathcal{M}
	&\llbracket \phi \to \psi \rrbracket^\mathcal{M} &= &\llbracket \phi \rrbracket^\mathcal{M} \to \llbracket \psi \rrbracket^\mathcal{M}.
	\end{array}
	\end{equation*} Moreover, if $\mathcal{M}$ is a dependence model and $\phi\in\langInqI$ , then its interpretation is defined by the clauses above together with the following one:
	$$\llbracket \phi \otimes \psi \rrbracket^\mathcal{M} = \llbracket\phi \rrbracket^{\mathcal{M}} \otimes  \llbracket\psi \rrbracket^{\mathcal{M}}. $$
\end{definition}

\noindent If $\phi$ is a formula and $\mathcal{M}$ is an (inquisitive or dependence) algebraic model, we also abbreviate $\phi^{\mathcal{M}}$ for the interpretation of $\phi$ in the model $\mathcal{M}$.     We  write $\mathcal{M}\vDash^c \phi$ and say that $\phi$ is true in $\mathcal{M} $ if $ \phi^{\mathcal{M}}=1$. We say that $\phi$ is valid in the inquisitive (or dependence) algebra $\mathcal{A} $ and write $\mathcal{A}\vDash^c \phi$ if $\phi$ is true in every model $\mathcal{M}=(\mathcal{A},\mu)$ over $\mathcal{A}$. If $\mathsf{C}$ is a class of inquisitive (or dependence) algebras, then we say $\phi$ is \textit{valid in }$\mathsf{C}$, and write $\mathsf{C}\vDash^c \phi$, if  $\mathcal{A}\vDash^c\phi$ for all $\mathcal{A}\in \mathsf{C}$. Finally, we say that $\phi$ is  an \textit{algebraic validity} of inquisitive (or dependence) logic if it is true in all inquisitive (or dependence) models.

Before we define arbitrary intermediate inquisitive and dependence algebras, we shall first explain how to relate validity under core-semantics to standard first-order validity. Recall that an \textit{equation} in a language $\mathcal{L}$ is an atomic first-order formula of the form $\epsilon= \delta$, where $\epsilon$ and $\delta$ are terms in $\mathcal{L}$. Notice in particular that any propositional formula in $\langInt$ is a term in $\langInt\cup\{A_c\}$ and \textit{vice versa} -- and clearly the same holds for formulas in $\langInqI$ and terms in $\langInqI\cup\{A_c\}$. We can thus associate every formula $\phi$ to a corresponding equation $\Theta(\phi):= \phi = 1$. Similarly, to any equation $\epsilon= \sigma$ -- with $\epsilon,\sigma$ terms in $\langInt$ or $\langInqI$ -- we associate the formula $\Delta( \epsilon,\sigma):=\epsilon\leftrightarrow\sigma$, where $\epsilon\leftrightarrow\sigma$ is a shorthand for $\epsilon\to\sigma\land \sigma\to\epsilon$.

\begin{proposition}\label{Trans}
	Let $\mathcal{A}$ be an inquisitive (or dependence) algebra. Then, for all formulas $\phi(x_0,\dots,x_n)$, and all equations $\epsilon =\delta$, we have that:
	\begin{align}
	\mathcal{A}\vDash^c \phi & \Longleftrightarrow \mathcal{A}\vDash \forall x_0,\dots,\forall x_n \Big(  \bigwedge_{i\leq n}A_c(x_i) \supset \Theta(\phi)  \Big) ;\\
	\mathcal{A}\vDash^c \Delta(\epsilon,\delta) & \Longleftrightarrow \mathcal{A}\vDash \forall x_0,\dots,\forall x_n \Big( \bigwedge_{i\leq n}A_c(x_i) \supset \epsilon(x_0,\dots,x_n)= \delta(x_0,\dots,x_n) \Big).		
	\end{align}
\end{proposition}
\begin{proof}
	The proof is analogous for inquisitive and dependence algebras. Both (1) and (2) follow by a straightforward induction. We prove the base case only and leave the rest to the reader.
	
	Consider first (1).	Let $\phi=p\in \at$, then we have that if $\mathcal{A}\vDash^c p$, then for for all \textit{core valuations} $\mu:\at\to \mathcal{A}_c$, we have that    $\mu(p)=1$. Hence, for all $x\in A_c$, we have that $x=1$, therefore $ \mathcal{A}\vDash \forall p \Big(  A_c(p) \supset p= 1  \Big)  $. The other direction follows analogously.
	
	Consider (2). Let $\epsilon=p$ and $\delta=q$. Then if $\mathcal{A}\vDash^c p\leftrightarrow q$, then for all core-valuations $\mu:\at\to \mathcal{A}$, we have that $\mathcal{M}=(\mathcal{A},\mu)\vDash^c p\leftrightarrow q$, which means that $\llbracket p \leftrightarrow q \rrbracket^\mathcal{M}$. Since $p^{\mathcal{M}},q^{\mathcal{M}}\in A_c$ and $\langle A_c \rangle $ is a Heyting algebra, this means that $p^{\mathcal{M}}\leq q^{\mathcal{M}}$ and $q^{\mathcal{M}}\leq p^{\mathcal{M}}$, which together entail $p^{\mathcal{M}}= q^{\mathcal{M}}$. Finally, this means that for all $x,y\in A_c$, we have that $x=y$, hence $ \mathcal{A}\vDash \forall x\forall y \Big(  A_c(x)\land A_x(y) \supset x= y  \Big)  $. The other direction follows analogously.
\end{proof}

\noindent This theorem gives us an important bridge between standard and core semantics. In particular, it shows that the truth of inquisitive (and dependence) formulas is equivalent to the validity under the standard Tarski semantics of a corresponding universal Horn formula.

We now use core-semantics to define arbitrary intermediate inquisitive and dependence algebras.

\begin{definition}$\:$
	\begin{itemize}
		\item Let $\Lambda\subseteq \langIntsta$ be a set of formulas closed under standard substitution, then 	an  inquisitive algebra  $\mathcal{A}$ is said to be an $\inqint$-algebra if  $ \mathcal{A} \vDash^c \Lambda$.
		
		\item Let $\Lambda\subseteq \langInqst$ be a set of formulas closed under standard substitution,	then a dependence algebra  $\mathcal{A}$ is said to be an $\inqint^\otimes$-algebra if  $ \mathcal{A} \vDash^c \Lambda$.
	\end{itemize}
\end{definition}

\noindent We then say that $\mathcal{A}$ is an intermediate inquisitive algebra if it is an $\inqint$-algebra for some $\Lambda\subseteq\langIntsta$. Similarly, $\mathcal{A}$ is an intermediate dependence algebra if it is an $\inqint^\otimes$-algebra for some $\Lambda\subseteq\langInqst$.   In particular, if $\Lambda= \{\neg\neg \alpha \to \alpha\}_{\alpha\in \langIntsta}$ and $ \mathcal{A} \vDash^c \Lambda$ then we say that  $\mathcal{A}$ is a $\inqB$-algebra. Similarly, if $\Lambda= \{\neg\neg \alpha \to \alpha\}_{\alpha\in \langInqst}$ and $ \mathcal{A} \vDash^c \Lambda$ then we say that  $\mathcal{A}$ is a $\inqB^\otimes$-algebra. It is straightforward to verify that  $\inqB^\otimes$-algebras are those dependence algebras whose core is a Boolean algebra with $\otimes$ as their join operator.  The class of $\inqB$-algebras defined here strictly extends the class of (classical) inquisitive algebras considered in \cite{Quadrellaro.2019B}. However, as we shall see later, the subclass of core-generated  $\inqB$-algebras coincides with the so-called regular inquisitive algebras of \cite{Quadrellaro.2019B}.

By Proposition \ref{Trans} we have that $\inqint$-algebras and  $\inqint^\otimes$-algebras are elementary classes of structures. In particular, since every formula in their axiomatisation is a Horn formula, it also follows that $\inqint$-algebra are Horn-axiomatisable. The following proposition is then easy to prove.

\begin{proposition}[Soundness]\label{soundness}
	$\:$
	\begin{itemize}
		\item[(i)] If $\mathcal{A}$ is a $\inqint$-algebra, then $\mathcal{A}\vDash^c \inqint$. 
		\item[(ii)] If $\mathcal{A}$ is a $\inqint^\otimes$-algebra, then $\mathcal{A}\vDash^c \inqint^\otimes$.
	\end{itemize}
	
\end{proposition}
\begin{proof}
	(i) Let  $\mathcal{A}$ is a $\inqint$-algebra, then by Proposition \ref{Trans} and  axioms (1)--(14) from Proposition \ref{firstorderaxiomat} it follows immediately that $ \mathcal{A} \vDash^c \inqI$, since the axioms in Definition \ref{inqI} correspond to those of Proposition \ref{firstorderaxiomat}. We leave to the reader to check this correspondence. Moreover, since $\inqI$-algebras are closed under modus ponens and by assumption  $\mathcal{A}\vDash^c\Lambda$ it follows  $ \mathcal{A} \vDash^c \inqint$.
	
	(ii) If  $\mathcal{A}$ is a $\inqint^\otimes$-algebra then our claim follows in the same way, by the fact that the axioms (1)--(16) of Proposition \ref{firstorderaxiomat} are equivalent by Proposition \ref{Trans} to  the propositional axioms of Definition \ref{inqIotimes}. 	
\end{proof}

\noindent  We  denote by $\mathsf{InqAlg\Lambda}$ the category (and the class) of $\inq\Lambda$-algebras with $\inqI$-homomorphisms. Similarly, we  denote by $\mathsf{InqAlg\Lambda^\otimes}$ the category (and the class) of $\inqI^\otimes$-algebras with $\inqI^\otimes$-homomorphisms.

We conclude this section proving some closure properties of the validities of formulas under core-semantics. If $\mathcal{B}$ is an inquisitive (dependence) substructure of $\mathcal{A}$ we also say that $\mathcal{B}$ is a \textit{subalgebra} of $\mathcal{A}$ and we write $\mathcal{B}\preceq \mathcal{A}$ -- notice that this does not mean that $\mathcal{B}$ is an elementary substructure of $\mathcal{A}$. Finally, if $\mathcal{A}\preceq \mathcal{B}$ and $\mathcal{A}_c= \mathcal{B}_c$, then we say that $\mathcal{B}$ is a \textit{core-superstructure} of $\mathcal{A}$. We now prove the two following closure properties.

\begin{proposition}\label{closure1}
	Let $\mathcal{A}, \mathcal{B}$ be inquisitive (or dependence) algebras, then:
	\begin{itemize}
		\item[(i)] $\mathcal{A}\vDash^c \phi$ and $\mathcal{B}\preceq \mathcal{A}$ entail $\mathcal{B}\vDash^c \phi$;
		\item[(ii)] $\mathcal{A}\vDash^c \phi$ and $\mathcal{A}\preceq \mathcal{B}$, $\mathcal{A}_c= \mathcal{B}_c$ entail $\mathcal{B}\vDash^c \phi$.
	\end{itemize}
\end{proposition}
\begin{proof}
	(i) Suppose $\mathcal{B}\nvDash^c  \phi(\vec{p})$, then $(\mathcal{B},\mu)\nvDash^c  \phi(\vec{p})$ for some core-valuation $\mu$. Since $\mathcal{B}\preceq \mathcal{A}$ and $\mu[\at]\subseteq \mathcal{B}_c$, we  have $\llbracket \phi(\vec{p}) \rrbracket^{(\mathcal{A},\mu)})=\llbracket \phi(\vec{p}) \rrbracket^{(\mathcal{B},\mu)}$, therefore $(\mathcal{A},\mu)\nvDash^c  \phi(\vec{p})$ and $\mathcal{A}\nvDash^c  \phi(\vec{p})$.
	
	(ii) Suppose $\mathcal{B}\nvDash^c  \phi(\vec{p})$ and also $\mathcal{A}\preceq \mathcal{B}$, $\mathcal{A}_c= \mathcal{B}_c$. Then we have $(\mathcal{B},\mu)\nvDash^c p$ for some core-valuation $\mu$ and, by $\mathcal{A}_c= \mathcal{B}_c$ it follows that $\mu$ is a core-valuation over $\mathcal{A}$ as well. Finally, since $\mathcal{A}\preceq \mathcal{B}$, we have that $\llbracket \phi(\vec{p}) \rrbracket^{(\mathcal{A},\mu)}=\llbracket \phi(\vec{p}) \rrbracket^{(\mathcal{B},\mu)}$, therefore $(\mathcal{A},\mu)\nvDash^c  \phi(\vec{p})$ and $\mathcal{A}\nvDash^c  \phi(\vec{p})$.
\end{proof}

\noindent The following closure properties of $\inqint$-algebras and $\inqint^\otimes$-algebras follow directly from the previous proposition.

\begin{corollary}\label{closure2}
	$\:$
	\begin{itemize}
		\item For every $\Lambda\subseteq \langIntsta$, $\inqint$-algebras are closed under subalgebras and core-superstructures.
		\item For every $\Lambda\subseteq \langInqst$, $\inqint^\otimes$-algebras are closed under subalgebras and core-superstructures.
	\end{itemize}
	
\end{corollary}

If $h:\mathcal{A}\to\mathcal{B}$ is a surjective $\inqI$-homomorphism (or $\inqI^\otimes$-homomorphism), we say that $\mathcal{B}$ is a homomorphic image of $\mathcal{A}$ and we write $h:\mathcal{A}\twoheadrightarrow\mathcal{B}$.  Since $\inqI$-homomorphisms and $\inqI^\otimes$-homomorphisms may map non-core elements to core elements, algebraic validities are not closed under homomorphic images. Let $\mathcal{A}$ and $\mathcal{B}$ be the following inquisitive algebras:

\medskip

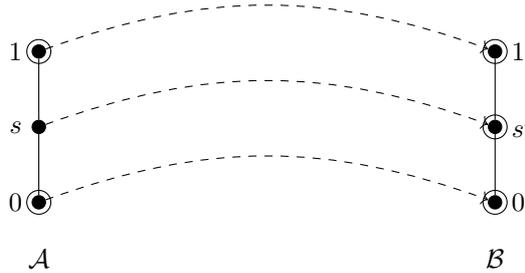
\begin{center}
	
	\begin{tikzpicture}
	
	\node[highlight] at ( 0,0) {};
	\node[highlight] at ( 0,2) {};

	\node[dot,label=left:$0$] (0) at (0,0) {};
	\node[dot,label=left:$s$] (s) at (0,1) {};
	\node[dot,label=left:$1$] (1) at (0,2) {};
	
	\draw (0) -- (s);
	\draw (s) -- (1);
	
	\node[below] at (0,-0.5) {$\mathcal{A}$};
	
	\node[highlight] at ( 6,0) {};
	\node[highlight] at ( 6,1) {};
	\node[highlight] at ( 6,2) {};
	
	\node[dot,label={0:$0$}] (0') at ( 6,0) {};
	\node[dot,label={0:$s'$}] (s') at (6,1) {};
	\node[dot,label={0:$1$}] (1') at (6,2) {};
	
	\draw (0') -- (s');
	\draw (s') -- (1');
	
	\node[below] at (6,-0.5) {$\mathcal{B}$};
	
	\draw[dashed, ->] (0) to [out=20,in=160] (0');
	\draw[dashed,->] (s) to [out=20,in=160] (s');
	\draw[dashed,->] (1) to [out=20,in=160] (1');

	\end{tikzpicture}
	
\end{center}

%
%\begin{center}
%	
%	\begin{tikzpicture}
%	
%	\node[highlight] at ( 0,0) {};
%	\node[highlight] at ( 1,1) {};
%	\node[highlight] at ( 0,3) {};
%	\node[highlight] at ( -1,1) {};
%	
%	
%	\node[dot,label={0:$0$}] (0) at ( 0,0) {};
%	\node[dot,label={180:$a$}] (a) at (-1,1) {};
%	\node[dot,label={0:$b$}] (b) at ( 1,1) {};
%	\node[dot,label={0:$s$}] (s) at (0,2) {};
%	\node[dot,label={0:$1$}] (1) at (0,3) {};
%	
%	\draw (0) -- (a);
%	\draw (0) -- (b);
%	\draw (a) -- (s);
%	\draw (b) -- (s);
%	\draw (s) -- (1);
%	
%	\node[below] at (0,-0.5) {$\mathcal{A}$};
%	
%	\node[highlight] at ( 6,0) {};
%	\node[highlight] at ( 5,1) {};
%	\node[highlight] at ( 7,1) {};
%	\node[highlight] at ( 6,3) {};
%	\node[highlight] at ( 6,2) {};
%	
%	\node[dot,label={0:$0$}] (0') at ( 6,0) {};
%	\node[dot,label={180:$a'$}] (a') at (5,1) {};
%	\node[dot,label={0:$b'$}] (b') at ( 7,1) {};
%	\node[dot,label={0:$s'$}] (s') at (6,2) {};
%	\node[dot,label={0:$1$}] (1') at (6,3) {};
%	
%	\draw (0') -- (a');
%	\draw (0') -- (b');
%	\draw (a') -- (s');
%	\draw (b') -- (s');
%	\draw (s') -- (1');
%	
%	\node[below] at (6,-0.5) {$\mathcal{B}$};
%	
%	\draw[dashed, ->] (0) to [out=20,in=160] (0');
%	\draw[dashed, ->] (b) to [out=20,in=160] (b');
%	\draw[dashed,->] (s) to [out=20,in=160] (s');
%	\draw[dashed,->] (a) to [out=20,in=160] (a');
%	\draw[dashed,->] (1) to [out=20,in=160] (1');
%	
%	
%	\end{tikzpicture}
%	
%\end{center}

\noindent where $\mathsf{core}(\mathcal{A})=\{0,1\}$ and $\mathsf{core}(\mathcal{B})=\mathcal{B}= \{0,s',1\}  $. It is easy to verify that with such core both $\mathcal{A}$ and $\mathcal{B}$ satisfy the conditions of Definition \ref{inqui} and are inquisitive algebras. Then, the function $h:\mathcal{A}\to\mathcal{B}$ depicted in the picture above is clearly a $\inqI$-homomorphism, since it is the identity function over the algebraic reduct of $\mathcal{A}$ and $\mathcal{B}$ and moreover $h[\mathcal{A}_c]\subseteq\mathcal{B}_c$. Finally, one can readily check that $\mathcal{A}\vDash^c \neg\neg p \to p$ but $\mathcal{B}\nvDash^c \neg\neg p \to p$, as it is witnessed by the core assignment $\mu:p\mapsto s'$.

Hence, the validity of core formulas is not preserved by $\inqI$-homomorphisms and $\inqI^\otimes$-homomorphisms. However, we shall prove later that closure under homomorphic image holds in a restricted class of cases.

%RIPARTO DA QUI 12/01

\subsection{Free $\inqint$-Algebras and Algebraic Completeness}\label{free}

We introduce in this section free inquisitive and dependence algebras -- i.e. Lindenbaum-Tarski algebras for inquisitive and dependence logics -- in order to prove the completeness of the algebraic semantics that we presented above. We refer the reader to \cite{Burris.1981} and \cite{Font.2016} for the standard construction of free algebras.

Recall that $\at$ is countable set of atomic formulas, $\langInt$ is the set of all formulas of inquisitive logic and $\langInqI$ is the set of all formulas of dependence logic. It is useful here to think of $\langInt$ and $\langInqI$ as \textit{term algebras}, whose elements are formulas and whose operations are respectively $\{\land, \lor, \rightarrow, 0 \}$ and $\{\land, \lor, \rightarrow, \otimes, 0 \}$.

Free algebras are generally obtained by quotienting term algebras by suitable congruences. Here, we introduce a congruence relation for every intermediate inquisitive logic $\inqint$, and one for every intermediate dependence logic $\inqint^\otimes$. We define the relations $ \equiv_{\inqint} \subseteq \langInt\times \langInt$ and $ \equiv_{\inqint^\otimes}\subseteq \langInqI\times\langInqI $ as follows:
\begin{align*}
\phi \equiv_{\inqint} \psi &\Longleftrightarrow \phi \leftrightarrow \psi \in \inqint;\\
\phi \equiv_{\inqint^\otimes} \psi &\Longleftrightarrow \phi \leftrightarrow \psi \in \inqint^\otimes.
\end{align*}

\noindent It is easy to verify that these are equivalence relations. Moreover, since intermediate inquisitive (and dependence) logics are closed under \textit{modus ponens} one can also verify that $\equiv_{\inqint}$  and $\equiv_{\inqint^\otimes}$ are congruences over the term algebras $\langInt$ and $\langInqI$. 

%See \cite[\S 7.2]{Zakharyaschev.1997}.

Since our setting is non-standard, we need to define free algebras as first-order structures with a core.  Free algebras for $\inqint$ are obtained by first quotienting the term algebra $\langInt$ by the congruence relation $\equiv_{\inqint}$ -- in this way we obtain the Heyting algebra $\mathcal{F}_{\inqint}=(\langInt/\equiv_{\inqint}, \land, \lor, \rightarrow, 0)$. To turn such structure into an inquisitive algebra, we then need to specify its core. To this end, we say that an equivalence class $[\phi]\in \mathcal{F}_{\inqint}=\langInt/\equiv_{\inqint}$ is \textit{classical} if there is some $\alpha\in \langIntsta$ such that $\alpha\in[\phi]$. We denote by  $\mathcal{F}_{\Lambda}$ the set of classical equivalence classes in $\mathcal{F}_{\inqint}$. We proceed similarly for dependence algebras:  we say that an  equivalence class $[\phi]\in \mathcal{F}^\otimes_{\inqint}=\langInt^\otimes/\equiv_{\inqint^\otimes}$ is \textit{classical} if there is some $\alpha\in \langInqst$ such that $\alpha\in[\phi]$ and we denote this set by $\mathcal{F}^\otimes_{\Lambda}$. Then, to obtain the free algebra for $\inqint^\otimes$ we quotient the term algebra $\langInqI$ by the congruence relation $\equiv_{\inqint^\otimes}$ and  we let the set of classical equivalence classes $\mathcal{F}^\otimes_{\Lambda}$  be its core. We can thus define free inquisitive and dependence algebras as follows.	

\begin{definition}
	The \textit{Free Inquisitive Algebra} $ \mathcal{F}_{\inqint} $ of the intermediate inquisitive logic $\inqint$ is the first order structure $\mathcal{F}_{\inqint}=(\langInt/\equiv_{\inqint}, \mathcal{F}_{\Lambda}, \land, \lor, \rightarrow, 0)$, where $\mathcal{F}_{\Lambda}$ is the set of classical equivalence classes in $\langInt/\equiv_{\inqint}$.
	The \textit{Free Dependence Algebra} $ \mathcal{F}_{\inqint^\otimes} $ of the intermediate inquisitive logic $\inqint$ is the first order structure $\mathcal{F}_{\inqint^\otimes}=(\langInqI/\equiv_{\inqint^\otimes}, \mathcal{F}^\otimes_{\Lambda}, \land, \lor, \otimes, \rightarrow, 0)$, where $\mathcal{F}^\otimes_{\Lambda}$ is the set of classical equivalence classes in $\langInqI/\equiv_{\inqint}$.
\end{definition}

\begin{proposition}\label{correct}
	Let $\inqint$ be an intermediate inquisitive logic and $\inqint^\otimes$ an intermediate dependence logic, then $\mathcal{F}_{\inqint}$ is  an inquisitive algebra and $\mathcal{F}_{\inqint^\otimes}$  a dependence algebra.
\end{proposition}
\begin{proof}
	(i) Since $\ipc\subseteq \inqint$ it follows immediately that $\mathcal{F}_{\inqint}$ is a Heyting algebra. Moreover, since standard formulas are closed under meet and implication, it follows that if $[\alpha],[\beta]\in\mathcal{F}_{\Lambda}$, then $[\alpha\land\beta]\in\mathcal{F}_{\Lambda}$ and $[\alpha\to\beta]\in\mathcal{F}_{\Lambda}$. It is clear that $[\bot]\in\mathcal{F}_{\Lambda}$. Hence it follows that $ \mathcal{F}_{\Lambda} $ is a Brouwerian semilattice. Also, for every $\alpha\in\langIntsta$ and all $\phi,\psi\in\langInt$ we have by Axiom (A10) of Definition \ref{inqI} that $(\alpha \rightarrow (\phi \lor \psi)) \rightarrow ((\alpha \rightarrow \phi)\lor (\alpha \rightarrow \psi))\in \inqint$. Moreover, since $(\alpha \rightarrow (\phi \lor \psi)) \leftarrow ((\alpha \rightarrow \phi)\lor (\alpha \rightarrow \psi))\in \ipc$, it follows that:
	$$(\alpha \rightarrow (\phi \lor \psi)) \leftrightarrow ((\alpha \rightarrow \phi)\lor (\alpha \rightarrow \psi))\in \inqint.$$
	\noindent Hence,
	$$\alpha \rightarrow (\phi \lor \psi)\equiv_{\inqint} (\alpha \rightarrow \phi)\lor (\alpha \rightarrow \psi).$$
	\noindent  Therefore, for every $[\alpha]\in\mathcal{F}_{\Lambda}$ and all $[\phi],[\psi]\in\mathcal{F}_{\inqint}$ we have that:
	$$[\alpha] \rightarrow ([\phi] \lor [\psi]) = ([\alpha] \rightarrow [\phi])\lor ([\alpha] \rightarrow [\psi]),$$
	\noindent proving that $\mathcal{F}_{\inqint}$ satisfies \textit{Split}. Thus $\mathcal{F}_{\inqint}$ is an inquisitive algebra.
	
	(ii) It is proven analogously to (i), by checking that \textit{Mon} and \textit{Dist} hold in $\mathcal{F}_{\inqint^\otimes}$.
\end{proof}

\noindent The following proposition shows that  every inquisitive logic $\inqint$ is the logic of a free algebra.

\begin{proposition}\label{freeli}
	Let $\inqint$ be an intermediate inquisitive logic and $\inqint^\otimes$ an intermediate dependence algebra, then:
	\begin{align*}
	\phi\in\inqint& \Longleftrightarrow \mathcal{F}_{\inqint}\vDash^c \phi ; \\
	\phi\in\inqint^\otimes& \Longleftrightarrow \mathcal{F}_{\inqint^\otimes}\vDash^c \phi.
	\end{align*}
\end{proposition}
\begin{proof}
	We prove the claim only for intermediate inquisitive logics, as the proof for dependence logics is the same. $(\Leftarrow)$  Suppose $ \phi(p_0, \dots, p_n)\notin \inqint $ where $p_0,\dots,p_n$ are the propositional variables occurring in $\phi$. Then we have that $\phi \leftrightarrow \top \notin \inqint$, hence $\phi\not\equiv_{\inqint} \top$ and so $1_{\mathcal{F}_{\inqint}}\neq[\phi]$. Since $[\phi]=\phi([p_0],\dots,[p_n]  )$ and $p_0,\dots,p_n$ are standard formulas, we can define the (canonical) core-valuation $\mu:\at\to \mathcal{F}_{\Lambda}$ such that for all $p_i$ with $i\leq n$,   $\mu (p_i)= [p_i]$. It follows immediately:
	$$\phi(p_0,\dots,p_n)^{\mathcal{F}_{\inqint}} = \phi([p_0],\dots,[p_n]  )  = [\phi(p_0,\dots,p_n)] \neq  1_{\mathcal{F}_{\inqint}}.   $$
	\noindent Which means that  $\mathcal{F}_{\inqint}\nvDash^c\phi$ and thus proves our claim. $(\Rightarrow)$  Analogously to the previous direction.
\end{proof}

\begin{corollary}\label{correct2}
	Let $\inqint$ be an intermediate inquisitive logic and $\inqint^\otimes$ an intermediate dependence logic, then $\mathcal{F}_{\inqint}$ is  a $\inqint$-algebra and $\mathcal{F}_{\inqint^\otimes}$  a $\inqint^\otimes$-algebra.
\end{corollary}
\begin{proof}
	Immediate by Propositions \ref{freeli} and \ref{correct} together with the definition of  $\inqint$-algebras and $\inqint^\otimes$-algebras.
\end{proof}

Free $ \inqint $-algebras thus witness the validity of every formula. By this fact, the algebraic completeness of intermediate inquisitive logics follows immediately.

\begin{theorem}[Algebraic Completeness]\label{algcomple}
	Every intermediate inquisitive logic $\inqint$ is complete with respect to the class of $\inqint$-algebras and every intermediate dependence logic $\inqint^\otimes$ is complete with respect to the class of $\inqint^\otimes$-algebras:
	\begin{align*}
	\phi\in \inqint &\Longleftrightarrow \mathsf{InqAlg\Lambda} \vDash^c \phi; \\
	\phi\in \inqint^\otimes &\Longleftrightarrow \mathsf{InqAlg\Lambda^\otimes} \vDash^c \phi.
	\end{align*}
\end{theorem}
\begin{proof}
	We prove the claim for inquisitive algebras only, as the case for dependence algebra is exactly the same.
	$ (\Rightarrow)$ Suppose $ \phi\in \inqint $ and let $\mathcal{A}\in \mathsf{InqAlg\Lambda}$. Then by Proposition \ref{soundness} we immediately have that  $ \mathcal{A} \vDash^c \phi $ and therefore $\mathsf{InqAlg\Lambda} \vDash^c \phi$.
	($\Leftarrow$) 	Suppose by contraposition that $\phi\notin \inqint$, then by Proposition \ref{freeli} we have that $\mathcal{F}_{\inqint}\nvDash^c \phi$ and then, since by Corollary \ref{correct2} we also have that $\mathcal{F}_{\inqint} \in  \mathsf{InqAlg\Lambda}$, it follows  that $  \mathsf{InqAlg\Lambda}\nvDash^c \phi$.
\end{proof}

\noindent We thus have shown that the algebraic semantics we introduced for inquisitive and dependence logics is both sound and complete. Inquisitive and dependence logic can be then investigated not only from the point of view of team semantics, but also from an algebraic perspective.

\section{Properties of Inquisitive and Dependence Algebras}\label{three}

In this section we study several properties of inquisitive and dependence algebras. In particular, we try to find suitable subclasses of $\mathsf{InqAlg\Lambda}$ and $\mathsf{InqAlg\Lambda^\otimes}$ that witness the validity of inquisitive and dependence formulas. To this end, we introduce and investigate finite, core-generated and well-connected inquisitive and dependence algebras, and we prove several results concerning such structures. %We will later use the results of this Section to prove Representation Theorems for finite, core-generated and well-connected inquisitive and dependence algebras .

\subsection{Core-Generated and Well-Connected Inquisitive and Dependence Algebras}

If $\mathcal{A}$ is an arbitrary inquisitive (dependence) algebra, there is not much we can say about its structure with full generality, for the axioms of inquisitive (dependence) algebras characterise only the substructure $\langle \mathcal{A}_c\rangle$ of $\mathcal{A}$. Therefore, it is useful to focus our attention on \enquote{small} inquisitive (dependence) algebras, namely to those structures $\mathcal{A}$ which are generated by their core $\mathcal{A}_c$.

\begin{definition}
	An inquisitive or dependence algebra $\mathcal{A}$ is  \textit{core-generated} if $\mathcal{A}=\langle \mathcal{A}_c\rangle$.
\end{definition}

\noindent We shall see in this section that core-generated inquisitive algebras play an important role in the algebraic semantics of $\inqI$. We denote by $ \mathsf{InqAlg_{CG}} $ the category of core-generated inquisitive algebras with $\inqI$-homomorphisms and by $ \mathsf{InqAlg^\otimes _{CG}}$ the category of core-generated dependence algebras together with $\inqI^\otimes$-homomorphisms. The categories  $ \mathsf{InqAlg\Lambda_{CG}} $ and $ \mathsf{InqAlg\Lambda_{CG}^\otimes} $ are defined analogously.

We first prove the following Normal Form Theorem, which allows us to express every element of $\langle \mathcal{A}_c\rangle$ in the form of a disjunction of core elements. This theorem is really an algebraic counterpart of the normal form result for $\inqI$ proved in \cite{ciardelli2020} and recalled earlier in Section \S 1.3. 

\begin{theorem}[Disjunctive Normal Form]\label{normalformalgebra}
	Let $\mathcal{A}$ be any inquisitive (dependence) algebra, then for all $x\in \langle \mathcal{A}_c\rangle$  there are pairwise incomparable elements $a_0,...,a_n\in \mathcal{A}_c$ such that $x=\bigvee_{i\leq n} a_i$. 
\end{theorem}
\begin{proof}
	Firstly, we notice that if $x\in \langle \mathcal{A}_c\rangle$, then $x$ can be expressed as a polynomial over core elements of $\mathcal{A}$. We thus have that $x=\phi(y_0,...,y_m)$, where $y_0,...,y_m\in \mathcal{A}_c$. It thus suffices to show by induction on the complexity of $\phi$ that $\phi= \bigvee_{i\leq n}a_i $,	 for some $a_0,...a_n\in \mathcal{A}_c$.
	\begin{itemize}
		\item If $\phi=a$, then obviously $\phi= \bigvee \{a \}$.
		%\item If $\phi=\bot$, then $\phi = \bigvee\{0 \}$ and $0\in \mathsf{core}(\mathcal{A})$ by definition.
	\end{itemize}
	\noindent By the induction hypothesis we have $\psi=\bigvee_{i\leq n} a_i$ and $\chi=\bigvee_{j\leq m}b_j$, then:
	\begin{itemize}			
		\item If $\phi=\psi\land \chi$, then:
		\begin{align*}
		\phi=\bigvee_{i\leq n} a_i \land \bigvee_{j\leq m}b_j = \bigvee_{i\leq n}\bigvee_{j\leq m} (a_i\land b_j).
		\end{align*}
		\noindent And since $ \mathcal{A}_c $ is closed under conjunction, $a_i\land b_j\in  \mathcal{A}_c $ for all $i\leq n, j\leq m$.
		
		\item If $\phi=\psi\lor \chi$, then:
		\begin{align*}
		\phi=\bigvee_{i\leq n} a_i \lor \bigvee_{j\leq m}b_j =\bigvee \{a_0,...,a_n,b_0,...,b_m\}.
		\end{align*}
		
		\item  If $\phi= \psi \rightarrow \chi$, then:
		\begin{align*}
		\phi &= \bigvee_{i\leq n}a_i\rightarrow \bigvee_{j\leq m}  b_j \\
		&= \bigwedge_{i\leq n}\Big(a_i \rightarrow \bigvee_{j\leq m}  b_j  \Big) & (\text{by } \langle \mathcal{A}_c \rangle  \text{ being a Heyting algebra)}\\
		&= \bigwedge_{i\leq n}\Big(\bigvee_{j\leq m} ( a_i \rightarrow   b_j)  \Big) & (\text{by \textit{Split}}) \\
		& = \bigvee_{f:[n]\rightarrow[m]} \Big(\bigwedge_{i\leq n}(a_i \rightarrow b_{f(j)})  \Big).
		\end{align*}
		\noindent  Where $ f:[n]\rightarrow[m] $ means that $f\in (m+1)^{n+1}$. Now, since $\mathcal{A}_c$ agrees with $\langle \mathcal{A}_c \rangle$ with respect to the reduct $\{\bot, \land,\rightarrow \}$, it  follows that  $\bigwedge_{i\leq n}(a_i \rightarrow b_{f(j)})  \in \mathcal{A}_c $.

	\end{itemize}
	
	\noindent We have  obtained that every $x\in\langle \mathcal{A}_c\rangle$ has a disjunctive representation $x= \bigvee_{i\leq n}a_i$ with $a_i\in \mathcal{A}_c$ for all $i\leq n$. Let $A= \{a_0, ..., a_n \}$, then  to obtain a non-redundant representation of $x$ it suffices to take the set $I=\{ m\leq n : a_m \text{ is maximal in } A \} $. Then clearly $x= \bigvee_{i\leq n}a_i= \bigvee_{i\in I}a_i$ and by construction  $a_i\nleq a_j$ for $i, j\in I$ such that $i\neq j$.
	
	Finally, if $\mathcal{A}$ is a dependence algebra, it suffices to supplement the previous reasoning with the following case:	
	\begin{itemize}
		\item If $\phi=\psi\otimes \chi$, then by the \textit{Dist} axiom we have:
		\[  \phi=\bigvee_{i\leq n} a_i \otimes \bigvee_{j\leq m}b_j=\bigvee_{i\leq n}\bigvee_{j\leq m} (a_i\otimes b_j). \]
	\end{itemize}
	
	\noindent Hence, since $ \mathcal{A}_c $ is closed under the tensor disjunction $\otimes$, it follows that $a_i\otimes  b_j\in \mathcal{A}_c $ for every $i\leq n, j\leq m$. This completes the proof of our claim.
\end{proof}

\noindent If $\mathcal{A}$ is a core-generated inquisitive or dependence algebra, then the following result follows immediately.

\begin{corollary}
	Let $\mathcal{A}$ be a core-generated inquisitive (or dependence) algebra, then for all $x\in \mathcal{A}$  there are pairwise incomparable elements $a_0,...,a_n\in \mathcal{A}_c$ such that $x=\bigvee_{i\leq n} a_i$. 
\end{corollary}

Core-generated algebras play a special role in the theory of inquisitive and dependence algebras, as they are algebras for which the inquisitive and dependence axioms hold for all elements of the underlying universe. In particular, core-generated structures also have the important role of generators of the class of all inquisitive and dependence algebras. This is made precise by the following proposition.

\begin{proposition}\label{Super}
	Every inquisitive (and dependence) algebra is a core-superstructure of a core-generated inquisitive algebra.
\end{proposition}
\begin{proof}
	Let $\mathcal{A}$ be an arbitrary inquisitive (dependence) algebra and consider the core-generated algebra $\langle \mathcal{A}_c\rangle$. Clearly $\langle \mathcal{A}_c\rangle_c =  \mathcal{A}_c$. Hence $ \langle \mathcal{A}_c\rangle \preceq \mathcal{A}$, which proves our claim.
\end{proof}

\noindent This gives us a first characterisation of the classes of inquisitive and dependence algebras. For any class $\mathsf{C}$ of $\inqI$-algebras or $\inqI^\otimes$-algebras, we let:
$$\mathsf{C}^\uparrow := \{ \mathcal{B} : \mathcal{A}\preceq \mathcal{B} \text{ and } \mathcal{A}_c=\mathcal{B}_c  \text{ for some } \mathcal{A}\in \mathsf{C} \}.$$
\noindent Moreover, we have the following Proposition.
\begin{proposition}\label{birkhofflike1}
	(i) If $\phi\in\langInt$ and $\mathcal{A}$ is an inquisitive algebra, then  $\mathcal{A}\nvDash^c\phi$ entails $\langle \mathcal{A}_c \rangle \nvDash^c\phi$. (ii) If $\phi\in\langInqI$ and $\mathcal{A}$ is a dependence algebra, then  $\mathcal{A}\nvDash^c\phi$ entails $\langle \mathcal{A}_c \rangle \nvDash^c\phi$.
\end{proposition}{}
\begin{proof}
	This is an immediate consequence of Proposition \ref{closure1}(ii).
\end{proof}

\noindent  It follows by the previous propositions that $\mathsf{InqAlg\Lambda}= (\mathsf{InqAlg\Lambda_{CG}})^\uparrow $ and $ \mathsf{InqAlg\Lambda^\otimes}= (\mathsf{InqAlg\Lambda^\otimes_{CG}})^\uparrow$. We thus obtain that core-generated $\inqint$-algebras generate the class of $\inqint$-algebras under the core superalgebra operator defined above. Similarly, core-generated $\inqint^\otimes$-algebras generate the class of $\inqint^\otimes$-algebras. This result is similar to what was obtained in \cite{Quadrellaro.2019B} for so-called $\dna$-varieties, though there core superalgebras of core-generated inquisitive algebras were assumed to be always Heyting algebras.

As we have seen, if an inquisitive or dependence algebra $\mathcal{A}$ is core-generated, then we can talk about arbitrary elements of $\mathcal{A}$ and describe its full structures. With a similar motivation, we introduce well-connected inquisitive and dependence algebras, as these are algebras for which we can give a characterisation of their core elements. 

A Heyting algebra $\mathcal{H}$ is \textit{well-connected} if for all $x,y\in \mathcal{H}$, if $x\lor y =1$ then $x=1$ or $y=1$. We say that an inquisitive (or dependence) algebra $\mathcal{A}$ is \textit{well-connected} if $\langle \mathcal{A}_c\rangle$ is well-connected.  We say that $x\in \mathcal{H}$ is  \textit{join-irreducible} if, for all $a,b \in \mathcal{H}$, $x= a\lor b$ entails $x= a$ or $x= b$.  Notice that in every distributive lattice, and thus in every Heyting algebra, join-irreducible elements coincide with the \textit{join-prime} elements, i.e. those elements $x\in \mathcal{H}$ such that for all $a,b \in \mathcal{H}$, if $x\leq a\lor b$ then $x\leq a$ or $x\leq b$. For any inquisitive (or dependence) algebra $\mathcal{A}$, we say that an element $x\in \langle \mathcal{A}_c \rangle$ is join-irreducible (join-prime) if $x$ is join-irreducible (join-prime) in $\langle \mathcal{A}_c \rangle$.  We denote by $ \mathcal{A}_{ji} $ the subset of join-irreducible members of $\langle \mathcal{A}_c \rangle$.

The following proposition provides a characterisation of core elements of well-connected, inquisitive and dependence algebras.

\begin{proposition}\label{wellconnectedsi}
	Let $\mathcal{A}$ be a well-connected, inquisitive or dependence algebra, then  $\mathcal{A}_c=\mathcal{A}_{ji}$
\end{proposition}
\begin{proof}
	$(\subseteq)$ Let $a\in\mathcal{A}_c $ and suppose that for some $x,y\in \langle \mathcal{A}_c \rangle$ we have that $a\leq x\lor y$. It follows that $a \rightarrow x\lor y =1$ and therefore, by \textit{Split}, $(a \rightarrow x)\lor (a \rightarrow y)=1$. Since $\mathcal{A}$ is well-connected, either  $a \rightarrow x=1$ or $a \rightarrow y=1$, which entails $a\leq x$ or $a\leq y$. Hence $a\in  \mathcal{A}_{ji}$. 		$(\supseteq)$   Suppose $x\in \mathcal{A}_{ji}$, then by definition $x\in \langle \mathcal{A}_{c}\rangle$. Hence by Theorem \ref{normalformalgebra} we have that $x=\bigvee_{i\leq n}a_i $ with $a_i\in\mathcal{A}_{c} $ for all $i\leq n$. Since $x$ is join-irreducible there is some $i\leq n$ for which $x= a_i$, which yields $x\in \mathcal{A}_c$.
\end{proof}

\begin{corollary}\label{joingenerated}
	Let $\mathcal{A}$ be a well-connected, core-generated inquisitive or dependence algebra, then $\mathcal{A}$ is generated by its subset of join-irreducible elements.
\end{corollary}

We have thus obtained an important characterisation of core elements of well-connected inquisitive and dependence algebras. This fact will be important in our duality results of Section \ref{four}. Here, we can immediately prove an important result, showing that the validity of formulas is preserved under homomorphic images of well-connected, core-generated inquisitive and dependence algebras.

\begin{proposition}\label{closure3}
	Let $\mathcal{A}$ be a core-generated, well-connected inquisitive or dependence algebra such that $\mathcal{A}\vDash^c \phi$ and $h:\mathcal{A}\twoheadrightarrow\mathcal{B}$, then $\mathcal{B}\vDash^c \phi$.
\end{proposition}
\begin{proof} 
	The proof is the same for inquisitive and dependence algebras.  We firstly show that $\mathcal{B}_c\subseteq h[\mathcal{A}_c]$. Suppose $y\in \mathcal{B}_c$ and let $h(x)=y$, then since $\mathcal{A}$ is core-generated we have that $x=\bigvee_{i\in I}a_i$ and $a_i\in \mathcal{A}_c$ for all $i\leq n$. It follows that  $y=\bigvee_{i\leq n} h(a_i)$ and since $y\in \mathcal{B}_c$ we have that $y$ is join-irreducible, which yields $y=h(a_i)$ for some $i\leq n$. %Therefore, for every $y\in \mathcal{B}_c$ there is some $a\in \mathcal{A}_c$ such that $h(a)=y$.
	
	Now, let $\vec{p}:=p_0,\dots,p_n$ and suppose towards contradiction that  $\mathcal{A}\vDash^c \phi(\vec{p})$ and $\mathcal{B}\nvDash^c \phi(\vec{p})$. Then $(\mathcal{B},\mu)\nvDash^c \phi(\vec{p})$ for some core-valuation $\mu$. Let $\nu:\at\to \mathcal{A}_c$ be such that for all $p\in\at$, $\nu(p)\in h^{-1}(\mu(p))\cap \mathcal{A}_c$. Notice that $\mu$ is well-defined by the considerations of the previous paragraph. Since $(\mathcal{A},\nu)\vDash^c \phi(\vec{p})$ then $\llbracket \phi(\vec{p}) \rrbracket^{(\mathcal{A},\nu)} =1_{\mathcal{A}}$ and since $h$ is a $\inqI$-homomorphism:
	\begin{align*}
	1_{\mathcal{B}}=h (1_{\mathcal{A}})=h (\llbracket \phi(\vec{p}) \rrbracket^{(\mathcal{A},\nu)})=\llbracket \phi(\vec{p}) \rrbracket^{(\mathcal{B},\mu)}.
	\end{align*}
	\noindent which entails $(\mathcal{B},\mu)\vDash^c \phi(\vec{p})$, contradicting our assumption.
\end{proof}

\subsection{A Birkhoff-Like Theorem for Inquisitive Algebras}

An important result proven by Birkhoff for varieties of algebras states that an equation is true in a variety if and only if it is true in its subclass of subdirectly irreducible algebras. Here we prove a similar result for a suitable subclass of inquisitive algebras. We shall deal separately in the next section with dependence algebras, as that case involves further complications. 

Firstly, we prove some preliminary results concerning finite inquisitive algebras. In the standard setting of universal algebra \cite[p. 69]{Burris.1981}, we say that an algebra $\mathcal{A}$ is \textit{locally finite} if for every $X\subseteq \mathsf{dom}(\mathcal{A})$ such that $|X|< \omega$ we have $|\langle X\rangle |< \omega$. A class of algebras $\mathsf{C}$ is locally finite if every $\mathcal{A}\in\mathsf{C}$ is locally finite. We recall the following well-known facts. We refer the reader to  \cite{DiegoAntonio1966Slad} and \cite{Grtzer2011} for a proof of the following statements.

\begin{theorem}[Diego, Folklore]
	$\:$
	\begin{itemize}
		\item[(i)] The class of Browerian semilattices is locally finite;
		\item[(ii)] The class of bounded distributive lattices is locally finite.
	\end{itemize}
\end{theorem}

To make sense of this property in our context, we should consider only subsets $X\subseteq \mathsf{core}(\mathcal{A})$. We say that an inquisitive algebra $\mathcal{A}$ is \textit{locally finite} if  for every $X\subseteq \mathsf{core}(\mathcal{A})$ such that $|X|< \omega$ we have $|\langle X\rangle| < \omega$.  A class of $\inqI$-algebras $\mathsf{C}$ is locally finite if every $\mathcal{A}\in\mathsf{C}$ is locally finite.  The following theorem shows that $\inqint$-algebras are locally finite.

\begin{theorem}\label{localfiniteness}
	$\inqint$-algebras are locally finite.
\end{theorem}
\begin{proof}Let $\mathcal{A}$ be a $\inqint$-algebra  and let $X\subseteq \mathcal{A}_c$. We first close $X$ under the Brouwerian semilattice operations $\{  \land, \rightarrow, 0\}$, we obtain the subalgebra $Y=\langle X\rangle_{(\land, \rightarrow, 0)}$ of $\mathcal{A}_c$. It follows by Diego's Theorem that $Y$ is finite. Secondly, we close $Y$ under meet and join, and we obtain the set $Z= \langle  Y \rangle_{(\land,\lor)}$. Since $( \mathcal{A},\land,\lor)$ is a bounded distributive lattice, it follows from the previous theorem that $Z$ is finite as well.
	
	To obtain an inquisitive algebra, we need to supplement $Z$ with a Heyting implication.  Notice that, for all $x,y\in Z$, we have that $x=\phi(\vec{a})$ and $y=\psi(\vec{b})$ for some tuples $\vec{a}=a_0,\dots,a_n\in X$ and $\vec{b}=b_0,\dots,b_m\in X$. By reasoning as in Theorem \ref{normalformalgebra}, we can put $\phi$ and $\psi$ in disjunctive form and show that for all $x,y\in Z$, $ x=\bigvee_{i\leq n}c_i$ and $y= \bigvee_{j\leq m}  d_j  $ for $c_i,d_j\in Y $. We define:	 
	\begin{align*}
	x\dot{\to} y &:= \bigvee_{f:[n]\rightarrow[m]} \Big[\bigwedge_{i\leq n}(c_i \rightarrow d_{f(j)})  \Big].
	\end{align*}
	\noindent By proceeding again as in the proof of the Normal Form Theorem, we see that for all $x,y\in Z$, $x\dot{\to}y=x\to y$, hence $(Z, \land,\lor,\rightarrow,0)$ is a well-defined Heyting subalgebra of $\langle \mathcal{A}_c\rangle.$
	
	Since $X\subseteq \mathcal{A}_c$ and $\mathcal{A}_c$ is a Brouwerian semilattice, it follows that $Y\preceq \mathcal{A}_c$ and that  $\mathcal{B}=(Z,Y, \land,\lor,\dot{\to},0)$ is an inquisitive subalgebra of $\mathcal{A}$.  Since by construction $Z= \langle X \rangle$, it follows that $\mathcal{B}$ is the smallest subalgebra of $\mathcal{A}$  containing $X$, i.e. $\mathcal{B}$ is the subalgebra of $\mathcal{A}$ generated by $X$. By the closure of $\inqint$-algebras under subalgebras we then have that $\mathcal{B}$ is a $\inqint$-algebra. Finally, by what we have argued  above, $|Z|<\omega$, hence $\mathcal{A}$ is locally finite and thus $\inqint$-algebras are locally finite.
\end{proof}

We then obtain the algebraic version of the finite model property.

\begin{theorem}[Finite Model Property]\label{finiteness}
	Suppose $\mathcal{A}$ is a $\inqint$-algebra and $\mathcal{A} \nvDash^c \phi$, then there is a finite $\inqint$-algebra $ \mathcal{B} $ such that $\mathcal{B} \nvDash^c \phi$.
\end{theorem}
\begin{proof}
	Let $\mathcal{A}$ be a $\inqint$-algebra such that $\mathcal{A} \nvDash^c \phi$, then for some  core valuation $\mu$ and $\mathcal{M}=(\mathcal{A},\mu)$ we have $\mathcal{M}\nvDash^c \phi $.	Let $At(\phi)=\{p_0,\dots,p_n  \}$, then by Theorem \ref{localfiniteness} above  $\langle At(\phi) \rangle$ is finite. The structure  $\mathcal{B}=(\langle At(\phi) \rangle, \langle At(\phi) \rangle_{(\land,\to, 0)}, \land,\lor,\dot{\to},0)$ is then a finite inquisitive subalgebra of $\mathcal{A}$. By  Corollary \ref{closure2}, $\mathcal{B}$ is a $\inqint$-algebra and since $(\mathcal{B},\mu)\nvDash^c \phi $ it follows that $\mathcal{B} \nvDash^c \phi$.
\end{proof}

We use the former results to prove a version of Birkhoff's Theorem for inquisitive algebras. In the standard setting, Birkhoff's result \cite[Thm. 9.6]{Burris.1981} says that every algebra in a variety is a subdirect product of subdirectly irreducible algebras. As a consequence, this means that an equation holds in a variety of algebras if and only if it holds in its subclass of subdirectly irreducible elements.  Here we prove a similar result for the class of $\inqint$-algebras: the next theorem specifies a class of representatives which witness the truth and falsity of formulas in inquisitive algebras. 

Recall that a Heyting algebra $\mathcal{H}$ is subdirectly irreducible if and only if it has a second greatest element. Also, if $\mathcal{H}$ is finite, then $\mathcal{H}$ is subdirectly irreducible if and only if $\mathcal{H}$ is well-connected. Finally, we also recall the following fact, originally due to Wronski.

\begin{proposition}[Wronski]\label{wronski}
	Let $\mathcal{H}$ be a Heyting algebra and $1\neq x\in\mathcal{H}$. Then there is a surjective Heyting homomorphism $h: \mathcal{H} \twoheadrightarrow  \mathcal{B} $ such that $\mathcal{B}$ is a subdirectly irreducible Heyting algebra and $h(x)= s_B$, where $s_B$ is the second greatest element in $\mathcal{B}$.
\end{proposition}

\noindent We refer the reader to \cite{9783030120955} and \cite{Wronski} for the proofs of the previous claims.

\begin{theorem}\label{birkhoff2}
	Suppose $\mathcal{A}$ is a $\inqint$-algebra and $\mathcal{A} \nvDash^c \phi$, then there is a finite, core-generated, well-connected $\inqint$-algebra $\mathcal{B}$ such that $\mathcal{B} \nvDash^c \phi$.
\end{theorem}
\begin{proof}
	Suppose $\mathcal{A}$ is a $\inqint$-algebra such that $\mathcal{A} \nvDash^c \phi$. By Theorem \ref{finiteness} above there is a finite $\inqint$-algebra $\mathcal{D}$ such that  $\mathcal{D} \nvDash^c \phi$ and, by Proposition \ref{birkhofflike1}, it follows that $\langle \mathcal{D}_c \rangle \nvDash^c \phi$. Hence there is  some core-valuation $\mu$ such that $\mathcal{M} \nvDash^c \phi$, where $\mathcal{M}=(\langle \mathcal{D}_c\rangle, \mu)$.   By Proposition \ref{wronski}, there is a surjective homomorphism $h:\langle \mathcal{D}_c \rangle \twoheadrightarrow  \mathcal{C} $ such that $\mathcal{C}$ is a subdirectly irreducible Heyting algebra and $h(\phi^{\mathcal{M}})= s_C$, where $s_C$ is the second greatest element in $\mathcal{C}$. Now let $\mathcal{B}=(h[\langle \mathcal{D}_c\rangle], h[\mathcal{D}_c], \land,\lor, \rightarrow, 0   )$ and let $\nu:\at\to\mathcal{B}_c$ be the core assignment $\nu := h\circ \mu$. Then from the fact that $\langle \mathcal{D}_c \rangle \nvDash^c \phi$, we obtain  $(\mathcal{B},\nu)\nvDash^c\phi$ and therefore $\mathcal{B} \nvDash^c \phi$. It thus suffices to verify that $\mathcal{B}$ is a finite, core-generated, well-connected $\inqint$-algebra.
	
	Clearly $\mathcal{B}$ is a Heyting algebra and by the fact that $\mathcal{B}_c$ is homomorphic image of $\mathcal{D}_c$ under $h$, it follows that $\mathcal{B}_c=h[\mathcal{A}_c]$ is a Brouwerian semilattice.  Since $\mathcal{D}$ is finite and $\mathcal{B}=h[\langle \mathcal{D}_c\rangle]$, we have that $ \mathcal{B} $ is finite.  Moreover, since $\mathcal{B}=h[\langle \mathcal{D}_c\rangle]$ and $\mathcal{B}_c=h[ \mathcal{D}_c]$,  $\mathcal{B}$ is clearly generated by $\mathcal{B}_c$.
	
	%	Since $h$ is a Heyting  homomorphism we have that, if $a,b\in \mathcal{B}_c$, then $a=h(a')$ and $b=h(b')$, which entails:
	%	$$a\land b=h(a')\land h(b')=h(a'\land b').$$
	%	\noindent Hence, $a\land b\in \mathcal{B}_c$. By a similar argument, it follows that $ \mathcal{B}_c $ is closed under $\rightarrow$ and contains $1,0$. Moreover, since $h$ preserves all operations in $\{1,0,\rightarrow, \land  \}$, and since $ \mathcal{A}_c{\upharpoonright} \{1,0,\rightarrow, \land  \}$ is a Brouwerian semilattice, it follows that $\mathcal{B}_c$ is also a Brouwerian semilattice.

	%	We now claim that $\mathcal{B}$ is generated by $\mathcal{B}_c$. Let $x\in \mathcal{B} $, then since $h$ is surjective, $x=h(x')$ for some $x'\in\langle \mathcal{D}_c\rangle$, hence $x'= \bigvee_{i\leq n}a_i$ and:
	%	$$x= h(x') = h(\bigvee_{i\leq n}a_i) = \bigvee_{i\leq n}h(a_i) . $$
	%	\noindent Then, for all $i\leq n$, we have that $h(a_i)\in \mathcal{B}_c$, hence $x\in \langle \mathcal{B}_c \rangle $, which shows $\mathcal{B} \subseteq \langle \mathcal{B}_c \rangle$. Since it is obvious that $\langle \mathcal{B}_c \rangle \subseteq  \mathcal{B}$, we obtain $\mathcal{B} = \langle \mathcal{B}_c \rangle$.
	
	We next claim that $\mathcal{B}$ is well-connected. By construction, $ \mathcal{B}{\upharpoonright}\{ \rightarrow, \land, \lor, 1,0 \}=\mathcal{C}$ is a subdirectly irreducible Heyting algebra and, since it is finite, it is also well-connected.

	To see that $\mathcal{B}$ validates the \textit{Split} axiom, let $a \in \mathcal{B}_c$, $x,y\in \langle \mathcal{B}_c \rangle$. Since $h[\langle \mathcal{D}_c\rangle]=\mathcal{B}$ and $h[ \mathcal{D}_c]=\mathcal{B}_c$, there are $a' \in \mathcal{D}_c$, $x',y'\in \langle \mathcal{D}_c\rangle$ such that $h(a')=a$, $h(x')=x$ and $h(y')=y$. Since $\langle \mathcal{D}_c\rangle$ is an inquisitive algebra  we then obtain:
	\begin{align*}
	\; & \; a' \rightarrow (x' \lor y') = (a' \rightarrow x') \lor (a' \to y')  \\
	\Longrightarrow \;  &\; h[a' \rightarrow (x' \lor y')]   = h[(a' \rightarrow x') \lor (a' \to y')]  \\
	\Longrightarrow \; &\; h(a') \rightarrow (h(x') \lor h(y'))  = ( h(a') \to h(x')) \lor (h(a') \rightarrow  h(y')) \\
	\Longrightarrow \; & \; a \rightarrow (x \lor y)  =(a \rightarrow x) \lor (a \to y).
	\end{align*}
	\noindent hence $\mathcal{B}$ is a $\inqI$-algebra.
	
	Finally, we can show in exactly the same way that $\Lambda\subseteq\langIntsta$, $ \mathcal{D}\vDash^c\Lambda$ entail  $\mathcal{B}\vDash^c\Lambda$, meaning that $\mathcal{B}$ is a $\inqint$-algebra. It follows that $\mathcal{B}$ is a finite, core-generated, well-connected $\inqint$-algebra and that $\mathcal{B}\nvDash \phi$, which  proves our theorem.	
\end{proof}

The previous theorem shows why finite, core-generated, well-connected inquisitive algebras are of special importance in the theory of inquisitive algebras. In fact, by the former result, they witness the validity of inquisitive formulas. We write  $ \mathsf{InqAlg_{FCGW}} $ for the category (and the class) of finite, core-generated, well-connected inquisitive algebras with $\inqI$-homomorphisms.

\subsection{A Birkhoff-Like Theorem for Dependence Algebras}\label{seven}

We prove in this section a theorem analogous to Theorem \ref{birkhoff2} for the class of dependence algebras. In this case we will show only a weaker version of our former result: We prove  that finite, core generated, well-connected dependence algebras witness the validity of formulas true in the class of well-connected dependence algebras. We conclude this section by proving a stronger version of this result for the case of locally tabular intermediate dependence logics.

Firstly, we prove the following lemma, which shows that every surjective map between core-generated algebras which preserves the Heyting operations also preserves the tensor disjunction.

\begin{lemma}\label{homlemma}
	Suppose $f:\mathcal{A}\rightarrow \mathcal{B}$ is a surjective Heyting homomorphism between two core-generated  $\inqI^\otimes$-algebras $\mathcal{A}$ and $\mathcal{B}$, then it follows that for all $x,y\in \mathcal{A}$, $f(x\otimes y)=f(x)\otimes f(y)$.
\end{lemma}
\begin{proof}
	
	Since  $f$ preserves the operations in $\{0,\land, \rightarrow  \}$, $f{\upharpoonright} \mathcal{A}_c$ is a Brouwerian semilattice homomorphism. By Lemma 2.4 in \cite{kohler1981brouwerian} (see also \cite[Lemma 2]{Bezhanishvili2015}) surjective Brouwerian semilattices homomorphisms preserve existing join, hence for all $a,b\in \mathcal{A}_c$, $f(a\otimes b)=f(a)\otimes f(b)$.
	
	We extend this result to arbitrary elements of $\mathcal{A}$. Let  $x,y\in \mathcal{A}$, then, by the Normal Form Result $x=\bigvee_{i\leq n}a_i$ and $y=\bigvee_{j\leq m}b_j$, where $a_i,b_j\in \mathcal{A}_c$, for all $i\leq n$ and $j\leq m$. We obtain:
	\begin{align*}
	f(x) \otimes f(y) &= f\Big(\bigvee_{i\leq n}a_i\Big) \otimes f\Big(\bigvee_{j\leq m}b_j\Big) &\\
	& = \bigvee_{i\leq n}f(a_i) \otimes \bigvee_{j\leq m}f(b_j) & \text{ (by \textit{f} homomorphism)}\\
	& = \bigvee \{ f(a_i)\otimes f(b_j) : i\leq n, j\leq m  \text{ and } a_i,b_j\in \mathcal{A}_c   \} &\text{ (by \textit{Dist})} \\ 
	& = \Big[\bigvee \{ f(a_i\otimes b_j) : i\leq n, j\leq m  \text{ and } a_i,b_j\in \mathcal{A}_c  \}\Big] & \text{ (by $ a_i,b_j\in\mathcal{A}_c$)}\\
	& = f\Big[\bigvee \{ a_i\otimes b_j : i\leq n, j\leq m  \text{ and } a_i,b_j\in \mathcal{A}_c  \}\Big] & \text{ (by \textit{f} homomorphism)}\\
	& = f\Big(\bigvee_{i\leq n}a_i \otimes \bigvee_{j\leq m}b_j\Big) & \text{ (by \textit{Dist})}  \\
	& = f(x\otimes y). &  
	\end{align*}
	\noindent Hence for all $x,y\in \mathcal{A}$, $f(x\otimes y)=f(x)\otimes f(y)$, proving our claim.
\end{proof}

%\noindent It thus follows immediately that $\inqI$-homomorphisms between core-generated dependence algebras are also $\inqI^\otimes$-homomorphisms.

\begin{corollary}\label{homlemma2}
	Suppose $h:\mathcal{A}\rightarrow \mathcal{B}$ is a surjective $\inqI$-homomorphism and that $\mathcal{A}$ and $\mathcal{B}$ are core-generated dependence algebras, then $h$ is a $\inqI^\otimes$-homomorphism.
\end{corollary}

While Brouwerian semilattices are locally finite, Heyting algebras are not, hence the proof we gave to Theorem \ref{localfiniteness} cannot be replicated in the setting of dependence algebras. We can then prove only a limited version of the Finite Model Property: we show that if a formula is falsified by a well-connected dependence algebra, then it is falsified by a finite dependence algebra.

\begin{theorem}\label{finiteness3}
	Suppose $\mathcal{A}$ is a well-connected $\inqI^\otimes$-algebra and $\mathcal{A} \nvDash^c \phi$, then there is a finite $\inqI^\otimes$-algebra $\mathcal{B}$ such that $\mathcal{B} \nvDash^c \phi$.
\end{theorem}
\begin{proof}
	We adapt to our context the strategy of the proof of the finite model property for the variety of Heyting algebras. Given a well-connected dependence algebra $\mathcal{A}$ such that $\mathcal{A} \nvDash^c \phi$, the main idea of this proof is to generate a finite distributive lattice $\mathcal{B}\preceq \mathcal{A}$ such that the tensor is well-defined over $\mathcal{B}_c$ and $\mathcal{B} \nvDash^c \phi$. Then, using the fact that $\mathcal{A}$ is well-connected, we define a ``fake'' heyting implication and we lift the tensor join to the whole of $\mathcal{B}$, so that we turn $\mathcal{B}$ into a suitable dependence algebra. 
	
	Suppose that $\mathcal{A} \nvDash^c \phi$ and $\mathcal{A}$ is well-connected. Let $\mu$ be a core valuation  such that $(\mathcal{A},\mu)\nvDash^c \phi $ and let $\mathcal{M}=(\mathcal{A},\mu)$.	By the Normal Form Theorem \ref{disjnf} for $\inqI^\otimes$, we can assume without loss of generality that $\phi= \bigvee_{i\leq n} \alpha_i$, where $\alpha_i$ is a standard formula for all $i\leq n$. We denote by $Sub(\phi)$ the set of subformulas of $\phi$. 
	
	Consider the set $ X= \{ \tau^\mathcal{M} : \tau \in Sub(\alpha) \}$. Since every $\alpha_i$ is standard, we clearly have that $X\subseteq \mathcal{A}_c$. By closing $X$ under all operations in $\{ \land,\otimes,0 \}$ we obtain $Y= \langle X \rangle_{  (\land,\otimes,0 ) }$. Since $\mathcal{A}_c  {\upharpoonright}\{ \land,\otimes,0 \}$ is a bounded distributive lattice it follows that $Y$ is a finite $\{ \land,\otimes,0 \}$--subalgebra of $\mathcal{A}_c$. We then close $Y$ under the operations $\{ \land,\lor \}$ and  obtain  $Z= \langle  Y \rangle_{(\land,\lor,0)}$. Since $Y$ is finite and $\mathcal{A}  {\upharpoonright}\{ \land,\lor,0 \}$ is a bounded distributive lattice it follows that $Z$ is finite. 
	
	\begin{claim}
		Let $x\in Z$, then $x=\bigvee_{i\leq n} b_i$ and $b_i\in Y$ for all $i\leq n$.
	\end{claim}
	\begin{proof}\renewcommand{\qedsymbol}{$\rfloor$}
		The claim follows by induction over $\{ \land, \lor,0  \}$ as in the Normal Form Theorem \ref{normalformalgebra}. 
	\end{proof}
	
	Now, in order to obtain an $\inqI^\otimes$-algebra $\mathcal{B}$ such that $\mathsf{dom}(\mathcal{B})=Z$ and $\mathsf{core}(\mathcal{B})=Y$, it suffices to extend the tensor operator to arbitrary elements of $Z$ and define a Heyting implication over such structure. To this end, we define:
	$$ x \otimes y := \bigvee \{ a\otimes b : a\leq x, b\leq y \text{ and }  a,b\in Y   \}.  $$
	\noindent  For all $x,y\in \mathcal{B}$, we have by the previous claim that  $ x=\bigvee_{i\leq n}a_i$ and $y= \bigvee_{j\leq m}  \beta_j  $ with $a_i,\beta_j\in Y $, hence by the finiteness of $Y$, the former definition is equivalent to $x \otimes y=\bigvee_{i\leq n}\bigvee_{j\leq m} (a_i\otimes b_j ) $ with $a_i,b_j\in Y$. Similarly, we proceed by defining a new \enquote{fake} implication $\dot{\to}$ as follows:
	$$x\dot{\to}y :=  \bigvee   \{ c\in Z : c\land x \leq y    \}.$$ 
	\noindent Notice that, since $Z$ is finite, $ \dot{\to}  $ is well-defined. We now prove the following claims.
	
	\begin{claim}
		The structure $\mathcal{B}= ( Z, Y\cap \mathcal{A}_c, \otimes, \lor, \land, \dot{\to}, 0 )$ is an $\inqI^\otimes$-algebra.		
	\end{claim}
	\begin{proof}\renewcommand{\qedsymbol}{$\rfloor$}
		Firstly, since $\mathsf{dom}(\mathcal{B})=Z$ and $(Z, \land,\lor, 0 )$ is a bounded distributive lattice, we have that $(\mathsf{dom}(\mathcal{B}), \land, \lor, \dot{\to}, 0)$ is a bounded distributive lattice together with a well-defined Heyting implication, hence it is a Heyting algebra. Similarly, $\mathsf{core}(\mathcal{B})=Y\subseteq \mathcal{A}_c$, hence since $\mathcal{A}_c$ is a bounded distributive lattice and $Y$ is closed under $\{ \land, \otimes \}$, it follows that $\mathsf{core}(\mathcal{B})$ is also a bounded distributive lattice. We have by the previous claim:
		\begin{align*}
		x\dot{\to}y & :=  \bigvee   \{ c\in \mathcal{B} : c\land x \leq y    \} \\
		&= \bigvee  \Big \{ \bigvee_{i\leq n} (a_i) : \bigvee_{i\leq n} (a_i)\land x \leq y  \text{ and } a_i \in \mathsf{core}(\mathcal{B}) \text{ for all } i\leq n \Big \} \\
		&= \bigvee   \{ a_i \in \mathsf{core}(\mathcal{B})  : a_i \land x \leq y \text{ for all } i\leq n  \}.\textit{}
		\end{align*}
		\noindent Hence $\dot{\to}$ is a well-defined Heyting implication and $\mathsf{core}(\mathcal{B})$ is a Heyting algebra.
		
		It then remains to check the axioms \textit{Dist}, \textit{Split} and \textit{Mon}. We check \textit{Dist} and \textit{Split} only as the case for \textit{Mon} is analogous. 
		
		(\textit{Dist}). Let $x,y,z\in \mathcal{B}$, then there are $a_i,b_j,c_k\in Y$ for all $i\leq n, j\leq m, k\leq l$ such that $ x=\bigvee_{i\leq n}a_i$,  $y= \bigvee_{j\leq m}  b_j  $ and $ z=\bigvee_{k\leq l}c_k$. By the definition of the tensor, we then have: 
		\begin{align*} 
		x \otimes (y\lor z)  & =  \Big(\bigvee_{i\leq n}a_i )\otimes ( \bigvee_{j\leq m}  b_j \lor  \bigvee_{k\leq l}c_k \Big)  \\
		& =\Big(\bigvee_{i\leq n}a_i)\otimes ( \bigvee_{i\leq n}\bigvee_{j\leq m} (b_j\lor c_k ) \Big)  \\
		&   =\Big(\bigvee_{i\leq n} \bigvee_{i\leq n}\bigvee_{j\leq m} (a_i \otimes (b_j\lor c_k )) \Big) .
		\end{align*}	
		\noindent And since $a_i , b_j, c_k \in Y $, it follows that:
		$a_i \otimes^{\mathcal{B}} (b_j\lor c_k ) = a_i \otimes^{\mathcal{A}} (b_j\lor c_k )  $. Thus, since the \textit{Dist} axiom holds in $\mathcal{A}$ we then obtain the following:
		\begin{align*} 
		\bigvee_{i\leq n} \bigvee_{i\leq n}\bigvee_{j\leq m} (a_i \otimes (b_j\lor c_k ))  &=  \bigvee_{i\leq n} \bigvee_{j\leq m}\bigvee_{k\leq l} ((a_i \otimes b_j)\lor (a_i \otimes c_k) ) \\
		& =  \Big(\bigvee_{i\leq n}a_i \otimes  \bigvee_{j\leq m}  b_j\Big) \lor \Big(\bigvee_{i\leq n} a_i \otimes \bigvee_{k\leq l}c_k\Big)   \\
		& = (x \otimes y)\lor   (x \otimes z).
		\end{align*}		
		
		(\textit{Split}). Let $a \in \mathcal{B}_c$ and $y,z\in \mathcal{B}$, then by our previous claim we have that $y= \bigvee_{j\leq m}  b_j  $ and $ z=\bigvee_{k\leq l}c_k$ such that $b_j,c_k\in Y$ for all $j\leq m, k\leq l$. By   $\mathcal{A}$ well-connected, we have:
		\begin{align*} 
		a \dot{\to} (  y \lor z   )
		& = \bigvee   \{ c\in \mathcal{B} : c\land a \leq (y\lor z)    \} \\  
		&= \bigvee   \{ c \in \mathsf{core}(\mathcal{B})  : c\land a \leq \Big(\bigvee_{j\leq m}  b_j\lor \bigvee_{k\leq l}c_k\Big)   \} \\
		& =  \bigvee   \{ c \in \mathsf{core}(\mathcal{B})  : c \land a 	\leq \bigvee_{j\leq m}  b_j \} \lor  \bigvee\{ c \in \mathsf{core}(\mathcal{B})  : c \land a  				\leq \bigvee_{k\leq l}c_k \} \\ 
		& =  (a \dot{\to}   y) \lor (a \dot{\to}z  ).
		\end{align*}
		It follows that $\mathcal{B}$ is an $\inqI^\otimes$-algebra.
	\end{proof}	
	
	Recall that $\mu$ is  a core valuation  such that $(\mathcal{A},\mu)\nvDash^c \phi $ and that by the Normal Form Theorem \ref{disjnf} we assume without loss of generality that  $\phi= \bigvee_{i\leq n} \alpha_i$ with each $\alpha_i$ being standard.
	\begin{claim}
		Let $\nu:\at \to \mathcal{B}$ be a core-valuation such that $\nu{\upharpoonright}\{p_0,\dots,p_n  \}=\mu{\upharpoonright}\{p_0,\dots,p_n  \} $, then  $(\mathcal{B},\nu)\nvDash^c \phi$.
	\end{claim}
	\begin{proof}\renewcommand{\qedsymbol}{$\rfloor$}
		We first prove by induction that for any standard formula $\beta(p_0,\dots,p_n) $, $\llbracket \beta \rrbracket^{(\mathcal{A},\mu)}=\llbracket \beta \rrbracket^{(\mathcal{B},\nu)}$.
		\begin{itemize}
			\item If $\beta=p$ is atomic, then  $\llbracket p \rrbracket^{(\mathcal{B},\nu)}=\llbracket p \rrbracket^{(\mathcal{A},\mu)}$ by definition of $\nu$.

			\item If $\beta=\psi\otimes \chi$, then by the fact that $\psi$ and $\chi$ are standard formulas, we have that $\llbracket \psi \rrbracket^{(\mathcal{A},\mu)}\in Y$ and $\llbracket \chi \rrbracket^{(\mathcal{A},\mu)}\in Y$, hence $\psi \otimes^{\mathcal{B}} \chi=\psi \otimes^{\mathcal{A}} \chi $. We obtain:
			\begin{align*}
			\llbracket \psi\otimes \chi \rrbracket^{(\mathcal{A},\mu)}	=  \llbracket \psi \rrbracket^{(\mathcal{A},\mu)} \otimes^{\mathcal{A}} \llbracket \chi \rrbracket^{(\mathcal{A},\mu)} 
			=  \llbracket \psi \rrbracket^{(\mathcal{B},\nu)} \otimes^{\mathcal{B}} \llbracket \chi \rrbracket^{(\mathcal{B},\nu)} 
			= \llbracket \psi\otimes \chi \rrbracket^{(\mathcal{B},\nu)}.
			\end{align*}
			
			\item If $\beta=\psi\land \chi$,  we proceed analogously.
			
			%			\begin{align*}
			%			\llbracket \psi\land \chi \rrbracket^{(\mathcal{A},\mu)}	&=  \llbracket \psi \rrbracket^{(\mathcal{A},\mu)} \land_{\mathcal{A}} \llbracket \chi \rrbracket^{(\mathcal{A},\mu)} \\
			%			&=  \llbracket \psi \rrbracket^{(\mathcal{B},\nu)} \land_{\mathcal{B}} \llbracket \chi \rrbracket^{(\mathcal{B},\nu)} \\
			%			& = \llbracket \psi\land \chi \rrbracket^{(\mathcal{B},\nu)}.
			%			\end{align*}
			
			\item If $\beta=\psi\to \chi$, we have:
			\begin{align*}
			\llbracket \psi\to \chi \rrbracket^{(\mathcal{B},\nu)}
			&= \bigvee   \{ a_i \in \mathsf{core}(\mathcal{B})  : a_i \land \llbracket \psi \rrbracket^{(\mathcal{B},\nu)} \leq \llbracket  \chi \rrbracket^{(\mathcal{B},\nu)}   \} \\
			&\leq  \bigvee   \{ a_i \in \mathsf{core}(\mathcal{A})  : a_i \land \llbracket \psi \rrbracket^{(\mathcal{A},\mu)} \leq \llbracket  \chi \rrbracket^{(\mathcal{A},\mu)}   \} \\
			& = \llbracket \psi\to \chi \rrbracket^{(\mathcal{A},\mu)}.
			\end{align*}
			\noindent Since by construction $\llbracket \psi\to \chi \rrbracket^{(\mathcal{A},\mu)}\in \mathcal{A}_c \cap X\subseteq \mathcal{B}_c$, it follows that:
			\begin{align*}
			& \; \llbracket \psi\to \chi \rrbracket^{(\mathcal{A},\mu)}\in \{ a_i \in \mathsf{core}(\mathcal{B})  : a_i \land \llbracket \psi \rrbracket^{(\mathcal{A},\mu)} \leq \llbracket  \chi \rrbracket^{(\mathcal{A},\mu)}   \} \\
			\Longrightarrow & \; \llbracket \psi\to \chi \rrbracket^{(\mathcal{A},\mu)} \leq \bigvee   \{ a_i \in \mathsf{core}(\mathcal{B})  : a_i \land \llbracket \psi \rrbracket^{(\mathcal{A},\mu)} \leq \llbracket  \chi \rrbracket^{(\mathcal{A},\mu)}   \}.
			\end{align*}
			\noindent and therefore $\llbracket \psi\to \chi \rrbracket^{(\mathcal{B},\nu)}= \llbracket \psi\to \chi \rrbracket^{(\mathcal{A},\mu)} $.
		\end{itemize}
		
		\noindent Now, we have that for all $i\leq n$,  $\llbracket \alpha_i \rrbracket^{(\mathcal{A},\mu)}=\llbracket \alpha_i \rrbracket^{(\mathcal{B},\nu)}$ and since $\vee^{\mathcal{A}}=\vee^{\mathcal{B}}$, it follows that
		$ \Big \llbracket \bigvee\limits_{i\leq n}\alpha_i \Big\rrbracket^{(\mathcal{A},\mu)}=\Big\llbracket \bigvee\limits_{i\leq n}\alpha_i \Big\rrbracket^{(\mathcal{B},\nu)}. $ Therefore, since $(\mathcal{A},\mu)\nvDash^c \phi$, we have  $(\mathcal{B},\nu)\nvDash^c \phi$.
	\end{proof}
	\noindent Finally, we have obtained a finite $\inqI^\otimes$-algebra $\mathcal{B}$ such that $\mathcal{B}\nvDash^c \phi$, which completes the proof of our theorem.
\end{proof}

The next proposition integrates the previous one and it allows us to obtain, starting from a finite $\inqI^\otimes$-algebra $\mathcal{A}$ such that $\mathcal{A}\nvDash^c\phi$, a finite, core-generated and well-connected $\inqI^\otimes$-algebra  $\mathcal{B}$ such that $\mathcal{B}\nvDash^c\phi$. 

\begin{proposition}\label{birkhoff3}
	Suppose $\mathcal{A}$ is a finite $\inqint^\otimes$-algebra and $\mathcal{A}\nvDash^c\phi$, then there is a finite, core-generated, well-connected $\inqint^\otimes$-algebra $\mathcal{B}$ such that $\mathcal{B} \nvDash^c \phi$.
\end{proposition}
\begin{proof}	
	The proof follows the same strategy of the proof of Theorem \ref{birkhoff2}.  Suppose $\mathcal{A} \nvDash^c \phi$ where $\mathcal{A}$ is a finite dependence algebra. By Proposition \ref{birkhofflike1}, we obtain that $\langle \mathcal{A}_c \rangle \nvDash^c \phi$, hence there is  some core-valuation $\mu$ for which $\mathcal{M} \nvDash^c \phi$, where $\mathcal{M}=(\langle \mathcal{A}_c\rangle, \mu)$.  Now, by Proposition \ref{wronski}, there is a surjective Heyting homomorphism $h:\langle \mathcal{A}_c \rangle \twoheadrightarrow  \mathcal{B} $ such that $\mathcal{B}$ is a subdirectly irreducible Heyting algebra and $h(\phi^{\mathcal{M}})= s_B$, where $s_B$ is the second greatest element in $\mathcal{B}$.
	
	To prove our claim, we specify a tensor operator and a subset of core elements of $\mathcal{B}$, thus obtaining a dependence algebra. Let $\mathcal{B}_c:= h[\mathcal{A}_c]$, then since $h$ is a Heyting  homomorphism in the signature $\{\land, \lor, \to, 0 \}$ we have  that $\mathcal{B}_c$ is a homomorphic image of $ \mathcal{A}_c {\upharpoonright}\{\land, \to, 0 \} $, hence it is a Brouwerian semilattice. Also, notice that $\mathcal{B}_c$ is an ordered structure, where the order $a\leq_{\mathcal{B}_c} b$  is defined by:
	$$a\leq_{\mathcal{B}_c} b \Longleftrightarrow   a\land^{\mathcal{B}_c} b = a .$$
	
	%	we have that, if $a,b\in \mathcal{B}_c$, then $a=h(a')$ and $b=h(b')$, which entails:
	%	$$a\land b=h(a')\land h(b')=h(a'\land b').$$
	%	\noindent Hence, $a\land b\in \mathcal{B}_c$. By a similar argument, it follows that $ \mathcal{B}_c $ is closed under $\rightarrow$ and contains $1,0$. Moreover, since $h$ preserves all operations in $\{1,0,\rightarrow, \land  \}$, and since $ \mathcal{A}_c{\upharpoonright} \{1,0,\rightarrow, \land  \}$ is a Brouwerian semilattice, it follows that $\mathcal{B}_c$ is also a Brouwerian semilattice. 
	
	We augment $\mathcal{B}_c$ with a join operator so that it becomes a Heyting algebra: we define, for all $a,b\in \mathcal{B}_c$, $ a \otimes b := \bigwedge\{ x \in \mathcal{B}_c : x\geq a,b  \}$ and we write $\mathcal{B}_c^\otimes$ for the expansion of $\mathcal{B}_c$ with this new operation. Notice that $a \otimes b$ always exists by the finiteness of $\mathcal{B}_c$. It is  easy to verify that $\mathcal{B}_c^\otimes$ is a Heyting algebra.
	
	We can now use the tensor operator defined over $\mathcal{B}_c^\otimes$ to extend $\mathcal{B}$ to a dependence algebra $\mathcal{B}^\otimes$.  We let $\mathcal{B}^\otimes=(\mathsf{dom}(\mathcal{B}),  \mathcal{B}_c^\otimes, \land, \lor, \otimes, \to, 0 )$, such that $\mathcal{B}^\otimes {\upharpoonright}\{\land,\lor,\to,0 \}=\mathcal{B}$ and where the tensor operator is interpreted in $\mathcal{B}^\otimes$ as follows:
	\begin{itemize}
		\item	For all $a,b\in \mathcal{B}^\otimes_c  $ we let $ a \otimes b := \bigwedge\{ x \in \mathcal{B}^\otimes_c : x\geq a,b  \}$,  as we defined above.
		\item  For all $x,y\in \mathcal{B}^\otimes\setminus \mathcal{B}^\otimes_c  $ we let $ x \otimes y := \bigvee \{ a\otimes b : a\leq x, b\leq y  \text{ and } a,b\in \mathcal{B}^\otimes_c  \}$.
	\end{itemize}	
	
	\begin{claim}
		$\mathcal{B}^\otimes$  is a finite, core-generated and well-connected $\inqint^\otimes$-algebra.
	\end{claim}
	\begin{proof}\renewcommand{\qedsymbol}{$\rfloor$}
		Since $\mathcal{A}$ is finite and $\mathcal{B}^\otimes=h[\langle \mathcal{A}\rangle ]$, it follows that $ \mathcal{B}^\otimes $ is finite.  Also, since $\mathcal{B}^\otimes=\mathsf{dom}(\mathcal{B})=h[\langle\mathcal{A}\rangle]$ and $\mathcal{B}^\otimes_c=h[\mathcal{A}_c]$, we have that $\mathcal{B}^\otimes=\langle \mathcal{B}^\otimes_c\rangle $. That $\mathcal{B}^\otimes$ is well-connected follows from the fact that $\mathcal{B}$ is finite and subdirectly irreducible.
		
		We next verify that $\mathcal{B}$ satisfies the conditions of Definition \ref{depe}. By construction, both $\mathcal{B}^\otimes$ and $\mathcal{B}^\otimes_c$ are Heyting algebras, so it suffices to verify that $\mathcal{B}^\otimes$ also satisfies the additional axioms \textit{Dist}, \textit{Split} and \textit{Mon}. We only check \textit{Dist}, as the proof that \textit{Split} holds is the same of Theorem \ref{birkhoff2} and \textit{Mon} follows similarly.
		%		
		%		\begin{itemize}
		%			\item We check the \textit{Split} axiom. Let $a \in \mathcal{B}^\otimes_c$, $x,y\in \langle \mathcal{B}^\otimes_c \rangle$. Since $h$ is surjective and $h[\mathcal{A}_c]=\mathcal{B}_c^\otimes$, there are $a' \in \mathcal{A}_c$, $x',y'\in  \mathcal{A}=\langle \mathcal{A}_c\rangle$ such that $h(a')=a$, $h(x')=x$ and $h(y')=y$. We proceed as follows. Since $\mathcal{A}$ is an dependence algebra  we have that:
		%			$$ a' \rightarrow (x' \lor y')  = a' \rightarrow (x' \lor y')  .$$
		%			\noindent Therefore
		%			$$ h[a' \rightarrow (x' \lor y')]  = h[ a' \rightarrow (x' \lor y')] . $$
		%			\noindent Hence, by $h$ being a Heyting homomorphism:
		%			$$ h(a') \rightarrow (h(x') \lor h(y'))  = h(a') \rightarrow (h(x') \lor h(y')).  $$
		%			\noindent And finally:
		%			$$ a \rightarrow (x \lor y)  = a \rightarrow (x 	\lor y) , $$
		%			\noindent which proves our claim.
		
		Since $\mathcal{B}$ is generated by its core $\mathcal{B}_c$, we have for $x,y,z\in \langle \mathcal{B}^\otimes_c \rangle$  that $ y=\bigvee_{i\leq n}k_i$ and $ z=\bigvee_{j\leq m}l_j  $. We then obtain:
		\begin{align*}
		x\otimes (y \lor z) & = \bigvee \{ a\otimes b : a\leq x , b \leq y \lor z \text{ and } a,b\in \mathcal{B}^\otimes_c   \} \\ 
		& = \bigvee \{ a\otimes b : a\leq x, b \leq \bigvee_{i\leq n}k_i\lor \bigvee_{j\leq m}l_j \text{ and } a,b\in \mathcal{B}^\otimes_c   \}.  &  
		\end{align*}
		\noindent Now, by the \textit{Split} axiom, together with the fact that $\mathcal{B}$ is well-connected, it follows as in Proposition \ref{wellconnectedsi} that $b$ is join-irreducible, hence  $b \leq \bigvee_{i\leq n}k_i\lor \bigvee_{j\leq m}l_j $ if and only if $b\leq k_i$ or $b\leq l_j$ for some $i\leq n, j\leq m$. Proceeding from the former equalities, we obtain:			
		\begin{align*}
		&= \bigvee  \{ a\otimes b : a\leq x , b \leq \bigvee_{i\leq n}k_i, a,b\in \mathcal{B}^\otimes_c   \}    \lor   \bigvee  \{ a\otimes b : a\leq x , b \leq \bigvee_{j\leq m}l_j, a,b\in \mathcal{B}^\otimes_c   \} \\ 
		&= \bigvee \{ a\otimes b : a\leq x , b \leq y  \text{ and } a,b\in \mathcal{B}^\otimes_c   \}    \lor   \bigvee \{ a\otimes b : a\leq x , b \leq z  \text{ and } a,b\in \mathcal{B}^\otimes_c   \} \\ 
		&= (x\otimes y) \lor (x\otimes z). 
		\end{align*}
		
		\noindent Thus $\mathcal{B}^\otimes$ validates \textit{Split}. Now, since $h:\mathcal{A}\rightarrow \mathcal{B}^\otimes_c$ is a surjective Heyting homomorphism and  $h[ \mathcal{A}_c]\subseteq \mathcal{B}_c$, it follows from Corollary \ref{homlemma2} that $h$ is a $\inqI^\otimes$-homomorphism. Then, by Proposition \ref{closure3}, we have that $h$ preserves the validity of every  $\phi\in \Lambda$ and that $\mathcal{B}^\otimes$ is a $\inqint^\otimes$-algebra.
	\end{proof}
	
	Finally, let $\nu=h\circ \mu$, then  since $h(\phi^\mathcal{A})=s_B$ we have $(\mathcal{B}^\otimes, \nu) \nvDash^c  \phi$, which shows $\mathcal{B}^\otimes\nvDash^c \phi$. Since $\mathcal{B}^\otimes$ is a finite, core-generated and well-connected $\inqI^\otimes$-algebra, this proves our theorem.	
\end{proof}

By combining the two previous results we obtain the following theorem, which is a restricted version of Theorem \ref{birkhoff2} for dependence algebras. We denote by $\mathsf{InqAlg_{W}^\otimes}$ the category of well-connected dependence algebras and $\inqI^\otimes$-homomorphisms. 

\begin{theorem}\label{birkhoff4}
	Suppose $\mathsf{InqAlg_{W}^\otimes} \nvDash^c \phi$,  then there is a finite, core-generated and well-connected dependence algebra $\mathcal{D}$ such that $\mathcal{D} \nvDash^c \phi$.
\end{theorem}
\begin{proof}
	Suppose $\mathsf{InqAlg_{W}^\otimes} \nvDash^c \phi$, then there is some well-connected $\inqI^\otimes$-algebra $\mathcal{A}$ such that  $\mathcal{A} \nvDash^c \phi$. Then, by Theorem \ref{finiteness3}, there is a finite $\inqI^\otimes$-algebra $\mathcal{B}$ such that  $\mathcal{B} \nvDash^c \phi$, and by Proposition \ref{birkhoff3} we can then find a finite, core-generated and well-connected $\inqI^\otimes$-algebra $\mathcal{D}$ such that $\mathcal{D} \nvDash^c \phi$.
\end{proof}

%As we have remarked  $\ipc^\otimes + \Lambda =\inqint^\otimes{\upharpoonright} \langInqst  $ is an intermediate logic with $\otimes$ as disjunction. We say that $\ipc^\otimes + \Lambda $ is \textit{locally tabular} if its corresponding variety of Heyting algebras $Var(\ipc^\otimes + \Lambda   )$ is locally finite. 

We say that a dependence algebra $\mathcal{A}$ is \textit{locally finite} if  for every $X\subseteq \mathsf{core}(\mathcal{A})$ such that $|X|< \omega$ we have $|\langle X\rangle| < \omega$.  A class of $\inqI^\otimes$-algebras $\mathsf{C}$ is locally finite if every $\mathcal{A}\in\mathsf{C}$ is locally finite. We say that an intermediate dependence logic $\inqint^\otimes$ is  \textit{locally tabular} if $\mathsf{InqAlg\Lambda^\otimes}$ is locally finite. When $\inqint^\otimes $ is locally tabular, we can replicate the proof that $\inqint$-algebras are locally finite for the setting of $\inqint^\otimes$-algebras.

\begin{theorem}[Finite Model Property]\label{finiteness2}
	Suppose $\mathcal{A}$ is a locally finite $\inqint^\otimes$-algebra such that $\mathcal{A} \nvDash^c \phi$, then there is a finite $\inqint^\otimes$-algebra $ \mathcal{B} $ such that $\mathcal{B} \nvDash^c \phi$.
\end{theorem}
\begin{proof}
	The proof of this theorem is the same to  that of Theorems \ref{localfiniteness} and \ref{finiteness}, by using the fact that $\mathcal{A}$ is locally finite in place of Diego's theorem.
\end{proof}

\noindent Together with Proposition \ref{birkhoff3}, this provides us a stronger result for locally tabular intermediate dependence logic.

\begin{theorem}\label{birkhoff6}
	Suppose $\inqint^\otimes$ is locally tabular and $\mathsf{InqAlg\Lambda}^\otimes \nvDash^c \phi$,  then there is a finite, core-generated, well-connected inquisitive algebra $\mathcal{D}$ such that $\mathcal{D} \nvDash^c \phi$.
\end{theorem}
\begin{proof}
	By Proposition \ref{birkhoff3} and Theorem \ref{finiteness2}.
\end{proof}

\section{Canonical Constructions and Representation of Inquisitive and Dependence Algebras}\label{four}

In the previous section we have proved several model-theoretic properties of inquisitive and dependence algebras and we studied finite, core-generated inquisitive and dependence algebras. In this section we focus on these classes and we prove some duality results.

Firstly, we prove that the category of finite posets with Köhler morphisms is dually equivalent to the category of finite, core-generated, well-connected inquisitive algebras. Then, we prove that the category of finite posets with p-morphisms is equivalent to the category of finite, core-generated, well-connected dependence algebras. These two results give an important representation of finite, core-generated and well-connected inquisitive and dependence algebras in terms of appropriate downset algebras.

In the second part of this section we then extend these categorical equivalences to teams over Kripke frames and models over inquisitive and dependence algebras. We conclude by remarking that, in this way, we obtain an alternative proof of algebraic completeness for some intermediate inquisitive and dependence logics.

\subsection{Birkhoff, Köhler, and Esakia Duality}

We say that two categories $\mathsf{C}$ and $\mathsf{D}$ are \textit{equivalent} if there are functors $F:\mathsf{C}\to \mathsf{D}$ and $G:\mathsf{D}\to \mathsf{C}$ such that $F\circ G \cong 1_{\mathsf{D}}$ and $G\circ F \cong 1_{\mathsf{C}}$. If $\mathsf{C}$ and $\mathsf{D}$ are  equivalent, we also write $\mathsf{C}\cong\mathsf{D}$. We say that two categories $\mathsf{C}$ and $\mathsf{D}$ are \textit{dually equivalent} and we write $\mathsf{C}\cong^{op}\mathsf{D}$ if $\mathsf{C}\cong\mathsf{D}^{op}$, namely if $\mathsf{C}$ is equivalent to the dual category of $\mathsf{D}$. We refer the reader to \cite{MacLane1978} for a precise definition of these categorical notions and to \cite{pries} for a discussion of several duality results.

We recall some important representation results that we shall use later in our proofs.  Firstly, we recall the following theorem by Birkhoff, which allows us to represent finite bounded distributive lattice in terms of suitable downsets algebras.  

\begin{theorem}[Birkhoff]\label{Birkhoffduality}
	Every finite bounded distributive lattice $\mathcal{L}$ is isomorphic to the  algebra of nonempty downsets of some finite poset: $\mathcal{L}\cong (\mathsf{Dw^+}(P), \land , \lor, 0 )$ for some poset $P$.  In particular, $P$ is the poset of all join-irreducible elements of $\mathcal{L}$ and $h:x\mapsto \{ y\in\mathcal{L}_{ji} : y\leq x  \}$ is the underlying isomorphism.
\end{theorem}

\noindent We shall not discuss here how to extend the previous representation result to suitable morphisms, as we will only need the previous version of Birkhoff's result in the subsequent on this section. We refer the interested reader to \cite{Burris.1981, pries} for more details about Birkhoff's theorem. See also the notes by Morandi \cite{morandi} on dualities in lattice theory for the extension of the previous result to a full categorical equivalence.

Let $\mathfrak{F}$ and $\mathfrak{G}$ be posets (i.e. intuitionistic Kripke frames) with orders $R,R'$ respectively, then we say that a partial function $p:\mathfrak{F}\to\mathfrak{G}$ is a Köhler map if the following hold:
\begin{align*}
(i) \; & \; \forall x,y \in \mathsf{dom}(p) \; [ xR y \Rightarrow p(x)R' p(y)   ]; \\
(ii) \; & \; \forall x\in \mathsf{dom}(p), \; \forall y'\in \mathfrak{G} \; [ p(x) R' y' \Rightarrow \exists y\in\mathfrak{F} \; \text{ such that } xR y \text{ and }   p(y)=y'    ].
\end{align*}

\noindent Thus a Köhler map is essentially a partial p-morphism. Notice that here we followed Bezhanishvili and Jansana \cite{Bezhanishvili2011} and we rephrased Köhler's original conditions with their equivalent conditions for the dual of his order. The reason is that, like Bezhanishvili and Jansana, we prefer to work with upsets of posets rather than downsets. We then denote by \textsf{Pos} the category of all posets (i.e. intuitionistic Kripke frames) with Köhler morphisms. Although posets and Kripke frames are the same objects, we usually talk of posets when the underlying morphisms are Köhler maps, and of Kripke frames when the underlying morphisms are p-morphisms.

\begin{theorem}[Köhler]\label{kohlerrepresentation}
	The category $ \mathsf{Pos}_{F} $ of finite posets and Köhler maps is dually equivalent to the category $\mathsf{BS}_{F}$ of finite Brouwerian semilattices and Brouwerian homomorphisms.
\end{theorem}

\noindent We denote by $\mathsf{O}:  \mathsf{Pos}_{F} \to  \mathsf{BS}_{F}$ and $\mathsf{Pf}: \mathsf{BS}_{F} \to  \mathsf{Pos}_{F}$ the underlying functors of Köhler's duality. For any poset $\mathfrak{F}$, we write $\mathsf{Up}(\mathfrak{F})$ for its collection of upward-closed sets. It is easy to see that $\mathsf{Up}(\mathfrak{F})$ forms a Bouwerian semilattice under the subset ordering $\subseteq$. Then, for any objects $\mathfrak{F},\mathfrak{G}\in  \mathsf{Pos}_{F} $ and any Köhler map $p: \mathfrak{F}\to \mathfrak{G}$, we have: 
\begin{align*}
\mathsf{O}:&\; \mathfrak{F}\mapsto \mathsf{Up}(\mathfrak{F});  \\
\mathsf{O}(p):&\;   \mathsf{Up}(\mathfrak{G})\to \mathsf{Up}(\mathfrak{F}); \; U\mapsto   R[p^{-1}(U)] .
\end{align*}

\noindent Notice that, since $p$ is a partial map, $p^{-1}(U)$ does not need to be upward-closed and thus we explicitly close it under R.

Conversely, we let $\mathsf{Pf}: \mathsf{BS}_{F} \to  \mathsf{Pos}_{F}  $ be the functor from finite Brouwerian semilattices to finite posets which acts as follows.  We say that a subset of a Brouwerian semilattice $F\subseteq \mathcal{B}$ is a \textit{prime filter} if $F$ is a proper filter such that, for all filters $P,Q$, $P\cap Q\subseteq F$ entails $P\subseteq F$ or $Q\subseteq F$ (See \cite{10.2307/41475171}).  Then, if $\mathcal{A}$ is a brouwerian semilattice, we denote by $\mathsf{Pf}(\mathcal{A})$ the poset of its prime filters. For all Brouwerian semilattices $ \mathcal{A},\mathcal{B}\in \mathsf{BS}_{F} $ and for all Brouwerian homomorphisms $h: \mathcal{A}\to\mathcal{B}$, we then  have:
\begin{align*}
\mathsf{Pf}:&\; \mathcal{A}\mapsto \mathsf{Pf}(\mathcal{A}) ;\\
\mathsf{Pf}(h):& \; \mathsf{Pf}(\mathcal{B})\to \mathsf{Pf}(\mathcal{A}); \; :x\mapsto  h^{-1}[G].
\end{align*}

%:x\mapsto\text{max}\{ y\in \mathsf{M}(\mathcal{B}) : x \in h( y^\uparrow  )  \}.

\noindent We refer the reader to the original presentation by Köhler in \cite{kohler1981brouwerian} for the proof that these maps are both well-defined and that they prove that $\mathsf{BS_F}\cong\mathsf{Pos_F}$. However, notice that our maps differ from Köhler's one as we are working with their duals and following \cite{Bezhanishvili2011}. At the same time, our maps are simpler than Bezhanishvili's and Jansana's maps as we are simply adapting Köhler's duality, while they generalise it to infinite Brouwerian semilattices and infinite posets. See also \cite[\S 6]{Bezhanishvili2011}. %inally, we  remark that the underlying isomorphism of Köhler's duality are $\alpha:\mathcal{A}\to \mathsf{Up}\mathsf{Pf}(\mathcal{A}) $

Finally, let us recall the following finite version of Esakia duality.

\begin{theorem}[Esakia]\label{esakiaduality}
	The category $ \mathsf{KF}_{F} $ of finite Kripke frames and p-morpisms is dually equivalent to the category $\mathsf{HA}_{F}$ of finite Heyting algebras and Heyting homomorphisms.
\end{theorem}

\noindent In one direction, we send finite Kripke frames (i.e. posets) to the Heyting algebra of their upsets with reverse ordering, exactly as we did for Köhler's duality. In the case of maps, however, the functor has a simpler description, since the preimage $p^{-1}(U)$ of some up-sets $U$ under a p-morphism $p$ is always an upset. For any objects $\mathfrak{F},\mathfrak{G}\in  \mathsf{KF}_{F} $ and any p-morphism $p: \mathfrak{F}\to \mathfrak{G}$, we have that: 
\begin{align*}
O:&\; \mathfrak{F}\mapsto \mathsf{Up}(\mathfrak{F}) ;  \\
O(p):&\;  \mathsf{Up}(\mathfrak{G})\to \mathsf{Up}(\mathfrak{F}); \; U\mapsto p^{-1}[U] .
\end{align*}

\noindent In the converse direction, we define the functor in the following way: if $\mathcal{L}$ is a bounded lattice, then a subset $F\subseteq \mathsf{dom}(\mathcal{L})$ is a prime filter if $F$ is a proper filter and, for all $x\lor y \in F$, either $x\in F$ or $y\in F$.  We denote by $\mathsf{Pf}(\mathcal{H})$ the set of all prime filters over $\mathcal{H}$. It is easy to check that this forms a poset under the subset ordering. We then define the functor $\mathsf{Pf}: \mathsf{BS}_{F} \to  \mathsf{Pos}_{F}$  by letting, for all $\mathcal{A},\mathcal{B}\in \mathsf{HA_F}$, and for all Heyting homomorphisms $h:\mathcal{A}\to\mathcal{B} $:
\begin{align*}
\mathsf{Pf}:&\; \mathcal{A}\mapsto \mathsf{Pf}(\mathcal{A}); \\
\mathsf{Pf}(h):& \; \mathsf{Pf}(\mathcal{B})\to \mathsf{Pf}(\mathcal{A}); \; G \mapsto h^{-1}[G] .
\end{align*}

\noindent We refer the reader to \cite{9783030120955} for a proof that these functors do indeed provide the category equivalence $\mathsf{KF_F}\cong\mathsf{HA_F}$. 

\subsection{Duality between Posets and Inquisitive Algebras}

We prove in this section that  $\mathsf{Pos_F}\cong^{op}\mathsf{InqAlg_{FCGW}}$. We first describe how, from a finite intuitionistic Kripke frame, one can obtain a finite, core-generated, well-connected inquisitive algebra. Our construction builds on Köhler's and Esakia's dualities. Given a poset $\mathfrak{F}=(W,R)$, we first build the set of all $R$-upsets $\mathsf{Up}(\mathfrak{F})$ and we then consider the set $\mathsf{Dw^+}(\mathsf{Up}(\mathfrak{F}))$ of all non-empty downsets of $\mathsf{Up}(\mathfrak{F})$ ordered by inclusion. % Notice that the first part of such construction is the standard functor of finite Esakia duality \cite{9783030120955}.

Given any poset $\mathfrak{F}$, we define the set $\mathsf{Up}(\mathfrak{F})$ as follows:
$$  \mathsf{Up}(\mathfrak{F})  :=   \{ t\subseteq W  : \text{ if }  x\in t \text{ and } xRy \text{ then } y\in t   \} .$$
\noindent Namely, $\mathsf{Up}(\mathfrak{F})$ is the set of $R$-upsets over $W$. If we think of subsets of $W$ as teams, then $\mathsf{Up}(\mathfrak{F})$ can be viewed as the set of all $R$-closed teams over $W$. One can then check by a routine argument that  $(\mathsf{Up}(\mathfrak{F}), \cap, \cup, \varnothing)$ forms a bounded distributive lattice, where the underlying order is the subset relation. Also, if we add to this structure a Heyting implication in the usual way, we obtain a Heyting algebra $(\mathsf{Up}(\mathfrak{F}), \land, \lor, \rightarrow, 0)$. %Similarly, one can also verify that $(\mathsf{Up}(\mathfrak{F}), \land, \rightarrow, 0)$ is a Brouwerian semilattice

To obtain an inquisitive algebra, we now construct two algebras starting from $\mathsf{Up}(\mathfrak{F})$. Firstly, we consider the following set:	
$$  \mathsf{Dw^+}(\mathsf{Up}(\mathfrak{F}))  :=   \{ x\subseteq \mathsf{Up}(\mathfrak{F})  :  x\neq\varnothing \text{ and if }    s\subseteq t\in x \text{ then } s\in x   \} .$$
\noindent Elements of $\mathsf{Dw^+}(\mathsf{Up}(\mathfrak{F}))$ are downward closed collections of $R$-upsets over $\mathfrak{F}$. It is immediate to check that  the structure $(\mathsf{Dw^+}(\mathsf{Up}(\mathfrak{F})), \cup, \cap, \{\varnothing\})$ is bounded distributive lattice under the subset ordering. This algebra can then be turned into a Heyting algebra in the usual way, by defining $x\to y \leq z \Longleftrightarrow x \land y \leq z$, for all $x,y\in \mathsf{Dw^+}(\mathsf{Up}(\mathfrak{F}))$. We use the symbols $\land, \lor, \rightarrow, 0$ to refer to the underlying operations over $\mathsf{Dw^+}(\mathsf{Up}(\mathfrak{F}))$.

Secondly, we consider the set of principal downward closed collections of $R$-upsets. Recall that a downset $x$ over the poset $\mathsf{Up}(\mathfrak{F})$ is \textit{principal} if there is some $t\in \mathsf{Up}(\mathfrak{F})$ such that:
$$x = \{ s\in \mathsf{Up}(\mathfrak{F}) : s\leq t   \} .$$
\noindent If $x$ is principal and $x = \{ s\in \mathsf{Up}(\mathfrak{F}) : s\leq t   \}$, we then write $x=\{t\}^\downarrow$. Notice, in particular, that since the underlying order of $\mathsf{Up}(\mathfrak{F}) $ is the subset relation, we have that $ s\leq t $ if and only if $s\subseteq t$, hence $\{t\}^\downarrow=\wp(t)\cap \mathsf{Up}(\mathfrak{F})$. We define:
$$  \mathsf{Dw_p}(\mathsf{Up}(\mathfrak{F}))  :=   \{  \{t\}^\downarrow  :  t\in \mathsf{Up}(\mathfrak{F})  \} .$$
\noindent The poset $(\mathsf{Dw_p}(\mathsf{Up}(\mathfrak{F})),\subseteq ) $ forms a bounded distributive lattice and can be augmented by a Heyting implication in the same way as we did for the previous lattices. The next proposition shows that $\mathsf{Dw_p}(\mathsf{Up}(\mathfrak{F})) $ is isomorphic to $\mathsf{Up}(\mathfrak{F})$.

\begin{proposition} \label{latiso1}
	The following Heyting algebras are isomorphic: $  \mathsf{Dw_p}(\mathsf{Up}(\mathfrak{F})) \cong \mathsf{Up}(\mathfrak{F})$ under the map $h:x\mapsto\{x\}^\downarrow$.
\end{proposition}
\begin{proof}
	Consider the map $h: \mathsf{Up}(\mathfrak{F}) \rightarrow \mathsf{Dw_p}(\mathsf{Up}(\mathfrak{F})) $ such that $h(x)=\{x \}^\downarrow$. It follows by the definition of principal downsets that this map is both surjective and injective.  The, since $ s\subseteq t$ if and only if $   \{s\}^\downarrow \subseteq \{t \}^\downarrow$, it follows that $h$ is a isomorphism.		
\end{proof}

In particular, the previous proposition allows us to characterise operations in $\mathsf{Dw_p}(\mathsf{Up}(\mathfrak{F}))$ by operations in $\mathsf{Up}(\mathfrak{F})$. We have that, for all  $\{a \}^\downarrow, \{b \}^\downarrow\in \mathsf{Dw_p}(\mathsf{Up}(\mathfrak{F}))$:
\begin{align*}
\{a \}^\downarrow\land \{b \}^\downarrow = \{a \land b \}^\downarrow \hspace{15pt}
\{a \}^\downarrow\lor \{b \}^\downarrow = \{a \lor b \}^\downarrow \hspace{15pt}
\{a \}^\downarrow\rightarrow \{b \}^\downarrow = \{a \rightarrow b \}^\downarrow.
\end{align*}

%\begin{align*}
%\{a \}^\downarrow\land \{b \}^\downarrow &= \{a \land b \}^\downarrow\\
%\{a \}^\downarrow\lor \{b \}^\downarrow &= \{a \lor b \}^\downarrow\\
%\{a \}^\downarrow\rightarrow \{b \}^\downarrow &= \{a \rightarrow b \}^\downarrow.
%\end{align*}

\noindent It is clear that $\mathsf{Dw_p}(\mathsf{Up}(\mathfrak{F})) \subseteq \mathsf{Dw^+}(\mathsf{Up}(\mathfrak{F})) $.  We now claim that $(\mathsf{Dw_p}(\mathsf{Up}(\mathfrak{F}),
\land,\rightarrow,0)$ is a subalgebra of $ \mathsf{Dw^+}(\mathsf{Up}(\mathfrak{F})) $ with respect to the operations $\land, \rightarrow, 0 $. We prove the following proposition.

\begin{proposition}\label{SEmilattice}
	$\mathsf{Dw_p}(\mathsf{Up}(\mathfrak{F}))$ is a Brouwerian semilattice subalgebra of $\mathsf{Dw^+}(\mathsf{Up}(\mathfrak{F}))$.
\end{proposition}
\begin{proof}
	We first prove that $\mathsf{Dw_p}(\mathsf{Up}(\mathfrak{F}))$  is closed under $\land, \rightarrow,0$. By construction, $0=\{\varnothing\}$, and $\{\varnothing\} \in \mathsf{Dw_p}(\mathsf{Up}(\mathfrak{F})) $, since $\varnothing\in \mathsf{Up}(\mathfrak{F}) $ and $\{\varnothing \}= \{\varnothing\}^\downarrow$.
	
	For all $x,y\in \mathsf{Dw_p}(\mathsf{Up}(\mathfrak{F}))$ we have that $x=\{a\}^\downarrow$ and $y=\{b  \}^\downarrow$ for some $a,b\in \mathsf{Up}(\mathfrak{F})$. Since $\mathsf{Dw_p}(\mathsf{Up}(\mathfrak{F})) \cong \mathsf{Up}(\mathfrak{F})$, we  have $\{a  \}^\downarrow \land \{b  \}^\downarrow  = \{a \land  b \}^\downarrow $ and $\{a  \}^\downarrow \to \{b  \}^\downarrow = \{a \to  b\}^\downarrow$. Therefore,  $x\land y, x \to y \in \mathsf{Dw_p}(\mathsf{Up}(\mathfrak{F}))$, showing that $\mathsf{Dw_p}(\mathsf{Up}(\mathfrak{F}))$ is closed under $\land$ and $\rightarrow$.  Since $\mathsf{Dw^+}(\mathsf{Up}(\mathfrak{F})){\upharpoonright }\{\land,\to,0  \}$ is a Brouwerian semilattice, it follows that $\mathsf{Dw_p}(\mathsf{Up}(\mathfrak{F}))$ is a Browerian subalgebra of $\mathsf{Dw^+}(\mathsf{Up}(\mathfrak{F})$.
\end{proof}

We can represent the relations between the algebras that we have constructed by the following diagram:
\begin{center}
	\begin{tikzcd}
		\mathfrak{F} \arrow[r, "\mathsf{Up}(\cdot)"] & \mathsf{Up}(\mathfrak{F}) \arrow[r, "\mathsf{Dw^+}(\cdot)"] \arrow[dd, "\mathsf{Dw_p}(\cdot)"'] & \mathsf{Dw^+}(\mathsf{Up}(\mathfrak{F})) \\		
		& &\\
		& \mathsf{Dw_p}(\mathsf{Up}(\mathfrak{F})) \arrow[ruu, hook, "\subseteq"']    &                       
	\end{tikzcd}
\end{center}

We can now use the algebras constructed above to obtain a finite, core-generated, well-connected inquisitive algebra. We first prove that $\mathsf{Dw_p}(\mathsf{Up}(\mathfrak{F})) $ generates the algebra $\mathsf{Dw^+}(\mathsf{Up}(\mathfrak{F})) $.

\begin{proposition}\label{gene}
	$\mathsf{Dw^+}(\mathsf{Up}(\mathfrak{F}))$ is generated by its subset $\mathsf{Dw_p}(\mathsf{Up}(\mathfrak{F})) $.
\end{proposition}
\begin{proof}
	If $x\in \mathsf{Dw^+}(\mathsf{Up}(\mathfrak{F}))$ then, since $\mathfrak{F}$ is finite,  $x\subseteq \mathsf{Up}(\mathfrak{F}) $ is finite too. Let $a_0,\dots,a_n$ be maximal upsets in $x$ -- they exist by the finiteness of $x$. Then,  $\{a_0\}^\downarrow,\dots,\{a_n\}^\downarrow$ are principal downsets, whence $\{a_0\}^\downarrow,\dots,\{a_n\}^\downarrow\in \mathsf{Dw_p}(\mathsf{Up}(\mathfrak{F}))$. We obtain:
	\begin{align*}
	x = \bigcup	\{\{a_0\}^\downarrow,\dots,\{a_n\}^\downarrow\}  = \bigvee (\{a_0\}^\downarrow,\dots,\{a_n\}^\downarrow);
	\end{align*}
	\noindent which means that $x\in \langle \mathsf{Dw_p}(\mathsf{Up}(\mathfrak{F})) \rangle $ and thus $\mathsf{Dw^+}(\mathsf{Up}(\mathfrak{F}))=\langle \mathsf{Dw_p}(\mathsf{Up}(\mathfrak{F})) \rangle $.
\end{proof}

Finally, the next result shows that the structure:
$\mathcal{A}:=(\mathsf{Dw^+}(\mathsf{Up}(\mathfrak{F})), \mathsf{Dw_p}(\mathsf{Up}(\mathfrak{F})), \land,\lor, \rightarrow,0)$ is a finite, core generated, well-connected, inquisitive algebra.

\begin{proposition}\label{canonicalinqalgebra}
	The structure  $\mathcal{A}=(\mathsf{Dw^+}(\mathsf{Up}(\mathfrak{F})), \mathsf{Dw_p}(\mathsf{Up}(\mathfrak{F})), \land,\lor, \rightarrow,0)$ is a finite, core-generated, well-connected, inquisitive algebra.
\end{proposition}
\begin{proof}
	By construction, $(\mathsf{Dw^+}(\mathsf{Up}(\mathfrak{F})), \land,\lor, \rightarrow,0)$ is a Heyting algebra and, by Proposition \ref{SEmilattice}, $\mathcal{A}_c= ( \mathsf{Dw_p}(\mathsf{Up}(\mathfrak{F})), \land, \rightarrow,0  ) $ is a Brouwerian semilattice. Also, we have by Proposition \ref{gene} that $\mathsf{Dw^+}(\mathsf{Up}(\mathfrak{F}))=\langle \mathsf{Dw_p}(\mathsf{Up}(\mathfrak{F})) \rangle $, hence $\mathcal{A}$ is core-generated. It follows immediately by our construction, together with the fact that $\mathfrak{F}$ is finite, that $\mathcal{A}$ is finite as well. 
	
	We now prove that $\mathcal{A}$ is well-connected. Since $\mathcal{A}$ is a finite Heyting algebra, it suffices to show that it has a second greatest element. Let $W$ be the set of all worlds in $\mathfrak{F}$ and let $s=\{W\}^\downarrow \setminus W $. It is clear that $s\in \mathsf{Dw^+}(\mathsf{Up}(\mathfrak{F}))$ and that $s\neq  \{W\}^\downarrow =1$. Now suppose $x\in \mathsf{Dw^+}(\mathsf{Up}(\mathfrak{F}))$ and $x\neq 1$, then for all $R$-upsets $t\in x$, we have that $t\neq W$, hence $t\in s$ and $x\leq s$. Thus $s$ is the second greatest element of $\mathcal{A}$ and $\mathcal{A}$ is well-connected.
	
	Finally, we check that $\mathcal{A}$ verifies the \textit{Split} axiom. Let $a\in \mathcal{A}_c$, then $a = \{ t  \}^\downarrow$ for some $t\in \mathsf{Up}(\mathfrak{F})$. If $\{ t  \}^\downarrow  \subseteq x \lor y$, it follows from the fact that $x,y$ are downward closed that either $\{ t  \}^\downarrow  \subseteq x$ or $\{ t  \}^\downarrow  \subseteq y$, showing  that $a$ is join-irreducible. By reasoning as in the proof of Proposition \ref{wellconnectedsi}, it follows that  $a \rightarrow (x\lor y) = (a \to x)\lor (a \to y)$, which proves that $\mathcal{A}$ satisfies the \textit{Split} axiom and therefore that $\mathcal{A}$ is a finite, core-generated, well-connected, inquisitive algebra.
\end{proof}

Following the construction described so far, we have seen how to obtain, from a finite poset $\mathfrak{F}$, a finite, core-generated, well-connected inquisitive algebra $F(\mathfrak{F})$. To obtain a functor $F:\mathsf{Pos_F}\to\mathsf{InqAlg_{FCGW}}$, it remains to extend this assignment to morphisms between posets, i.e. to Köhler maps. The following proof follows easily from Köhler's duality.

\begin{proposition}\label{canonicalhomo}
	Suppose $p:\mathfrak{F}\rightarrow\mathfrak{G}$ is a Köhler map, then the function $F(p): F(\mathfrak{G})\rightarrow F(\mathfrak{F})  $ such that:
	$$F(p):  \bigcup \{ \{ t_0 \}^\downarrow,\dots \{ t_n \}^\downarrow \} \longmapsto \bigcup  \{ \{R(p^{-1}[t_0])\}^\downarrow \dots \{R(p^{-1}[t_n])\}^\downarrow     \} $$
	\noindent is a inquisitive homomorphism.
\end{proposition}
\begin{proof}  Firstly, notice that  since $F(\mathfrak{G})$ is core-generated, then every element $x\in F(\mathfrak{G})$ is of the form $x=\bigvee_{i\leq n}a_i$ where $a_i\in F(\mathfrak{G})_c$ for all $i\leq n$. Hence, since $F(\mathfrak{G})$ is a downset inquisitive algebra,  every element $x\in F(\mathfrak{G})$ is of the form $x=\bigcup \{ \{ t_0 \}^\downarrow,\dots \{ t_n \}^\downarrow \}$, which means that $F(p)$ is total.  Moreover, since by construction $R (p^{-1}[t] )$ is  an $R$-upset, it follows that $F(p)$ is well-defined.
	
	To see that $F(p)$ is core preserving it suffices to notice that, for all $ \{ t \}^\downarrow \in F(\mathfrak{G})_c$ :
	\begin{align*}
	F(p)( \{ t \}^\downarrow ) =  \{R(p^{-1}[t])\}^\downarrow \in \mathsf{Dw_p}(\mathsf{Up}(\mathfrak{F}))  = F(\mathfrak{F})_c.
	\end{align*}
	
	Now, since $p$ is a Köhler map, it follows by Köhler's duality that $\mathsf{O}(p):\mathsf{Up}(\mathfrak{G})\to\mathsf{Up}(\mathfrak{F})$ such that $\mathsf{O}(p):t\mapsto R[p^{-1}(t)]$ is a Brouwerian semilattice homomorphism. Since $\mathsf{Up}(\mathfrak{F})\cong \mathsf{Dw}_p(\mathsf{Up}(\mathfrak{F}))$ and $\mathsf{Up}(\mathfrak{G})\cong \mathsf{Dw}_p(\mathsf{Up}(\mathfrak{G}))$, it immediately follows that $F(p){\upharpoonright}\mathsf{Dw}_p(\mathsf{Up}(\mathfrak{G}))$ is a Brouwerian semilattice homomorphism. To verify that it is a $\inqI$-homomorphism, it suffices by Theorem \ref{normalformalgebra} to check that $F(p)$ preserves joins of core elements. Then, for any $ \{ t \}^\downarrow,  \{ s \}^\downarrow \in F(\mathfrak{G})_c$, we immediately have by our definition:
	\begin{align*}
	F(p) ( \{ t \}^\downarrow \lor   \{ s \}^\downarrow   ) =    \{ R[p^{-1}(t)] \}^\downarrow \cup   \{ R[p^{-1}(s)] \}^\downarrow =F(p) ( \{ t \}^\downarrow) \lor F(p) ( \{ s \}^\downarrow).
	\end{align*}
	\noindent Which proves our claim.
\end{proof}

\noindent In particular, it follows from the proof of the proposition above that $F(p)$ is always total, even if $p$ is not. In fact, if $R(p^{-1}[t_i])=\varnothing$ for all $i\leq n$, then $\{R(p^{-1}[t_i])\}^\downarrow=\{\varnothing\}$ for all $i\leq n$ and therefore $F(p) \Big (\bigcup \{ \{ t_0 \}^\downarrow,\dots \{ t_n \}^\downarrow \} \Big ) = \{\varnothing \}  $. It is routine to check that the map $F:\mathsf{Pos_F}\to \mathsf{InqAlg_{FCGW}}$  is functorial.

Now, in order to establish the equivalence between $\mathsf{Pos_F}$ and $\mathsf{InqAlg_{FCGW}}$, we need to proceed in the opposite direction, and  define a functor $G$ which associates a finite poset to every  finite, core-generated, well-connected inquisitive algebra. We employ Köhler's duality to prove a representation theorem for finite, core-generated, well-connected,  inquisitive algebras.

\begin{proposition}\label{representation1}
	Let $\mathcal{A}$ be a finite, core-generated, well-connected inquisitive algebra, then there is a finite Kripke frame $\mathfrak{F}$ such that: 	
	\begin{align*}
	\mathcal{A}&\cong(\mathsf{Dw^+}(\mathsf{Up}(\mathfrak{F})),\mathsf{Dw_p}(\mathsf{Up}(\mathfrak{F})), \land,\lor, \rightarrow, 0  ).
	\end{align*}
\end{proposition}
\begin{proof}
	Suppose $\mathcal{A}$ is a finite, core-generated, well-connected inquisitive algebra, and let $\mathcal{A}_{ji}$ be its subset of join-irreducible elements. By Proposition \ref{wellconnectedsi} we have that $\mathcal{A}_{ji}=\mathcal{A}_{c}$. Hence, by the fact that $\mathcal{A}$ is an inquisitive algebra, it follows $\mathcal{A}_{ji}$ is a Brouwerian semilattice in the signature $ \{\land,\rightarrow, 0 \} $. By Theorem \ref{kohlerrepresentation}, we have that $\mathcal{A}_{ji}\cong \mathsf{Up}(\mathfrak{F})$ for some finite poset $\mathfrak{F}$ -- here we shall think of $\mathfrak{F}$ as an intuitionistic Kripke frame, where the underlying ordering is the accessibility relation $R$ between worlds. We let $g:\mathcal{A}_{ji}\to \mathsf{Up}(\mathfrak{F})$ be the function witnessing such isomorphism.
	
	Now, let $\mathcal{B}=(\mathsf{Dw^+}(\mathsf{Up}(\mathfrak{F})),\mathsf{Dw_p}(\mathsf{Up}(\mathfrak{F})), \land,\lor, \rightarrow, 0  )$ be the canonical $\inqI$-algebra obtained from the frame $\mathfrak{F}$ using the construction outlined in the previous section. It follows by Proposition \ref{canonicalinqalgebra} that $\mathcal{B}$ is an inquisitive algebra. We show that $\mathcal{A}\cong \mathcal{B}$.

	%We let $\mathsf{J}(x):=\{ y\in \mathcal{A}_{ji} : y\leq x   \}$.
	
	By Theorem \ref{Birkhoffduality} the function   $h:\mathcal{A} \rightarrow \mathsf{Dw^+}(\mathcal{A}_{ji}  )$ such that $h(x)=\{ y\in \mathcal{A}_{ji} : y\leq x   \}$ is a lattice isomorphism. 	Moreover, since $ \mathsf{Dw^+}(\mathcal{A}_{ji} )$ is a finite bounded distributive lattice, we can expand it with a Heyting implication and obtain the Heyting algebra $(\mathsf{Dw^+}(\mathcal{A}_{ji}), \land,\lor,\to,0 )$. Since $h$ is an order-preserving bijection, it is also a Heyting algebra isomorphism and thus we obtain that $\mathcal{A}\cong \mathsf{Dw^+}(\mathcal{A}_{ji})$. %Therefore, we immediately have that $\mathcal{A}{\upharpoonright}\{\land,\lor,\rightarrow,0  \}\cong \mathcal{B}{\upharpoonright}\{\land,\lor,\rightarrow,0  \}.  $.
	
	Let $\hat{g}: \mathsf{Dw^+}(\mathcal{A}_{ji}) \rightarrow \mathcal{B} $ be defined by lifting $g$ to the algebra $ \mathsf{Dw^+}(\mathcal{A}_{ji} )$:
	$$\hat{g}:  \{a_0,\dots,a_n  \}^\downarrow\mapsto \bigcup \{ \{ g(a_0)\}^\downarrow,\dots,\{g(a_n)\}^\downarrow \}.   $$
	\noindent By the fact that $g$ is a isomorphism, together with the fact that  $\mathcal{B}$ is core-generated, it follows that $\hat{g}$ is a bijection. Moreover, if $\{a_0,\dots,a_n  \}^\downarrow\subseteq \{b_0,\dots,b_m  \}^\downarrow$ we then clearly have that $\bigcup \{ \{ g(a_0)\}^\downarrow,\dots,\{g(a_n)\}^\downarrow \}\subseteq \bigcup \{  g(b_0^\downarrow),\dots,g(b_m^\downarrow) \}$, since $g$ is order preserving. Then, $\hat{g}$ is order preserving and, since it is a bijection, it is a Heyting algebra isomorphism as well.
	
	We then let $f:=\hat{g}\circ h: \mathcal{A}\to \mathcal{B}$. Since $h$ and $\hat{g}$ are both Heyting algebra isomorphisms, it follows that $f$ is also a Heyting algebra isomorphism. Moreover, since $\mathcal{A}$ is core-generated and  well-connected, it follows by Proposition \ref{wellconnectedsi} that $\mathcal{A}_c = \mathcal{A}_{ji}$ and thus,  for all $a\in\mathcal{A}_c$, $h(a)=\{a\}^\downarrow$ and:
	$$f(a)= \hat{g}\circ h(a) =  \hat{g} \;  ( \{ a   \}^\downarrow  )=  \{ g(a)   \}^\downarrow  \in \mathsf{Dw_p}(\mathsf{Up}(\mathfrak{F})).  $$
	\noindent Therefore, $f$ is a bijective, core-preserving, Heyting algebra homomorphism between $\mathcal{A}$ and $\mathcal{B}$, which means that $f$ is a $\inqI$-homomorphism and thus that $\mathcal{A}\cong \mathcal{B}$.
\end{proof}

While the previous theorem gives a representation of finite,  core-generated, well-connected, inquisitive algebras, the following proposition provides a representation of the maps between them. Our proof follows easily from Köhler's duality (see in particular \cite[Lemma 3.2]{kohler1981brouwerian}).

\begin{proposition}\label{representation2}
	Let $h:\mathcal{A}\rightarrow\mathcal{B}$ be a $\inqI$-homomorphism between finite, core-generated, well-connected $\inqI$-algebras, then there is a unique Köhler map ${p:\mathfrak{G}\rightarrow \mathfrak{F} }$ such that $h=F(p)$, $\mathcal{A}= \mathsf{Dw^+}(\mathsf{Up}(\mathfrak{F}))$ and $\mathcal{B}= \mathsf{Dw^+}(\mathsf{Up}(\mathfrak{G})) $. 
\end{proposition}
\begin{proof} Suppose $h:\mathcal{A}\rightarrow\mathcal{B}$ is a $\inqI$-homomorphism between finite, core-generated, well-connected $\inqI$-algebras. By Proposition \ref{representation1} we assume without loss of generality  that $\mathcal{A}= \mathsf{Dw^+}(\mathsf{Up}(\mathfrak{F}))$ and $\mathcal{B}= \mathsf{Dw^+}(\mathsf{Up}(\mathfrak{G})) $ for some finite Kripke frames $\mathfrak{F}$ and $\mathfrak{G}$. 
	
	Let $k:\mathsf{Up}(\mathfrak{F})\to \mathsf{Up}(\mathfrak{G}) $ be the map such that $k(t)= s$ if and only if  $h(\{t\}^\downarrow)=\{s \}^\downarrow $. Since $h{\upharpoonright}\mathsf{Dw}_p\mathsf{Up}(\mathfrak{F})$ is a Brouwerian semilattice homomorphism, we obtain by Proposition \ref{latiso1} that $k$ is a Browuerian semilattice homomorphism as well. By Köhler's duality there is a unique Köhler's map $p:\mathfrak{G}\to\mathfrak{F}$ such that  $\mathsf{O}(p)=k$. In addition, for any $\bigcup \{ \{ t_0 \}^\downarrow,\dots \{ t_n \}^\downarrow \}\in \mathsf{Dw^+}(\mathsf{Up}(\mathfrak{F})) $:
	\begin{align*}
	F(p) \Big (  \bigcup \{ \{ t_0 \}^\downarrow,\dots \{ t_n \}^\downarrow \} \Big ) & =  \bigcup  \{ \{R(p^{-1}[t_0])\}^\downarrow \dots \{R(p^{-1}[t_n])\}^\downarrow     \}  \\
	& =  \bigcup  \{   \{  \mathsf{O}(p)(t_0)  \}^\downarrow \dots \{ \mathsf{O}(p)(t_n)  \}^\downarrow     \}\\
	& =  \bigcup  \{   \{ k(t_0)  \}^\downarrow \dots \{ k(t_n)  \}^\downarrow     \}\\
	& =  \bigcup  \{   h(\{ t_0 \}^\downarrow) \dots h(\{ t_n  \}^\downarrow )    \}\\
	& = h \bigcup  \{   \{ t_0 \}^\downarrow \dots \{ t_n  \}^\downarrow     \}.
	\end{align*}	
	\noindent which completes the proof of our claim.
\end{proof}

Now, let $G:\mathsf{InqAlg_{FCGW}}\rightarrow\mathsf{Pos_F}  $ be the functor defined as follows. On objects, we let $G(\mathcal{A})=\mathfrak{F}$, where $\mathfrak{F}$ is a finite poset such that $\mathcal{A}\cong F(\mathfrak{F})$.  On inquisitive homomorphisms, we let $G(h)=p:\mathfrak{F}\to\mathfrak{G}$ where $F(p):F(\mathfrak{G})\rightarrow F(\mathfrak{F})$. It follows by this very construction that $G$ is functorial and that together with  $F$ it forms a dual equivalence between  $\mathsf{Pos_F}$ and  $\mathsf{InqAlg_{FCGW}}$. 

\begin{theorem}\label{duality}
	The category of finite posets with Köhler maps is dually equivalent to the category of finite, core-generated, well-connected inquisitive algebras with $\inqI$-homomorphisms: $$\mathsf{Pos_F}\cong^{op}\mathsf{InqAlg_{FCGW}}.$$	
\end{theorem}
\begin{proof}
	Let $F:\mathsf{Pos_F}\rightarrow \mathsf{InqAlg_{FCGW}} $ and $G:\mathsf{InqAlg_{FCGW}}\rightarrow\mathsf{KF_F}  $ be the functors defined above. It follows from Proposition \ref{canonicalinqalgebra} and Proposition \ref{canonicalhomo} that $F$ is well-defined and from Proposition \ref{representation1} and Proposition \ref{representation2} that $G$ is well-defined. By our definitions of $F$ and $G$ we have that $F\circ G\cong id$ and $G\circ F \cong id$. Therefore, $F$ and $G$ describe a dual categorical equivalence between $\mathsf{Pos_F}$ and $\mathsf{InqAlg_{FCGW}}$.
\end{proof}

\subsection{Duality between Kripke Frames and Dependence Algebras}

In the previous section we have proved that the category of finite posets is dually equivalent to the category of finite, core-generated, well-connected inquisitive algebras.Here we prove that the category $\mathsf{InqAlg_{FCGW}^\otimes}$ of finite, core-generated, well-connected dependence algebras is equivalent to $\mathsf{KF_F}$, thus obtaining a similar result for dependence algebras.

We first describe the functor $F:  \mathsf{KF_F}\rightarrow\mathsf{InqAlg_{FCGW}^\otimes}$. For any finite kripke frame $\mathfrak{F}$, we let $F(\mathfrak{F})$ be the $\inqI^\otimes$-algebra obtained by adding a tensor operator $\otimes$ to the inquisitive algebra $\mathsf{Dw^+}(\mathsf{Up}(\mathfrak{F}))$. For all $a,b\in \mathsf{Dw_p}(\mathsf{Up}(\mathfrak{F}))$, we have $a=\{t\}^\downarrow$ and $b=\{s\}^\downarrow$ for some $t,s\in \mathsf{Up}(\mathfrak{F})$. We then let:
\[a\otimes b = \{ t \cup s  \}^\downarrow. \]
\noindent And we lift such operation to all $x,y \in \mathsf{Dw^+}(\mathsf{Up}(\mathfrak{F}))$ as follows:
$$x\otimes y := \bigvee \{ a \otimes b : a\leq x , b \leq y \text{ and } a,b\in \mathsf{Dw_p}(\mathsf{Up}(\mathfrak{F}))   \}. $$

\noindent The next result shows that $F$ is well-defined on objects.

\begin{proposition}\label{canonicaldepalgebra}
	The structure  $\mathcal{A}=(\mathsf{Dw^+}(\mathsf{Up}(\mathfrak{F})), \mathsf{Dw_p}(\mathsf{Up}(\mathfrak{F})), \land,\lor, \rightarrow,\otimes,0)$ is a finite, core-generated, well-connected, dependence algebra.
\end{proposition}
\begin{proof}
	
	Firstly, we have by Proposition \ref{canonicalinqalgebra} that $\mathcal{A}$ is a finite, core-generated, well-connected, inquisitive algebra. Hence it suffices to show that $\mathcal{A}$ is also a dependence algebra. 
	
	Now, by Proposition \ref{latiso1}, we have that $a\otimes b = \{ t \lor s  \}^\downarrow $ is a well-defined join operator over $\mathsf{Dw_p}(\mathsf{Up}(\mathfrak{F}))$, hence the core $(\mathsf{Dw_p}(\mathsf{Up}(\mathfrak{F})), \land, \otimes, \to, 0)$ is a Heyting algebra. 
	
	It remains to verify that $\mathcal{A}$ validates the axioms \textit{Dist} and \textit{Mon}. We only show that \textit{Dist} holds, as \textit{Mon} is easily checked in a similar way. By Theorem \ref{normalformalgebra} we have $y=\bigvee_{i\leq n}k_i$ and $z=\bigvee_{j\leq m}l_j$ with $k_i, l_j\in \mathsf{Dw_p}(\mathsf{Up}(\mathfrak{F}))$ for all $i\leq n, j\leq m$. We then obtain:
	\begin{align*}
	x\otimes (y \lor z) & = \bigvee \{ a\otimes b : a\leq x , b \leq y \lor z \text{ and } a,b\in \mathsf{Dw_p}(\mathsf{Up}(\mathfrak{F}))  \} \\ 
	& = \bigvee \{ a\otimes b : a\leq x, b \leq \bigvee_{i\leq n}k_i\lor \bigvee_{j\leq m}l_j \text{ and } a,b\in \mathsf{Dw_p}(\mathsf{Up}(\mathfrak{F}))  \}.  
	\end{align*}
	\noindent Now, since by Proposition \ref{wellconnectedsi} core elements of well-connected dependence algebras are join-irreducible, we have that $b \leq \bigvee_{i\leq n}k_i\lor \bigvee_{j\leq m}l_j $ if and only if $b\leq k_i$ or $b\leq l_j$ for some $i\leq n, j\leq m$. Then, proceeding from the former equalities:			
	\begin{align*}
	&= \bigvee  \{ a\otimes b : a\leq x , b \leq \bigvee_{i\leq n}k_i, a,b\in \mathcal{A}_c  \}    \lor   \bigvee  \{ a\otimes b : a\leq x , b \leq \bigvee_{j\leq m}l_j, a,b\in \mathsf{Dw_p}(\mathsf{Up}(\mathfrak{F}))   \} \\ 
	&= \bigvee \{ a\otimes b : a\leq x , b \leq y  \text{ and } a,b\in \mathcal{A}_c  \}    \lor   \bigvee \{ a\otimes b : a\leq x , b \leq z  \text{ and } a,b\in \mathsf{Dw_p}(\mathsf{Up}(\mathfrak{F}))   \} \\ 
	&= (x\otimes y) \lor (x\otimes z). 
	\end{align*}
	
	\noindent Hence $\mathcal{A}$ is a $\inqI^\otimes$-algebra. 
\end{proof}

Now, let $p:\mathfrak{F}\to\mathfrak{G}$ be a p-morphism. We proceed as in the case of inquisitive algebras and we show how to obtain a canonical dependence homomorphism $F(p)$.

\begin{proposition}\label{canonicaldepehomo}
	Suppose $p:\mathfrak{F}\rightarrow\mathfrak{G}$ is a p-morphism, then the function $F(p): F(\mathfrak{G})\rightarrow F(\mathfrak{F}) $ such that:
	$$F(p):  \bigcup \{ \{ t_0 \}^\downarrow,\dots \{ t_n \}^\downarrow \} \longmapsto \bigcup  \{ \{p^{-1}[t_0]\}^\downarrow \dots \{p^{-1}[t_n]\}^\downarrow     \} $$
	\noindent is a dependence homomorphism.
\end{proposition}
\begin{proof}
	The proof is analogous to that of Proposition \ref{canonicalhomo}, by using Esakia's duality in place of Köhler's duality.
\end{proof}

It is then easy to verify that $F$ is functorial. Then, to obtain a categorical equivalence, we prove the following representation results.

\begin{proposition}\label{representation3}
	Let $\mathcal{A}$ be a finite, core-generated, well-connected dependence algebra, then there is a finite Kripke frame $\mathfrak{F}$ such that $\mathcal{A}\cong F(\mathcal{A})$.
\end{proposition}	
\begin{proof}
	By proceeding exactly as in the proof of Proposition \ref{representation1}, using Esakia's duality instead of Köhler's duality, we obtain that there is a finite Kripke frame $\mathfrak{F}$ such that the following map is a isomorphism:
	$$f:
	\mathcal{A}\to(\mathsf{Dw^+}(\mathsf{Up}(\mathfrak{F})),\mathsf{Dw_p}(\mathsf{Up}(\mathfrak{F})), \land,\lor, \rightarrow, 0); \bigvee_{i\leq n} a_i \mapsto  \bigcup \{ \{g(a_0)\}^\downarrow,\dots,  \{g(a_n)\}^\downarrow  \} .$$
	\noindent where  $g:\mathcal{A}_{ji}\to \mathsf{Up}(\mathfrak{F})$ is a isomorphism of Heyting algebras. Therefore, we have that for all $a,b\in \mathcal{A}_{ji}=\mathcal{A}_c$: $ g(a\otimes b)=g(a)\otimes g(b)= g(a)\cup g(b)$, thus: 
	\begin{align*} f(a\otimes b)=\{ g(a\otimes b)\}^\downarrow=\{g(a)\}^\downarrow\cup \{g(b)\}^\downarrow.
	\end{align*}
	\noindent which proves that $f$ preserve the tensor operation for core elements. 
	
	Then, since for all $x,y\in \mathcal{A} $, we have that $x=\bigvee_{i\leq n}a_i$ and $y=\bigvee_{j\leq m}b_j$ for  $a_i,b_j\in \mathcal{A}_c$ for all $i\leq n, j\leq m$, one can proceed as in the proof of Lemma \ref{homlemma} and verify that $ f(x\otimes y)=f(x)\otimes f(y) $ for all $x,y\in \mathcal{A} $.
\end{proof}

\noindent Similarly, we have the following representation of p-morphisms. 

\begin{proposition}\label{representation4}
	Let $h:\mathcal{A}\rightarrow\mathcal{B}$ be a $\inqI^\otimes$-homomorphisms between finite, core-generated, well-connected $\inqI^\otimes$-algebras, then there is a unique p-morphism ${p:\mathfrak{F}\rightarrow \mathfrak{G} }$ such that $h=F(p)$. 
\end{proposition}
\begin{proof}
	The proof is the same as the proof of  Proposition \ref{representation2}, by using Esakia's duality instead of Köhler's duality. 
\end{proof}

Then, let  $G:\mathsf{InqAlg_{FCGW}^\otimes}\to\mathsf{KF_F}$ be the functor defined as follow. On objects, we let $G(\mathcal{A})=\mathfrak{F}$, where $\mathfrak{F}$ is a finite Kripke frame such that $\mathcal{A}\cong F(\mathcal{A})$.  On $\inqI^\otimes$-homomorphisms, we let $G(h)=p:\mathfrak{F}\to\mathfrak{G}$ where $F(p):F(\mathfrak{G})\rightarrow F(\mathfrak{F})$. It follows that $G$ is functorial and that together with  $F$ it forms a dual equivalence between  $\mathsf{KF_F}$ and  $\mathsf{InqAlg_{FCGW}^\otimes}$.

\begin{theorem} \label{equiv5} 
	The category of finite Kripke frames with p-morphisms is dually equivalent to the category of finite, core-generated, well-connected $\inqI^\otimes$-algebras with $\inqI^\otimes$-homomorphisms: $$\mathsf{KF_F}\cong^{op}\mathsf{InqAlg_{FCGW}^\otimes}.$$	
\end{theorem}{}
\begin{proof}
	By the definition of the functors $F$ and $G$, and by Propositions \ref{representation3} and \ref{representation4}, we have that $F\circ G\cong id $ and $G\circ F \cong id$. Hence, $F$ and $G$ together describe a dual categorical equivalence.
\end{proof}

\subsection{Equivalence of Team and Algebraic Semantics}

We use the categorical equivalences of the previous section to obtain some results on the equivalence of team and algebraic semantics of inquisitive and dependence logics. The equivalence of the two semantics can be proved from the former duality results in the case of inquisitive algebras, while in the case of dependence logics we only show a limited version for well-connected algebras. However, a full semantic equivalence for dependence algebras can be proved relying on our former algebraic completeness result. Notice that, since in this section we will not consider maps between algebras, we shall talk about Kripke frames both in the context of inquisitive and dependence algebras.

We start by providing canonical core-valuations to the canonical inquisitive and dependence algebras described above. With a slight abuse of  notation we indicate by $\mathcal{A}_\mathfrak{F}$ both the inquisitive and the dependence algebra dual to the finite Kripke frame $\mathfrak{F}$. If $\mathfrak{M}=(\mathfrak{F},V)$ is a finite Kripke model, we then obtain an inquisitive (dependence) model corresponding to $\mathfrak{M}$ by defining the \textit{canonical core-valuation} $ \mu^V:  \: \at \rightarrow \mathfrak{F} $ as follows:
\begin{align*}
\mu^V: & \: p \mapsto \{V^{-1}(p)\}^\downarrow.
\end{align*}

\noindent In this way we supplement $\mathcal{A}_\mathfrak{F}$ with a core-valuation, and we obtain a model $\mathcal{M}^V_\mathfrak{F}:=(\mathcal{A}_\mathfrak{F}, \mu^V)$ for inquisitive (dependence) logic. We say that $\mathcal{M}^V_\mathfrak{F}$ is the \textit{dual inquisitive (dependence) algebraic model} to the Kripke model $\mathfrak{M}=(\mathfrak{F},V)$.  We also recall that if $\mathfrak{F}$ is a Kripke frame and $s$ is an $R$-upset, then we denote by $\mathfrak{F}_{s}$ the subframe  $\mathfrak{F}{\upharpoonright} s=(W\cap s, R{\upharpoonright} s)$. We  proceed by first proving the following technical lemma.

\begin{lemma}\label{lemma21} Let $\mathfrak{M}=(\mathfrak{F},V)$ be a Kripke model and $\mathfrak{F}=(W,R)$. For all $s,t\in \mathsf{Up}(\mathfrak{F})$ the following facts hold:
	\begin{itemize}
		\item[(i)] $s\in \llbracket\phi \rrbracket^{\mathcal{M}^{V}_{\mathfrak{F}_{s}}}\Longleftrightarrow   \mathcal{M}^{V}_{\mathfrak{F}_{s}}\vDash^c \phi $. 
		\item[(ii)] If $t\subseteq s$, then $\llbracket\phi \rrbracket^{\mathcal{M}^{V}_{\mathfrak{F}_{t}} }=\llbracket\phi \rrbracket^{\mathcal{M}^{V}_{\mathfrak{F}_{s}} }\cap \{t \}^\downarrow$.
		\item[(iii)]  If  $s,t_i\in \mathsf{Up}(\mathfrak{F})$ for all $i\in I$, then $\{s \}^\downarrow =\bigcup \{ \{t_i \}^\downarrow : i\in I \}   \Longleftrightarrow t_i =s \text{ for some } i\in I.$
	\end{itemize}
\end{lemma}
\begin{proof}
	$\:$
	\begin{itemize}	
		\item[(i)] By construction, $\llbracket\phi \rrbracket^{\mathcal{M}^{V}_{\mathfrak{F}_{s}}}\subseteq \{s\}^\downarrow$ and, since $ \llbracket\phi \rrbracket^{\mathcal{M}^{V}_{\mathfrak{F}_{s}}} $ is downward closed, $s\in \llbracket\phi \rrbracket^{\mathcal{M}^{V}_{\mathfrak{F}_{s}}}$ entails $ \{s\}^\downarrow \subseteq \llbracket\phi \rrbracket^{\mathcal{M}^{V}_{\mathfrak{F}_{s}}}$. Then $s\in \llbracket\phi \rrbracket^{\mathcal{M}^{V}_{\mathfrak{F}_{s}}}$ if and only if  $\{s\}^\downarrow = \llbracket\phi \rrbracket^{\mathcal{M}^{V}_{\mathfrak{F}_{s}}}$ if and only if $\mathcal{M}^{V}_{\mathfrak{F}_{s}}\vDash^c \phi$. 
		
		\item[(ii)]  By induction on the complexity of $\phi$.
		\begin{itemize}
			\item If $\phi= p$, then $\llbracket p \rrbracket^{\mathcal{M}^{V}_{\mathfrak{F}_{t}}}= \{V^{-1}(p) \cap t\}^\downarrow   = \{V^{-1}(p) \}^\downarrow \cap \{t\}^\downarrow  = \llbracket p \rrbracket^{\mathcal{M}^{V}_{\mathfrak{F}_{s}} }\cap \{t\}^\downarrow.$
			\item  If $\phi= \bot$, then $\llbracket\phi \rrbracket^{\mathcal{M}^{V}_{\mathfrak{F}_{t}} }=\{\varnothing \}=\llbracket\phi \rrbracket^{\mathcal{M}^{V}_{\mathfrak{F}_{s}} }$.

			\item If $\phi=\psi \lor \chi$ or $\phi= \psi \land \chi$,  then our claim follows directly from the induction hypothesis.
			
			\item If $\phi= \psi \rightarrow \chi$, then: 
			\begin{align*}
			\llbracket \psi \rightarrow \chi \rrbracket^{\mathcal{M}^{V}_{\mathfrak{F}_{t}} }&= \bigcup\{  x\in \mathcal{M}^{V}_{\mathfrak{F}_{t}}  : x \cap  \llbracket \psi \rrbracket^{\mathcal{M}^{V}_{\mathfrak{F}_{t}} } \subseteq  \llbracket \chi \rrbracket^{\mathcal{M}^{V}_{\mathfrak{F}_{t}} }\}\\
			&= \bigcup\{  x\in \mathcal{M}^{V}_{\mathfrak{F}_{t}}  : x \cap  (\llbracket \psi \rrbracket^{\mathcal{M}^{V}_{\mathfrak{F}_{s}} }\cap \{t\}^\downarrow) \subseteq ( \llbracket \chi \rrbracket^{\mathcal{M}^{V}_{\mathfrak{F}_{s}} }\cap \{t\}^\downarrow)\} & \text{ (by i.h.) }\\
			& = \bigcup\{  x\in \mathcal{M}^{V}_{\mathfrak{F}_{s}}  : x  \cap  \llbracket \psi \rrbracket^{\mathcal{M}^{V}_{\mathfrak{F}_{s}} } \subseteq  \llbracket \chi \rrbracket^{\mathcal{M}^{V}_{\mathfrak{F}_{s}} }\} \cap \{t\}^\downarrow\\
			&=  (\llbracket \psi \rrbracket^{\mathcal{M}^{V}_{\mathfrak{F}_{s}} }\rightarrow  \llbracket \chi \rrbracket^{\mathcal{M}^{V}_{\mathfrak{F}_{s}} }) \cap  \{t\}^\downarrow\\
			& = \llbracket \psi \rightarrow \chi \rrbracket^{\mathcal{M}^{V}_{\mathfrak{F}_{s}} } \cap \{t\}^\downarrow.
			\end{align*}
			
			\item If $\phi= \psi \otimes \chi$, then we first notice that, since $\mathsf{Dw}_p(\mathsf{Up})(\mathfrak{F})$ is a Heyting algebra with the tensor as join, we have for all $t_0,t_1,t_2\in\mathsf{Up}(\mathfrak{F})$:
			\begin{align*}
			(*) \; &\; ( \{ t_0 \}^\downarrow \cap \{t_2 \}^\downarrow  )\otimes ( \{ t_1 \}^\downarrow \cap \{t_2 \}^\downarrow  ) = ( \{ t_0 \}^\downarrow \otimes \{t_1 \}^\downarrow  )\cap  \{t_2 \}^\downarrow.
			\end{align*}
			
			Now, we let $\llbracket \psi  \rrbracket^{\mathcal{M}^{V}_{\mathfrak{F}_{t}} } = \bigcup_{i\leq n} \{ k_i \}$ and $\llbracket \chi  \rrbracket^{\mathcal{M}^{V}_{\mathfrak{F}_{t}} }=\bigcup_{j\leq m} \{ z_j  \} $. We obtain:
			\begin{align*}
			\llbracket \psi \otimes \chi \rrbracket^{\mathcal{M}^{V}_{\mathfrak{F}_{t}} }&=  \llbracket \psi  \rrbracket^{\mathcal{M}^{V}_{\mathfrak{F}_{t}} } \otimes \llbracket  \chi \rrbracket^{\mathcal{M}^{V}_{\mathfrak{F}_{t}} } \\
			& = ( \llbracket \psi \rrbracket^{\mathcal{M}^{V}_{\mathfrak{F}_{s}} }\cap \{t\}^\downarrow) \otimes ( \llbracket \chi \rrbracket^{\mathcal{M}^{V}_{\mathfrak{F}_{s}} }\cap \{t\}^\downarrow) & \text{ (by i.h.) }\\
			& =  \bigcup_{i\leq n} (\{ k_i \}^\downarrow\cap \{t\}^\downarrow) \otimes  \bigcup_{j\leq m} (\{ z_j  \}^\downarrow\cap \{t\}^\downarrow ) & \text{}\\
			& =   \bigcup_{i\leq n, j\leq m} [(\{ k_i \}^\downarrow\cap \{t\}^\downarrow) \otimes  (\{ z_j  \}^\downarrow\cap \{t\}^\downarrow )] & \text{(by \textit{Dist})}\\
			& =   \bigcup_{i\leq n, j\leq m} [(\{ k_i \}^\downarrow \otimes  \{ z_j  \}^\downarrow)\cap \{t\}^\downarrow )] & \text{(by (*))}\\
			& =  \Big ( \bigcup_{i\leq n}\{ k_i \}^\downarrow \otimes  \bigcup_{j\leq m}\{ z_j \}^\downarrow \Big )\cap  \{t\}^\downarrow  & \text{(by \textit{Dist})}\\
			&= ( \llbracket \psi \rrbracket^{\mathcal{M}^{V}_{\mathfrak{F}_{s}} }\otimes  \llbracket \chi \rrbracket^{\mathcal{M}^{V}_{\mathfrak{F}_{s}} }) \cap  \{t\}^\downarrow & \text{} \\
			&= \llbracket \psi \otimes \chi \rrbracket^{\mathcal{M}^{V}_{\mathfrak{F}_{s}} }\cap \{t\}^\downarrow .
			\end{align*}
		\end{itemize}

		%\item If $\phi=\psi \lor \chi$, then:
		%\begin{align*}
		%\llbracket \psi \lor \chi \rrbracket^{\mathcal{M}^{V}_{\mathfrak{F}_{t}} }&= \llbracket \psi \rrbracket^{\mathcal{M}^{V}_{\mathfrak{F}_{t}} } \cup \llbracket \chi \rrbracket^{\mathcal{M}^{V}_{\mathfrak{F}_{t}} } \\
		%& = ( \llbracket \psi \rrbracket^{\mathcal{M}^{V}_{\mathfrak{F}_{s}} }\cap \{t\}^\downarrow) \cup ( \llbracket \chi \rrbracket^{\mathcal{M}^{V}_{\mathfrak{F}_{s}} }\cap \{t\}^\downarrow) & \text{ (by i.h.) }\\
		%&= ( \llbracket \psi \rrbracket^{\mathcal{M}^{V}_{\mathfrak{F}_{s}} }\cup  \llbracket \chi \rrbracket^{\mathcal{M}^{V}_{\mathfrak{F}_{s}} }) \cap  \{t\}^\downarrow\\
		%&= \llbracket \psi \lor \chi \rrbracket^{\mathcal{M}^{V}_{\mathfrak{F}_{s}} }\cap \{t\}^\downarrow .
		%\end{align*}
		%\item If $\phi= \psi \land \chi$, we proceed as in the previous case.
		
		\item[(iii)]  $(\Rightarrow)$ Suppose $t_i \neq s \text{ for all } i\in I $. If there is some $i\in I$ and some element $y_i\in t_i \setminus s$, then  $Ry_i\notin \{s\}^\downarrow$ but $Ry_i\in \bigcup \{ \{t_i\}^\downarrow : i\in I \}$, proving our claim. Otherwise, for each $t_i$ there is an element $x_i$ such that $x_i\in s \setminus t_i$, hence $Rx_i\in \{s\}^\downarrow \setminus \{t_i\}^\downarrow$. It follows that $\{Rx_i\}_{i\in I}\subseteq \{s\}^\downarrow$ but $\{Rx_i\}_{i\in I}\nsubseteq \bigcup \{ \{t_i\}^\downarrow : i\in I \}$, showing $\{s\}^\downarrow \neq \bigcup \{ \{t_i\}^\downarrow : i\in I \} $.	 $(\Leftarrow)$ Obvious. \qedhere
	\end{itemize}
\end{proof}

\noindent By using the previous lemma we obtain the following proposition, which we shall use later in the proof of Theorem \ref{equiv}.

\begin{proposition}\label{propo1}
	Let $\mathfrak{M}=(\mathfrak{F},V)$ be a Kripke model and $s,t\in \mathsf{Up}(\mathfrak{F})$ such that $t\subseteq s$. Then  $\mathcal{M}^{V}_{\mathfrak{F}_{t}} \vDash^c \phi $ if and only if $t\in \llbracket\phi \rrbracket^{\mathcal{M}^{V}_{\mathfrak{F}_{s}}} $.
\end{proposition}
\begin{proof}
	For any $s,t\in \mathsf{Up}(\mathfrak{F})$ such that $t\subseteq s$, we have:
	\begin{align*}
	t\in \llbracket\phi \rrbracket^{\mathcal{M}^{V}_{\mathfrak{F}_{s}}} & \Longleftrightarrow t\in \llbracket\phi \rrbracket^{\mathcal{M}^{V}_{\mathfrak{F}_{s}}} \cap \{ t \}^\downarrow \\
	&  \Longleftrightarrow t\in \llbracket\phi \rrbracket^{\mathcal{M}^{V}_{\mathfrak{F}_{t}}} &\text{(by Lemma \ref{lemma21}(ii))} \\
	&\Longleftrightarrow  \mathcal{M}^{V}_{\mathfrak{F}_{t}} \vDash^c \phi, &\text{(by Lemma \ref{lemma21}(i))}
	\end{align*}
	\noindent which proves our claim.
\end{proof}

The next theorem finally shows that any inquisitive or dependence formula is true in a finite Kripke model if and only if it is valid, under the canonical core-valuation, in its dual inquisitive algebraic model.

\begin{theorem}\label{equiv}
	Let $\mathfrak{M}=(\mathfrak{F},V)$ be a finite Kripke frame and $\mathcal{M}^V_\mathfrak{F}= (\mathcal{A}_\mathfrak{F}, \mu^V)$ its dual inquisitive (dependence) algebraic model. Then  $\mathfrak{M}\vDash \phi$ if and only if $\mathcal{M}^V_\mathfrak{F}\vDash^c \phi$.
\end{theorem}
\begin{proof}
	By induction on the complexity of $\phi$. 		
	\begin{itemize}
		\item For $p\in \at$ we have that:
		\begin{align*}
		\mathfrak{M} \vDash  p & \Longleftrightarrow \forall w \in W \ ( w(p)=1) \\
		& \Longleftrightarrow W=  \{ w\in W : w(p)=1   \} \\
		& \Longleftrightarrow \{W\}^\downarrow =  \{ w\in W : w(p)=1   \}^\downarrow \\
		& \Longleftrightarrow 1_{\mathcal{M}^V_\mathfrak{F}} = \mu^V(p) \\
		&\Longleftrightarrow\mathcal{M}^V_\mathfrak{F}\vDash^c p.
		\end{align*}
		
		\item 	For $\phi= \bot $ we have:
		\begin{align*}
		\mathfrak{M} \vDash  \bot \Longleftrightarrow W=\varnothing 
		\Longleftrightarrow \{ W \}^\downarrow=  \{ \varnothing   \}^\downarrow 			\Longleftrightarrow 1_{\mathcal{M}^V_\mathfrak{F}} = 0_{\mathcal{M}^V_\mathfrak{F}} \Longleftrightarrow\mathcal{M}^V_\mathfrak{F}\vDash^c \bot.
		\end{align*}

		\item 	For $\phi= \psi \lor \chi $ we have:
		\begin{align*}
		\mathfrak{M}\vDash \psi \lor \chi  &\Longleftrightarrow \mathfrak{M}\vDash \psi \text{ or } \mathfrak{M}\vDash \chi& \\
		&\Longleftrightarrow 1_{\mathcal{M}^V_\mathfrak{F}}=  \llbracket \psi \rrbracket^{\mathcal{M}^V_\mathfrak{F}} \text{ or }  1_{\mathcal{M}^V_\mathfrak{F}}=\llbracket \chi \rrbracket^{\mathcal{M}^V_\mathfrak{F}}& \text{ (by induction hypothesis) }  \\
		&\Longleftrightarrow 1_{\mathcal{M}^V_\mathfrak{F}}=  \llbracket \psi \rrbracket^{\mathcal{M}^V_\mathfrak{F}} \lor \llbracket \chi \rrbracket^{\mathcal{M}^V_\mathfrak{F}} & (\text{by well-connectedness of }  \mathcal{A}_{\mathfrak{F}}) \\
		& \Longleftrightarrow 1_{\mathcal{M}^V_\mathfrak{F}}=  \llbracket \psi \lor \chi \rrbracket^{\mathcal{M}^V_\mathfrak{F}} &\\
		& \Longleftrightarrow\mathcal{M}^V_\mathfrak{F}\vDash^c \psi \lor \chi.
		\end{align*}

		\item 	For $\phi= \psi \land \chi $ the claim follows by straightforward application of the induction hypothesis.
		
		%		\begin{align*}
		%		\mathfrak{M}\vDash \psi \land \chi  &\Longleftrightarrow \mathfrak{M}\vDash \psi \text{ and } \mathfrak{M}\vDash \chi \\
		%		&\Longleftrightarrow 1_{\mathcal{M}^V_\mathfrak{F}}=  \llbracket \psi \rrbracket^{\mathcal{M}^V_\mathfrak{F}} \text{ and }  1_{\mathcal{M}^V_\mathfrak{F}}=\llbracket \chi \rrbracket^{\mathcal{M}^V_\mathfrak{F}}   & \text{ (by induction hypothesis) } \\
		%		&\Longleftrightarrow 1_{\mathcal{M}^V_\mathfrak{F}}=  \llbracket \psi \rrbracket^{\mathcal{M}^V_\mathfrak{F}} \wedge \llbracket \chi \rrbracket^{\mathcal{M}^V_\mathfrak{F}} \\
		%		& \Longleftrightarrow 1_{\mathcal{M}^V_\mathfrak{F}}=  \llbracket \psi \land \chi \rrbracket^{\mathcal{M}^V_\mathfrak{F}} \\
		%		& \Longleftrightarrow\mathcal{M}^V_\mathfrak{F}\vDash^c \psi \land \chi.
		%		\end{align*}
		
		\item 	For $\phi= \psi \rightarrow \chi $ we have: 
		\begin{align*}
		\mathfrak{M}\vDash \psi \rightarrow \chi  & \Longleftrightarrow \forall t \ ( \text{if }t\subseteq W \text{ and } \mathfrak{M},t\vDash\psi \text{ then } \mathfrak{M},t\vDash \chi ) \\
		& \xLongLeftRightArrow[]{\text{Prop. \ref{upsetprop}}} \forall t \ ( \text{if }t\in \mathsf{Up}(\mathfrak{F}) \text{ and } \mathfrak{M},t\vDash\psi \text{ then } \mathfrak{M},t\vDash \chi ) \\
		& \xLongLeftRightArrow[]{ \text{ by i.h. }} \forall t \ ( \text{if }t\in \mathsf{Up}(\mathfrak{F}) \text{ and } \mathcal{M}^{V}_{\mathfrak{F}_{t}}\vDash^c \psi \text{ then } \mathcal{M}^{V}_{\mathfrak{F}_{t}}\vDash^c \chi  ) \\
		&\Longleftrightarrow  \{W\}^\downarrow =  \{ t\in \mathsf{Up}(\mathfrak{F}) : \mathcal{M}^{V}_{\mathfrak{F}_{t}}\vDash^c \psi \Rightarrow \mathcal{M}^{V}_{\mathfrak{F}_{t}}\vDash^c \chi    \}     \\
		&\Longleftrightarrow  \{W\}^\downarrow = \bigcup \{ x \in \mathcal{A}_{\mathfrak{F}} : t\in x \ \Rightarrow \ [ \mathcal{M}^{V}_{\mathfrak{F}_{t}}\vDash^c \psi \Rightarrow \mathcal{M}^{V}_{\mathfrak{F}_{t}}\vDash^c \chi ]   \} \\
		&  \xLongLeftRightArrow[]{\text{Prop. \ref{propo1}}}{}  \{W\}^\downarrow = \bigcup \{  x \in \mathcal{A}_{\mathfrak{F}} :  t\in x \ \Rightarrow \ [ t\in  \llbracket \psi \rrbracket^{\mathcal{M}_{\mathfrak{F}}} \Rightarrow t\in  \llbracket \chi \rrbracket^{\mathcal{M}^V_\mathfrak{F}} ]   \} \\
		&\Longleftrightarrow  \{W\}^\downarrow = \bigcup \{ x \in \mathcal{A}_{\mathfrak{F}} : x\cap \llbracket \psi \rrbracket^{\mathcal{M}^V_\mathfrak{F}} \subseteq \llbracket \chi \rrbracket^{\mathcal{M}^V_\mathfrak{F}}   \}    \\
		&\Longleftrightarrow 1_{\mathcal{M}^V_\mathfrak{F}}=  \llbracket \psi \rightarrow \chi \rrbracket^{\mathcal{M}^V_\mathfrak{F}} .
		\end{align*}	
	\end{itemize}
	
	If $\mathcal{M}^V_\mathfrak{F}$ is an algebraic dependence model, we need to check also the case for the tensor disjunction.
	
	$(\Rightarrow)$ Suppose $\mathfrak{M}\vDash \psi \otimes \chi$. then by Proposition \ref{upsetprop} there are two upsets $t,r\in \mathsf{Up}(\mathfrak{F})$ such that $ t\cup r = W$ and $\mathfrak{M}_t\vDash \psi,  \mathfrak{M}_r\vDash \chi$. By induction hypothesis, we obtain that $\mathcal{M}^{V}_{\mathfrak{F}_{t}}\vDash^c \psi$   and  $\mathcal{M}^{V}_{\mathfrak{F}_{r}}\vDash^c \chi$, thus by Proposition \ref{propo1} $t\in \llbracket  \psi \rrbracket^{\mathcal{M}^V_\mathfrak{F}}$ and $r\in \llbracket  \chi \rrbracket^{\mathcal{M}^V_\mathfrak{F}}$. Now, since $W=t\cup r$, we have:
	\begin{align*}
	\llbracket \psi\otimes \chi \rrbracket^{\mathcal{M}^V_\mathfrak{F}} &= \bigcup \{ \{u\cup v\}^\downarrow : u\in \llbracket  \psi \rrbracket^{\mathcal{M}^V_\mathfrak{F}} \text{ and } v\in \llbracket  \psi \rrbracket^{\mathcal{M}^V_\mathfrak{F}} \} \\
	& = \bigcup \{ \{t\cup r\}^\downarrow \} \\
	& = \{W \}^\downarrow.
	\end{align*}
	Hence $\mathcal{M}^V_\mathfrak{F} \vDash^c \psi\otimes \chi$.
	
	($\Leftarrow$) Now suppose $\mathcal{M}^V_\mathfrak{F} \vDash^c \psi\otimes \chi$, then we have that:
	$$\{W\}^\downarrow=\bigcup \{ \{u\cup v\}^\downarrow : u\in \llbracket  \psi \rrbracket^{\mathcal{M}^V_\mathfrak{F}} \text{ and } v\in \llbracket  \psi \rrbracket^{\mathcal{M}^V_\mathfrak{F}} \}  .$$
	\noindent By Lemma \ref{lemma21}(iii) there are $r,t\in \mathsf{Up}(\mathfrak{F})$ such that $\{W\}^\downarrow=\{t\cup r\}^\downarrow$, $t\in \llbracket  \psi \rrbracket^{\mathcal{M}^V_\mathfrak{F}}$ and $r\in \llbracket  \chi \rrbracket^{\mathcal{M}^V_\mathfrak{F}} $. Then, it follows that $t\in \llbracket  \psi \rrbracket^{\mathcal{M}^V_\mathfrak{F}}\cap \{ t\}^\downarrow$ and $r\in \llbracket  \chi \rrbracket^{\mathcal{M}^V_\mathfrak{F}} \cap \{r\}^\downarrow$, hence  by Lemma \ref{lemma21}(ii) we have  $t\in \llbracket  \psi \rrbracket^{\mathcal{M}_{\mathfrak{F}_t}}$ and $r\in \llbracket  \chi \rrbracket^{\mathcal{M}^{V}_{\mathfrak{F}_{r}}}$, thus by Lemma \ref{lemma21}(i)  $\mathcal{M}^{V}_{\mathfrak{F}_{t}}\vDash^c \psi$ and $\mathcal{M}_{\mathfrak{F}_r}\vDash^c \chi$.   By induction hypothesis, we have  $\mathfrak{M}_{t} \vDash \psi $ and $\mathfrak{M}_r\vDash \chi$, which together with $r\cup t = W$  yields $\mathfrak{M}\vDash \psi\otimes \chi$.
\end{proof}

The previous theorem establishes that a formula is true in a Kripke model if and only if it is true in its dual  algebraic model. Now we proceed in the converse direction and we prove that a formula is true in an algebraic model if it is true in its corresponding Kripke model. As we did above, we proceed by defining canonical valuations over Kripke frames dual to finite, core-generated, well-connected inquisitive (or dependence) algebras.

Given an algebraic model $\mathcal{M}=(\mathcal{A},\mu)$ such that  $\mathcal{A}$ is a finite, core-generated, well-connected inquisitive (dependence) algebra, we can find by Proposition \ref{representation1} (and Proposition \ref{representation4}),  a finite Kripke frame $\mathfrak{F}$ such that $\mathcal{A}\cong\mathcal{A}_{\mathfrak{F}}$. Let  $h:\mathcal{A}\to \mathcal{A}_{\mathfrak{F}}$ be a isomorphism and let $\mu'=h\circ\mu$.  We define the canonical valuation $V: \: W \rightarrow \wp(\at)$ over the frame $\mathfrak{F}$  as follows:
\begin{align*}
V: \: &\: w \mapsto  \{ p\in\at :  \{R[w]\}^\downarrow \subseteq \mu'(p)   \}.
\end{align*}

\noindent It is straightforward to verify that $V$ is a suitable valuation, i.e. that if $wRv$ and $p\in V(w)$, then $p\in V(v)$. We say that the Kripke model $\mathfrak{M}= (\mathfrak{F},  V)$ is the dual of $\mathcal{M}=(\mathcal{A},\mu)$. We show that $\mathfrak{M}$ validates exactly the same formulas of the original algebraic model.

\begin{proposition}\label{equiv2}
	Let $\mathcal{M}$ be a finite, core-generated, well-connected inquisitive (dependence) algebraic model and let $\mathfrak{M}$ be its dual Kripke model, then  $\mathcal{M}\vDash^c \phi$ if and only if $\mathfrak{M}\vDash \phi$.
\end{proposition}
\begin{proof}
	We prove this theorem for inquisitive algebraic models only, as the proof for dependence models is the same, using the corresponding duality result. 
	
	Let $\mathcal{M}=(\mathcal{A},\mu)$ be an inquisitive algebraic model with  $\mathcal{A}\in \mathsf{InqAlg_{FCGW}}$ and let $\mathfrak{M}=(\mathfrak{F}, V) $ be its dual Kripke model. By Theorem \ref{duality} we have that $\mathcal{A}\cong\mathcal{A}_{\mathfrak{F}}$.  Let $h:\mathcal{A}\to \mathcal{A}_{\mathfrak{F}}$ witness this isomorphism and consider the canonical algebraic downward team model $\mathcal{M}^V_{\mathfrak{F}}=(\mathcal{A}_\mathfrak{F}, \nu^V )$.  Then, since $\mu'(p)=h\circ\mu(p)=\{t\}^\downarrow$ for some $t\in\mathsf{Up}(\mathfrak{F})$, we have:
	\begin{align*}
	\nu^V(p) & =   \{V^{-1}(p)\}^\downarrow =     \{ w\in \mathfrak{F} : \{R[w]\}^\downarrow \subseteq \mu'(p)    \}^\downarrow  = \{ w\in \mathfrak{F} : w\in t    \}^\downarrow  =   h\circ\mu (p).
	\end{align*}
	
	\noindent Then, since $h$ is a isomorphism and for all $p\in \at$ \, $ \nu^V(p)=h\circ\mu (p)$, it follows that for all $\phi\in \langInt$, $h\llbracket \phi(\vec{x})\rrbracket^{(\mathcal{A},\mu)}= \llbracket \phi(\vec{x})\rrbracket^{(\mathcal{A}_\mathfrak{F},\nu)}$. Therefore, we have that for all $\phi\in \langInt$, $\mathcal{M}\vDash^c \phi$ if and only if $\mathcal{M}^V_{\mathfrak{F}}\vDash^c \phi$. By Theorem \ref{equiv} we have $\mathcal{M}^V_{\mathfrak{F}}\vDash^c \phi$ if and only if $\mathfrak{M}\vDash \phi$. Finally, this entails  $\mathcal{M}\vDash^c \phi$ if and only if $\mathfrak{M}\vDash \phi$, which proves our claim.
\end{proof}

The previous results show that the dual equivalences that we studied in the previous section also preserve the validity of inquisitive and dependence formulas. Now, for any intermediate inquisitive logic $\inqint$ we let $ \mathsf{KF\Lambda} $ be the class of Kripke frames $\mathfrak{F}$ such that $\mathfrak{F}\vDash \Lambda $, and we let $ \mathsf{KF\Lambda_F} $ be its subcollection of finite frames. We denote by $ \mathsf{InqAlg\Lambda_{FCGW}} $ the subcategory of $ \mathsf{InqAlg\Lambda} $ consisting finite, core-generated, well-connected intuitionistic inquisitive algebras.  We define analogously $ \mathsf{KF \Lambda^\otimes} $, $ \mathsf{KF\Lambda_F^\otimes} $ and $ \mathsf{InqAlg\Lambda^\otimes_{FCGW}} $ for every intermediate dependence logic $\inqint^\otimes$. We  obtain the following corollary of Theorem \ref{equiv} and Proposition \ref{equiv2}.

\begin{corollary}
	For every intermediate inquisitive logic $\inqint$, and every intermediate dependence logic $\inqint^\otimes$, we have the following:
	\begin{align*}
	\mathfrak{F}\in  \mathsf{KF\Lambda_F} & \Longleftrightarrow \mathcal{A}_{\mathfrak{F}}\in    \mathsf{InqAlg\Lambda_{FCGW}}; \\
	\mathfrak{F}\in  \mathsf{KF \Lambda_F^\otimes} & \Longleftrightarrow \mathcal{A}_{\mathfrak{F}}  \in  \mathsf{InqAlg\Lambda^\otimes_{FCGW}}.	\end{align*}
\end{corollary}

Finally, we use the results of this section to prove the following  theorems, which show the equivalence between team and algebraic semantics.

\begin{theorem}[Semantic Equivalence I]\label{semanticequivalence}
	$\:$
	\begin{itemize}
		\item[(i)] The class of $\inqI$-algebras is semantically equivalent to the class of all Kripke frames, i.e. for all $\phi\in\langInt$:
		\begin{align*}
		\mathsf{KF}\vDash \phi  &\Longleftrightarrow  \mathsf{InqAlg}\vDash^c \phi.
		\end{align*}
		
		\item [(ii)] The class of well-connected $\inqI^\otimes$-algebras is semantically equivalent to the class of all Kripke frames, i.e. for all $\phi\in\langInqI$:
		\begin{align*}
		\mathsf{KF}\vDash \phi  &\Longleftrightarrow  \mathsf{InqAlg_W^\otimes}\vDash^c \phi.
		\end{align*}
	\end{itemize}
\end{theorem}
\begin{proof}
	(i)	($\Rightarrow$) Suppose $\mathsf{InqAlg}\nvDash^c \phi$, then by Theorem \ref{birkhoff2} there is a finite, core-generated, well-connected  $\inqI$-algebra such that  $\mathcal{A}\nvDash^c \phi$. Hence, for some core-valuation $\mu$, we have that $(\mathcal{A},\mu)\nvDash^c\phi$ and then, by   Proposition \ref{equiv2}, it follows that for some finite Kripke model $\mathfrak{M}$ we have  that $\mathfrak{M} \nvDash \phi$ and thus  $ \mathsf{KF}\nvDash \phi $ .  ($\Leftarrow$) Suppose $\mathsf{KF}\nvDash \phi $ then, by Theorem \ref{FMP}, there is a finite Kripke frame $\mathfrak{F}=(\mathfrak{F},V)$ such that $\mathfrak{M}\nvDash \phi$. By Theorem \ref{equiv} we have $\mathcal{M}^V_\mathfrak{F}\nvDash^c \phi$, where $\mathcal{M}^V_\mathfrak{F}=(\mathcal{A}_\mathfrak{F},\mu^V)$ and $\mathcal{A}_\mathfrak{F}\in \mathsf{InqAlg_{FCGW}}$. Finally, this shows that  $\mathsf{InqAlg}\nvDash^c \phi$.
	
	(ii) Analogous to (i), by using Theorem \ref{birkhoff4}.
\end{proof}

We have already remarked in Section \ref{two} that a Kripke frame $\mathfrak{F}$ such that $\mathfrak{F}\vDash \inqB  $ or $\mathfrak{F}\vDash \inqB^\otimes  $ is classical and it can be viewed as a set of assignments. We can then prove the following result for the logics $\inqB$ and $\inqB^\otimes$. We denote by \textsf{Team} the class of all classical teams, by $ \mathsf{InqBAlg} $ the class of all $\inqB$-algebras and by $ \mathsf{InqBAlg^\otimes} $ the class of all $\inqB^\otimes$-algebras.

\begin{theorem}[Semantic Equivalence II]\label{semanticequivalence2}
	$\:$
	\begin{itemize}
		\item[(i)] The class of $\inqB$-algebras is semantically equivalent to the class of all teams, i.e. for all $\phi\in\langInt$:
		\begin{align*}
		\mathsf{Team}\vDash \phi  &\Longleftrightarrow  \mathsf{InqAlg}\vDash^c \phi.
		\end{align*}
		
		\item [(ii)] The class of $\inqI^\otimes$-algebras is semantically equivalent to the class of all teams, i.e. for all $\phi\in\langInqI$:
		\begin{align*}
		\mathsf{Team}\vDash \phi  &\Longleftrightarrow  \mathsf{InqAlg}^\otimes\vDash^c \phi.
		\end{align*}
	\end{itemize}
\end{theorem}
\begin{proof}
	(i) Analogous to the Proof of Theorem \ref{semanticequivalence}(i) using the finite model property of $\inqB$.
	(ii) Since $\inqB^\otimes{\upharpoonright}\langIntsta =\cpc $ is locally tabular we apply Theorem \ref{birkhoff6} and proceed otherwise as in the Proof of Theorem \ref{semanticequivalence}(i). 
\end{proof}

Notice that the previous theorems were proven using only the results of this section, together with the finite model property and the completeness of team semantics. In particular, they do not rely on the algebraic completeness theorem of Section \ref{two}. This is interesting, as we can provide an alternative proof of the algebraic completeness for some inquisitive and dependence logics. %By combining together the two previous Semantic Equivalence Theorems, the Completeness Theorems \ref{complteam} and \ref{classcomplteam}, we then obtain the following algebraic completeness results directly as corollaries of the results in this section.

\begin{corollary}[Algebraic Completeness]  \label{inqicomplete}
	$\:$
	\begin{itemize}
		\item[(i)] For all $\phi\in\langInt$, $\phi \in \inqI \Longleftrightarrow \mathsf{InqAlg}\vDash^c \phi$.
		
		\item[(ii)] For all $\phi\in\langInt$, $\phi \in \inqB \Longleftrightarrow \mathsf{InqBAlg}\vDash^c \phi $.
		
		\item[(iii)] For all $\phi\in\langInqI$, $\phi \in \inqB^\otimes \Longleftrightarrow \mathsf{InqBAlg}^\otimes\vDash^c \phi$.
	\end{itemize}
\end{corollary}
\begin{proof}
	(i) By Theorem \ref{complteam} we have that $\phi \in \inqI$ if and only if $\mathsf{KF}\vDash \phi $ and, by Proposition \ref{FMP}, this is the case if and only if $\mathsf{KF_F}\vDash \phi $. Finally, this is equivalent, by Theorem \ref{semanticequivalence}, to $\mathsf{InqAlg}\vDash^c \phi$. (ii) and (iii) are proven analogously using Theorems \ref{classcomplteam} and \ref{semanticequivalence2}.
\end{proof}

\noindent We then obtained an alternative  proof of the algebraic completeness for the logics $\inqI,\inqB$ and $\inqB^\otimes$. Whether this method could be used to prove the algebraic completeness of other intermediate inquisitive and dependence logics should be object of further investigations.

Finally, we conclude this section by remarking that, using the algebraic completeness theorem of Section \ref{two}, the following result follows:

\begin{theorem}[Semantic Equivalence III]\label{semanticequivalence3}
	$\:$
	The class of $\inqint^\otimes$-algebras is semantically equivalent to the class of all Kripke frames satisfying $\Lambda$, i.e. for all $\phi\in\langInqI$:
	\begin{align*}
	\mathsf{KF}\vDash \phi  &\Longleftrightarrow  \mathsf{InqAlg}^\otimes\vDash^c \phi.
	\end{align*}
\end{theorem}
\begin{proof}
	By Theorem \ref{complteam}, $ \mathsf{KF}\vDash \phi$ holds if and only if $\phi \in \inqI^\otimes $. By Theorem \ref{algcomple} the former is equivalent to $\mathsf{InqAlg}^\otimes\vDash^c \phi$.
\end{proof}

\section{Concluding Remarks and Open Problems} \label{ten}

In this article we studied intermediate inquisitive and dependence logics from an algebraic perspective. We presented in Section \ref{one} an axiomatisation of inquisitive and dependence logics, and we introduced in Section \ref{two} algebraic semantics  using so-called inquisitive and dependence algebras. We then adapted the standard method of free algebras to our setting, and we used it to prove that every intermediate inquisitive and dependence logic is algebraically complete.

We then considered, in Section \ref{three}, several model-theoretic properties of the classes of inquisitive and dependence algebras. We defined core-generated and well-connected inquisitive and dependence algebras and we proved several properties concerning them. Most importantly, we then focused on finite, core-generated, well-connected inquisitive and dependence algebras. In the inquisitive case we proved a version of Birkhoff's theorem, showing that $\mathsf{InqAlg\Lambda}\vDash^c \phi$ if and only if $\mathsf{InqAlg\Lambda_{FCGW}}\vDash^c \phi$. Differently, in the dependence case, we proved two weaker versions of this result: one stating that $\mathsf{InqAlg\Lambda_{W}}^\otimes\vDash^c \phi$ if and only if $\mathsf{InqAlg\Lambda_{FCGW}}^\otimes\vDash^c \phi$ and a second one stating that $\mathsf{InqAlg\Lambda}^\otimes\vDash^c \phi$ if and only if $\mathsf{InqAlg\Lambda_{FCGW}^\otimes}\vDash^c \phi$ holds whenever $\Lambda$ is locally tabular. 

Finally, in Section \ref{four}, we focused on the relation between frames and algebras, and the relation between algebraic and team semantics. To this end, we proved that the category $\mathsf{Pos_F}$ is dual to $ \mathsf{InqAlg_{FCGW}}$, and that $\mathsf{KF_F}$ is dual to $ \mathsf{InqAlg_{FCGW}}^\otimes $. We then derived several results concerning the equivalence of team and algebraic semantics and we provided an alternative proof of algebraic completeness for the logics $\inqI,\inqB,\inqB^\otimes$.

The main goal of this article was to provide a workable algebraic framework for inquisitive and dependence logics. The results we obtained show that standard algebraic methods can be adapted and used to study these logics, even despite the fact that they are non-standard systems where uniform substitution fails. The present work also suggests some possible directions for future investigations.

Firstly, as we have already remarked, the version of Birkhoff's theorem that we proved for dependence logics differs from the version we proved for inquisitive logics. Is it possible to prove a stronger result and show that, for any intermediate dependence logic $\inqint^\otimes$ and any formula $\phi\in \langInqI$, $\mathsf{InqAlg\Lambda}\vDash^c \phi$ if and only if $\mathsf{InqAlg\Lambda_{FCGW}}\vDash^c \phi$? %In fact, the fact that  both the class of dependence algebras and the class of well-connected dependence algebras are semantically equivalent to Kripke frames, by Theorems \ref{semanticequivalence} and \ref{semanticequivalence3}, seems to suggest this could be the case. 

It should also be considered for what intermediate inquisitive and dependence logics we can give a completeness proof using duality, as we did in Section \ref{four} for $\inqI, \inqB$ and $\inq^\otimes$. Interestingly, this problem relates to the question whether all intermediate inquisitive and dependence logics are complete with respect to some classes of Kripke frames, and whether they all have the finite model property. 

It is also natural to investigate whether the representation theorems for $ \mathsf{InqAlg_{FCGW}}$ and $ \mathsf{InqAlg_{FCGW}}^\otimes$  can be extended to the infinite case, namely to $ \mathsf{InqAlg_{CGW}}$ and $ \mathsf{InqAlg_{CGW}}^\otimes$. Since the representation of infinite algebraic structures also involves topological dualities, this issue also relates to the question whether it is possible to give a topological semantics to intermediate inquisitive and dependence logics. This problem  has been considered in \cite{grilletti} for $\inqB$, but it has not been investigated in the general case of all intermediate inquisitive and dependence logics.

Finally, we should consider whether it is possible to prove that the algebraic semantics outlined in this articles is any in any sense, unique -- as it happens for the standard algebraic semantics of $\cpc$ and $\ipc$ and for every standard algebraizable logics. This would require to develop a framework for algebraisability for logics without uniform substitution. We leave this and the previous problems to future investigations.

\printbibliography

%
%\section*{Acknowledgements}
%
%I am very grateful to Fan Yang for many helpful conversations and for her constant support and advice. I would also like to thank Georgi Nakov, Gianluca Grilletti, Tommaso Moraschini and Nick Bezhanishvili for helpful discussions. This research was supported by grant 336283 of the Academy of Finland and Research Funds of the University of Helsinki.

\end{document}